\documentclass[a4paper,11pt]{article}
\usepackage{amssymb}
\usepackage{epsfig}
\usepackage{graphics}
\usepackage{graphicx}
\usepackage{pgfplots}
\usetikzlibrary{patterns}
\usepackage{tikz-cd}
\usepackage{amsmath, amsthm}
\usepackage{authblk}
\usepackage{caption}
\usepackage{subcaption}
\usepackage{multirow}
\usepackage{longtable}
\usepackage{lscape}
\usepackage{booktabs}
\usepackage{lineno}
\newtheorem{thm}{Theorem}
\newtheorem{rem}{Remark}

\newtheorem{cor}{Corollary}
\newtheorem{prop}[thm]{Proposition}
\newtheorem{lem}[thm]{Lemma}

\usepackage[margin=1in]{geometry}
\usepackage{epstopdf}

\usepackage{scalerel,stackengine}
\stackMath
\newcommand\reallywidehat[1]{%
	\savestack{\tmpbox}{\stretchto{%
			\scaleto{%
				\scalerel*[\widthof{\ensuremath{#1}}]{\kern-.6pt\bigwedge\kern-.6pt}%
				{\rule[-\textheight/2]{1ex}{\textheight}}
			}{\textheight}%
		}{0.5ex}}%
	\stackon[1pt]{#1}{\tmpbox}%
}

\usepackage{pgfplots}
\allowdisplaybreaks
\usepackage{xcolor}
\usepackage{endnotes}
\let\footnote=\endnote

\def\bf{\mathbf}
\def\bb{\mathbb}
\def\cal{\mathcal}
\def\bsm{\boldsymbol}

\usepackage[numbers]{natbib}

\vspace{1cm}

\begin{document}
\title{Computing Two-Stage Robust Optimization with Mixed Integer Structures}
\author {Wei Wang and Bo Zeng\vspace{-10pt}}
\date{Industrial Engineering, University of Pittsburgh\\ w.wei@pitt.edu, \ \ bzeng@pitt.edu
}
\maketitle

\renewcommand{\baselinestretch}{1.5}

\abstract{Mixed integer sets have a strong modeling capacity to describe practical systems. Nevertheless, incorporating a mixed integer set often renders an optimization formulation drastically more challenging to compute. In this paper, we study how to effectively solve two-stage robust mixed integer programs built upon decision-dependent uncertainty (DDU). We particularly consider those with a mixed integer DDU set that has a novel modeling power on describing complex dependence between decisions and uncertainty, and those with a mixed integer recourse problem that captures discrete recourse decisions. 
	
	Note that there has been very limited research on computational algorithms for these challenging problems. To this end, several exact and approximation solution algorithms are presented, which address some critical limitations of existing counterparts (e.g., the relatively complete recourse assumption) and greatly expand our solution capacity. In addition to  theoretical analyses on those algorithms' convergence and computational complexities, a set of numerical studies are performed to showcase mixed integer structures' modeling capabilities and to demonstrate those algorithms' computational strength.  
}

\section{Introduction}
In the context of mathematical optimization, mixed integer structures are necessary to accurately describe physical systems characterized by discrete states, alternative choices, logic relationships or behavior discontinuities. Also they can be used to build  tight approximations for complicated functions that otherwise might be computationally infeasible.  Nevertheless, compared to its continuous counterpart, a formulation with a mixed integer structure typically incurs a substantially heavier computational burden \cite{NemWol.1988,wolsey2020integer}.  It is particularly the case when we are seeking an exact solution, which generally requires enumeration on all possible discrete choices.  Hence, rather than taking a brute-force fashion,  researchers have recognized that  enumeration should be performed with sophisticated strategies and advanced algorithm designs to mitigate the associated computation complexity.       

For the classical two-stage robust optimization (RO) defined on a fixed uncertainty set, referred to as a \textit{decision-independent uncertainty} (DIU) set, computing an instance with a mixed integer recourse is much more demanding than the one with a linear programming (LP) recourse. For example,  basic column-and-constraint generation (C\&CG) method is implemented in a nested fashion to handle mixed integer recourse exactly, which may drastically increase the computational time  \cite{zhao2012exact}. To address such a challenge, several approximation strategies  \citep{bertsimas2010power,hanasusanto2015k,zhang2018ambulance} have been studied to compute practical instances. They certainly have a great advantage on ensuring the computational efficiency. Yet, the quality of approximate solutions cannot be guaranteed. We mention that the case with a mixed integer structure appearing in the first stage or in the uncertainty set is comparatively much easier, as it typically does not require specialized algorithms other than those for regular mixed integer programs (MIPs).

Recently,  the traditional approach that employs a fixed set to model the concerned uncertainty has evolved to use a decision-dependent set to capture  non-fixed uncertainty that changes with respect to the choice of the decision maker (\cite{nohadani2018optimization,lappas2018robust,zhang2020unified,chen2023robust}). Indeed, \textit{decision-dependent uncertainty} (DDU) often occurs in our societal system,  e.g., consumers' demand or drivers' route selections are largely affected by  price or toll  set by the system operator~\cite{haghighat2023robust}. As another example, in system maintenance, components' reliabilities or failure rates depend on maintenance decisions \citep{kobbacy2008complex,zhu2021multicomponent}. 
To compute two-stage RO with DDU, several exact and approximation algorithms have appeared in the literature, including those based on decision rules 
\cite{zhang2020unified}, Benders decomposition \cite{zhang2022two}, parametric programming \cite{avraamidou2020adjustable}, and basic C\&CG~\cite{chen2022robust}. We note that parametric C\&CG, a new extension and generalization of basic C\&CG, demonstrates a strong and accurate solution capacity when the DDU set is a polyhedron and the recourse problem is an LP. Moreover, it can, with some rather simple modifications, approximately solve two-stage RO with mixed integer recourse. Especially, its approximation quality can be quantitatively evaluated, as both lower and upper bounds of the objective function value are generated during the execution of the algorithm. 

Nevertheless, two-stage RO with DDU that contains mixed integer structure(s) in the uncertainty set and/or the recourse problem has not been studied systematically. To the best of our knowledge, no exact algorithm has been developed to deal with the challenge arising from those types of mixed integer structures for two-stage RO with DDU.  Note that a mixed integer DDU set has a much more powerful capacity in capturing sophisticated interactions between the first-stage decision and the underlying randomness, which are often beyond that of any convex set. Also, with a DDU set that is varying, we anticipate that handling mixed integer recourse will be very demanding, compared to its counterpart defined with a fixed DIU set.   Hence, in this paper, we investigate mixed integer structures within DDU-based two-stage RO, and develop exact  and approximation solution algorithms. Note that they address some critical limitations of existing counterparts, including the relatively complete recourse assumption made for mixed integer recourse, and greatly expand our capacity to handle more complex robust optimization problems. 

Actually, with those algorithms and potentially more powerful ones in the future, it is worth mentioning a great opportunity to delve deeper into an interesting and promising  observation made in \cite{zeng2022two}. It suggests that a DIU set can be converted to a DDU one by taking advantage of domain expertise or structural properties, which allows us to compute RO with a significantly improved computational performance. So, in addition to DDUs that occur naturally, we probably can use a mixed integer set to explicitly represent  hidden yet more fundamental connection between decisions and worst case scenarios in a DIU/DDU set,  converting it to a (simpler) mixed integer DDU one. Subsequently, we can leverage the aforementioned algorithms to facilitate faster computations for large-scale instances.

The remainder of this paper is organized as follows. In Section \ref{sect_general_formulation}, we introduce  general mathematical formulations of two-stage RO with mixed integer DDU and recourse, and present a set of basic properties. In Section  \ref{sec_MIDDU}, we consider mixed integer DDU and present a new variant of parametric C\&CG to exactly solve this type of two-stage RO, along with analyses on its convergence and iteration complexity.  Similarly, Section \ref{sect_MIPrecourse} focuses on mixed integer recourse, and presents and analyzes a nested variant of parametric C\&CG algorithm. In Section \ref{sect_MIDDUMIP}, we also describe algorithm operations on top of  those variants of parametric C\&CG, to exactly and approximately handle the most complex formulation with mixed integer structures in both DDU set and recourse. Section \ref{sec_MIPC} demonstrates the modeling strength of mixed integer structures and the performances of those algorithms on variants of the robust facility location model. The whole paper is concluded by Section~\ref{sect_conclusion}.  

\section{Formulations and Basic Properties}
\label{sect_general_formulation}
In this section, we first introduce the general mathematical formulation of DDU-based two-stage RO with mixed integer structures. Then, we present a set of  properties that deepen our understanding and analyses on this modeling and optimization paradigm.

\subsection{The General Formulation}
In a two-stage decision making procedure, the decision maker initially determines the value of the first-stage decision variable $\mathbf x$ before the materialization of random factor $\mathbf u$. Then, after the uncertainty is cleared, she has an opportunity to make a recourse decision for mitigation, which, however, is restricted by her choice of $\mathbf x$  and $\mathbf{u}$'s realization. Specifically, let $\mathbf x=(\mathbf x_c,\mathbf x_d)$ denote the first-stage decision variable vector with $\mathbf x_c$ and $\mathbf x_d$ representing  its continuous   and discrete components, respectively. Similarly, $\mathbf y=(\mathbf y_c,\mathbf y_d)$ and $\mathbf{u}=(\mathbf{u}_c, \mathbf{u}_d)$ denote the recourse decision variable vector and the uncertainty variable vector, respectively, both of which could contain continuous and discrete variables. The general mathematical formulation of  two-stage RO with DDU is
\begin{align}
	\label{eq_2RO}
	\mathbf{2-Stg \ RO}: \ \ \ \mathit{w}^* = & \min_{\mathbf x\in \mathcal{X}}\mathbf c_1\mathbf x+\max_{\mathbf u\in\mathcal U(\mathbf x)} \ \min_{\mathbf y\in\mathcal Y(\mathbf x,\mathbf u)}\mathbf c_2\mathbf y,
\end{align}
\textrm{where}  $ \mathcal{X} =\{\mathbf{x}\in \mathbb{R}^{n_x}_+ \times \mathbb{Z}^{m_x}_+: \mathbf{Ax}\geq \mathbf{b}\}$,
$\mathcal U(\mathbf x)$ is the uncertainty set with
\begin{align}
	\label{eq_uncer_set}
	\mathcal U(\mathbf x)=\left\{\bf u\in \mathbb{R}^{n_u}_+ \times \mathbb{U}_d: \mathbf F_c(\mathbf x)\mathbf u_c+ \mathbf F_d(\mathbf x)\mathbf u_d\leq\mathbf h+\mathbf{Gx}\right\},
\end{align}
and $\mathcal Y(\mathbf x,\mathbf u)$ is the feasible set of the recourse problem as in the next.
\begin{align}
	\label{eq_recourse_set}
	\mathcal Y(\mathbf x,\mathbf u)=\left\{\mathbf y\in\mathbb{R}^{n_y}_+\times \mathbb{Y}_d: \mathbf B_{2,c}\mathbf y_c+\mathbf B_{2,d}\mathbf y_d\geq\mathbf d-\mathbf B_1\mathbf x-\mathbf E_c\mathbf u_c - \mathbf E_d\mathbf u_d\right\}
\end{align}
Note that $\mathbb U_d\subseteq \mathbb Z^{m_u}_+$ and $\mathbb Y_d\subseteq \mathbb Z^{m_y}_+$ are two bounded sets of discrete structures. Coefficient vectors $\mathbf{c}_1$, $\mathbf{c}_2$ (both are row vectors), $\mathbf{b}$, $\mathbf h$, $\mathbf d$, and matrices, $\mathbf{A}$, $\mathbf F(\mathbf x)=[\mathbf F_c(\mathbf x) \ \mathbf F_d(\mathbf x)]$,
$\mathbf{G}$, $\mathbf B_1$, $\mathbf B_2=[\mathbf B_{2,c} \ \mathbf B_{2,d}]$, and $\mathbf E=[\mathbf E_c \ \mathbf E_d]$ are all with appropriate dimensions.

Regarding $\mathcal U(\mathbf x)$, the DDU set defined by a point-to-set mapping in \eqref{eq_uncer_set}, we consider  two types of decision-dependence: right-hand-side (RHS) dependence and left-hand-side (LHS) dependence. The former one has $\mathbf x$ appeared in RHS of \eqref{eq_uncer_set} only, while the latter one has $\mathbf x$ in \eqref{eq_uncer_set}'s LHS  only. If $\mathcal{U}(\mathbf x)$  have both RHS and LHS dependence, the dependence can be converted into an LHS one by appending $-\mathbf {Gx}$ as a column to $\mathbf{F(x)}$ and extending $\mathbf u$ with one more dimension that takes value $1$. Obviously,  \eqref{eq_2RO} subsumes the classical DIU-based two-stage RO model  \citep{ben2004adjustable,bertsimas2010optimality} by setting $\mathcal U(\mathbf x)$ to a fixed set $\mathcal{U}^0$ for all $\mathbf x\in \mathcal{X}$.

We employ the definition of \textit{equivalence} presented in \cite{zeng2022two}, which says that two formulations are equivalent to each other if they share the same optimal value, and one's optimal first stage solution is also optimal to the other one. Also, we adopt three very mild assumptions from \cite{zeng2022two} for the remainder study, which practically do not impose restriction. \\
\noindent $(\textit {A1})$ For any $\mathbf{x}\in \mathcal{X}$, $\mathcal{U}(\mathbf x)\neq \emptyset$; \\
\noindent  $(\textit{A2})$ $\mathcal{U}(\mathbf x)$ is a bounded set, i.e., for any given $\mathbf{x}\in \mathcal{X}$, $u(\mathbf x)_j<\infty$ \ $\forall j$;\\
\noindent $(\textit{A3})$ The next monolithic MIP  has a finite optimal value.
\begin{eqnarray}
	\label{eq_MIP_RO}
	\mathit{w}_R=\min\Big\{\mathbf{c_1x+c_2y}: \mathbf{x}\in \mathcal{X}, \mathbf{u}\in \mathcal{U}(\mathbf x), \mathbf{y}\in \mathcal{Y}(\mathbf x, \mathbf u)\Big\}> -\infty
\end{eqnarray}

As noted in \cite{zeng2022two}, (\textit{A1}) substantiates the two-stage decision making framework, and can be simply achieved by the next formulation that again is in the form of \eqref{eq_2RO}.  
$$		\min_{\mathbf x\in\mathcal X, \mathbf{u}'\in \mathcal{U}(\mathbf x)}\mathbf c_1\mathbf x+\max_{\mathbf u\in\mathcal U(\mathbf x)} \ \min_{\mathbf y\in\mathcal Y(\mathbf x,\mathbf u)}\mathbf c_2\mathbf y$$

Also, (\textit{A2}) generally holds for practical random factor. We say $\mathbf x^0\in \mathcal{X}$ is \textit{infeasible} if the recourse problem is infeasible for some $\mathbf u\in \mathcal{U}(\mathbf x^0)$, and otherwise it is feasible. The whole model $\mathbf{2-Stg \ RO}$ is infeasible if no first-stage decision in $\mathcal{X}$ is feasible. 
Regarding (\textit{A3}), it is  clear that $w_R$ from \eqref{eq_MIP_RO} yields a lower bound to $w^*$,  since \eqref{eq_MIP_RO} is a simple relaxation of~\eqref{eq_2RO}. Hence, if \eqref{eq_MIP_RO} is infeasible, so is $\mathbf{2-Stg \ RO}$. 

\subsection{Basic Properties}

Let $w(\mathcal{X}, \mathcal{U}(\mathbf x), \mathcal{Y}(\mathbf{x},\mathbf{u}))$ denote the optimal value of \eqref{eq_2RO} defined on those sets.  Next, we present some simple results of $\mathbf{2-Stg \ RO}$ in \eqref{eq_2RO} by varying the underlying sets, which help us derive more sophisticated relaxations.   With a slight abuse of notation and when either $\mathcal{X}$, $\mathcal U(\mathbf x)$ or $\mathcal Y(\mathbf x, \mathbf u)$ is the original one without any modification, we simply use $\cdot$ to replace it. For example, by definition, $\mathit{w}^*=\mathit{w}\big(\mathcal{X}, \mathcal{U}(\mathbf x), \mathcal{Y}(\mathbf x,\mathbf u)\big) = \mathit w (\cdot, \cdot, \cdot)$.

\begin{prop}
	\label{prop_relax_1}
	\begin{description}
		\item[$(i)$] Consider two sets $ \emptyset\neq \mathcal{U}^1(\mathbf x)\subseteq \mathcal{U}^2(\mathbf x)$ $\forall \mathbf{x}\in \mathcal{X}$. We have\\
		$\mathit{w}\big(\cdot, \mathcal{U}^1(\mathbf x), \cdot\big)\leq \mathit{w}\big(\cdot, \mathcal{U}^2(\mathbf x), \cdot\big).$ 
		\item[$(ii)$]  Consider two sets $\mathcal{Y}^1(\mathbf x, \mathbf u)\subseteq \mathcal{Y}^2(\mathbf x, \mathbf u)$ $\forall \mathbf{x}\in \mathcal{X}, \mathbf{u}\in \mathcal{U}(\mathbf x)$. We have\\
		$\mathit{w}\big(\cdot, \cdot, \mathcal Y^2(\mathbf x,\mathbf u)\big)\leq \mathit{w}\big(\cdot, \cdot, \mathcal Y^1(\mathbf x,\mathbf u)\big).$ 
		\item[$(iii)$] Let $\mathcal{X}_r$, $\mathcal{U}_r(\mathbf x)$ and $\mathcal{Y}_r(\mathbf x,\mathbf u)$ denote some relaxations of $\mathcal{X}$, $\mathcal{U}(\mathbf x)$ and $\mathcal{Y}(\mathbf x,\mathbf u)$, respectively. We have \\
		$
		\mathit{w}\big(\mathcal{X}_r,\cdot,\mathcal{Y}_r(\mathbf x,\mathbf u)\big)\leq \mathit{w}\big(\cdot,\cdot,\mathcal{Y}_r(\mathbf x,\mathbf u)\big)\leq \mathit{w}^*\leq
		\mathit{w}\big(\cdot,\mathcal{U}_r(\mathbf x),\cdot\big).$ 
		\item[$(iv)$] Consider two point-to-set mappings  $\mathcal{U}^1(\mathbf x)\neq \emptyset$ and $\mathcal{U}^2(\mathbf x)\neq\emptyset$ such that
		$\mathcal{U}^1(\mathbf x)\cup \mathcal{U}^2(\mathbf x)\subseteq \mathcal{U}(\mathbf x)$ for all
		$\mathbf{x}\in \mathcal{X}$, and $\mathcal{U}^3(\mathbf x)$
		such that  $\mathcal{U}(\mathbf x)\subseteq \mathcal{U}^3(\mathbf x)$
		for all $\mathbf x\in \mathcal{X}$. Together with~\eqref{eq_MIP_RO}, we have \\	
		$
		w_R\leq \max\big\{\mathit{w}\big(\cdot,\mathcal{U}^1(\mathbf x),\cdot\big), \mathit{w}\big(\cdot,\mathcal{U}^2(\mathbf x),\cdot\big)\big\}  \leq \mathit{w}\big(\cdot,\mathcal{U}^1(\mathbf x)\cup \mathcal{U}^2(\mathbf x),\cdot\big) \leq  \mathit{w}^*
		\leq \mathit{w}\big(\cdot,
		\mathcal{U}^3(\mathbf x),\cdot\big). \hfill\square$
	\end{description}
\end{prop}

In this paper, the result in Proposition \ref{prop_relax_1}.$(iii)$ is specific to their continuous relaxations, although it is valid for any relaxation. Indeed, it can be strengthened under some special situation. Consider cases where $\mathcal U(\mathbf x)$ or $\mathcal Y(\mathbf x, \mathbf u)$ is an integer set and the associated constraint  matrix is totally unimodular (TU). Then, if its RHS is integral,  computing an integer program can reduce to solving its linear program relaxation.

\begin{prop}
	\label{prop_extremepoint}
	$(i)$ If $n_u=m_y=0$ and $\mathbf{F}_d(\mathbf x)$ is TU and $\bf h+\bf G\bf x$ is an integral vector for $\mathbf x\in \mathcal{X}$, we say the DDU has the TU property, and we have $\mathit{w}^*=\mathit{w}\big(\cdot, \mathcal{U}_r(\mathbf x), \cdot\big)$. \\
	$(ii)$ Suppose that $n_u=n_y=0$, $\mathbf{B}_{2,d}$  is TU, and $\mathbf d-\mathbf B_1\mathbf x - \mathbf E_d\mathbf u_d$ is an integral vector for $\mathbf x\in \mathcal{X}$ and $\bf u_d\in \cal U(\bf x)$, we have  $\mathit{w}^*=\mathit{w}\big(\cdot, \cdot, \mathcal Y_r(\mathbf x, \mathbf u)\big)$.  $\hfill\square$
\end{prop}
Result in Proposition \ref{prop_extremepoint}.$(i)$ can be proven by using the duality of the recourse problem, noting that  it is a linear program. Then, according to \cite{zeng2022two}, optimal $\mathbf u$ to the max-min substructure is always an optimal solution of integer program $\max\{\mathbf c^0\mathbf u: \mathbf u\in \mathcal U(\mathbf x)\}$ for some $\mathbf c^0$, which, by the TU property, simply reduces to a linear program with $\mathcal U(\mathbf x)$ replaced by $\mathcal U_r(\mathbf x)$. Result in Proposition \ref{prop_extremepoint}.($ii$) is straightforward by the TU property.

As a special case of the point-to-set mapping, we consider a point-to-point mapping $f^u$ satisfying $f^u(\mathbf x)\in \mathcal U(\mathbf x)$ for $\mathbf x\in \mathcal X$. Using this mapping and following Proposition~\ref{prop_relax_1}, we probably can construct a non-trivial but more tractable single-level relaxation for $\mathbf{2-Stg \ RO}$.

\begin{cor}
\label{cor_relaxation_f_mapping}
The following formulation is a relaxation of $\mathbf{2-Stg \ RO}$ in \eqref{eq_2RO} that is stronger than \eqref{eq_MIP_RO}, i.e.,
\begin{eqnarray*}
\mathit{w}_R \leq \min\big\{\mathbf{c_1x+c_2y}: \mathbf{x}\in \mathcal{X},\ \mathbf{u}=f^u(\mathbf x),\ \mathbf{y}\in \mathcal{Y}(\mathbf x, \mathbf u)\big\} \leq \mathit{w}^*, 
\end{eqnarray*} recalling $w_R$ is the optimal value of \eqref{eq_MIP_RO}.  $\hfill\square$
\end{cor}

The assumption made in Corollary \ref{cor_relaxation_f_mapping} might be restrictive to derive a strong relaxation. Instead of fixing $\mathbf u$ to $f^u(\mathbf x)$ completely, we can fix it partially according to that mapping, which yields a flexible strategy in deriving a strong relaxation. Consider $
\hat J\subseteq J=\{1,\dots, n_u+m_u\}$ and define a new DDU set that fixes $u_{j}$ for $j\in \hat J$, i.e.,   
\begin{equation*}
	\label{eq_U_J}
	\mathcal{U}^{f^u}(\mathbf x|\bf u_{\hat J})= \mathcal{U}(\mathbf x)\cap\{\mathbf u\in \mathbb{R}^{n_u}_+ \times \mathbb{U}_d: u_{j}= f^u(\bf x)_{j},\ \forall j\in \hat J\}.
\end{equation*}

Again, the next result follows from Proposition \ref{prop_relax_1}. 
\begin{cor}
	\label{cor_relaxation_f_mapping2}
	For a given $f^u(\mathbf x)$ and $\hat J_1\subseteq \hat J_2\subseteq J$, we have $\mathcal U^{f^u}(\mathbf x|\bf u_{\hat J_2})\subseteq \mathcal U^{f^u}(\mathbf x|\bf u_{\hat J_1})$ and 
	\begin{eqnarray*} w_R\leq \mathit w(\cdot,\mathcal U^{f^u}(\mathbf x|\bf u_{\hat J_2}), \cdot)\leq \mathit w(\cdot,\mathcal U^{f^u}(\mathbf x|\bf u_{\hat J_1}), \cdot) \leq  \mathit w(\cdot, \mathcal U^{f^u}(\mathbf x), \cdot) = \mathit w^*. \pushQED{\qed}\qedhere
	\end{eqnarray*} 
\end{cor}

Result in Corollary \ref{cor_relaxation_f_mapping2}, which subsumes Corollary \ref{cor_relaxation_f_mapping}, can be easily generalized by introducing multiple mappings. Suppose that we have two mappings, i.e.,
$f^u_i(\mathbf x)\in \mathcal U(\mathbf x)$ for $\mathbf x\in \mathcal X$ and $i=1,2$. Also, consider subsets $\hat J_{i,1}\subseteq \hat J_{i,2}\subseteq J$, $i=1,2$, noting that it is not necessary to have $\hat J_{1,j}\cap \hat J_{2,j}=\emptyset$ for $j=1,2$. 

\begin{cor}
	\label{cor_relaxation_2f_mapping2}
	For given $f^u_i(\mathbf x)$, $i=1, 2$ and $\hat J_{i,1}\subseteq \hat J_{i,2}\subseteq J$, we have
	\begin{eqnarray*} w_R\leq \mathit w(\cdot,\mathcal U^{f^u_1}(\mathbf x|\bf u_{\hat J_{1,2}})\cup U^{f^u_2}(\mathbf x|\bf u_{\hat J_{2,2}}), \cdot)\leq \mathit w(\cdot,\mathcal U^{f^u_1}(\mathbf x|\bf u_{\hat J_{1,1}})\cup \mathcal U^{f^u_2}(\mathbf x|\bf u_{\hat J_{2,1}}), \cdot)  \leq \mathit w^*. \pushQED{\qed}\qedhere
	\end{eqnarray*} 
\end{cor}

\begin{rem}	
Results in Proposition \ref{prop_relax_1} and Corollaries \ref{cor_relaxation_f_mapping}-\ref{cor_relaxation_2f_mapping2} provide a strategy to handle mixed integer DDU approximately by using simpler DDU sets. For example, as shown in Corollary~\ref{cor_relaxation_2f_mapping2}, different deep insights or domain knowledge can be used to design multiple $f^u(\mathbf x)$'s so that they can jointly yield a strong approximation. 
\end{rem}

\subsection{A Discussion on Mixed Integer DDU}  
\begin{figure}[ht!]
	\centering
	\includegraphics[scale=0.45]{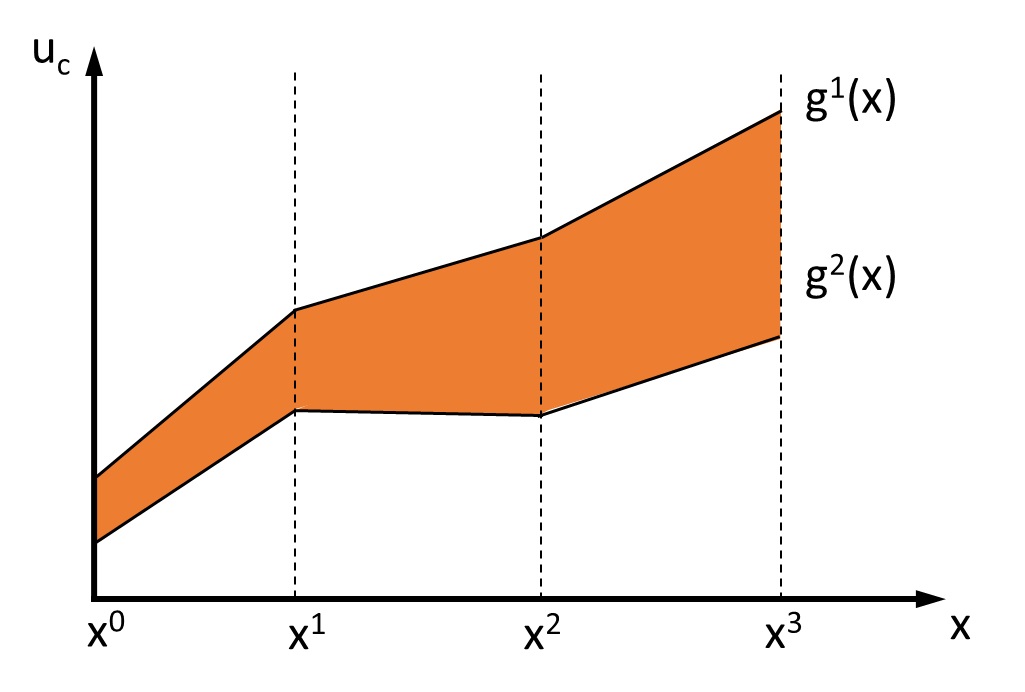}
	\caption{A DDU Set of \eqref{eq_piecewiseDDU}}
	\label{fig_MIP_DDU}
\end{figure}

  As noted earlier, introducing mixed integer structure drastically improves the modeling capacity of DDU. Sophisticated connection between the first-stage decision and behavior of the random factor, which is often impossible to be described by any convex DDU set, can be well captured by employing discrete, especially binary, and continuous variables.  Nevertheless, as shown in the next section, a mixed integer DDU set could be more complex to analyze and compute. Hence, it should be preferred to reduce the number of discrete variables needed to describe $\mathcal{U}(\mathbf x)$. Some special cases may even allow us to eliminate  discrete variables from $\mathcal{U}(\mathbf x)$. For example, consider the uncertainty set depicted in Figure \ref{fig_MIP_DDU} that is a function of $x$, which is represented by the following formulation. 
\begin{subequations}
	\label{eq_piecewiseDDU}
	\begin{align}
		\mathcal{U}(x) = \Big\{&(u_c, \mathbf u_d, \mathbf{t}, \tilde g_1, \tilde g_2): \ \tilde g_1= \sum_{j=0}^3 g^1(x^j)t_j; \ \tilde g_2= \sum_{j=0}^3 g^2(x^j)t_j; \ x= \sum_{j=0}^3 x^j t_j; \  \label{eq_ddu-1stg1}  \\
		& \tilde g_1\leq u_c\leq \tilde g_2; \  t_0\leq u_{d,1}; \ t_j\leq u_{d,j}+u_{d,j+1}, j=1,2; \
		t_3\leq u_{d,3}; \\
		& \sum_{j=1}^3u_{d,j}=1; \ \sum_{j=0}^3t_j=1; \ u_{d,j}\in \{0,1\}, j=1,2,3; \  t_j\geq 0,  j=0,1,2,3 \label{eq_ddu-1stg3}\Big\}
	\end{align}
\end{subequations}

This DDU set is bounded by two piece-wise linear functions, i.e., $g^1(x)$ and $g^2(x)$. Although $(\bf u_d,\bf t, \tilde g_1, \tilde g_2)$ is not part of the first-stage decision, it can be seen that $(\bf u_d, \bf t,\tilde g_1, \tilde g_2)$ is uniquely determined by $x$. Under such a situation, we can augment $x$ with $(\bf u_d,\bf t,\tilde g_1, \tilde g_2)$ (and therefore annex all constraints, except $\tilde g_1\leq u_c\leq \tilde g_2$, to $\mathcal{X}$). By using  simpler DDU set $\mathcal{U}(\mathbf x)=\{\tilde g_1\leq u_c\leq \tilde g_2\}$, an equivalent but much less complex 2-stage RO can be obtained.  As a rule of thumb, if some variables, especially discrete ones, can be uniquely determined by $\mathbf x$, they can be safely moved from $\mathcal U(\mathbf x)$ to $\mathcal{X}$. We next formalize this result by rewriting $\mathbf u$ as $(\mathbf u_0, \mathbf u_1)$, assuming that those sub-vectors may contain continuous or/and discrete variables.  

\begin{prop}
\label{prop_DDUto1stg}
Suppose that $$\mathcal U(\mathbf x)=\Big\{(\mathbf u_0, \mathbf u_1)\in (\mathbb{R}^{n^0_u}_+ \times \mathbb{U}_{d_0})\times(\mathbb{R}^{n^1_u}_+ \times \mathbb{U}_{d_1}): \mathbf u_0 = \tilde f(\mathbf x), \mathbf F_0(\mathbf x) \mathbf u_0+\mathbf{F}_1(\mathbf x)\mathbf u_1\leq \mathbf h+\mathbf G\mathbf x\Big\},$$ where $\tilde f(\mathbf x)$ returns a unique value for $\mathbf u_0$ for any fixed $\mathbf x$.  Then, $\mathbf{2-Stg \ RO}$ in \eqref{eq_2RO} is equivalent to the following two-stage RO formulation
\begin{eqnarray*}
	\label{eq_2RO_reduced}
	\min_{\mathbf x\in \mathcal{X},\mathbf u_0=\tilde f(\mathbf x)}\mathbf c_1\mathbf x+\max_{\bf u_1\in\mathcal U_1(\mathbf x,\bf u_0)} \ \min_{(\mathbf y_c,\mathbf y_d)\in\mathcal Y(\mathbf x,\mathbf u)}\mathbf c_2\mathbf y, 
\end{eqnarray*}
where $\mathcal U_1(\mathbf x,\bf u_0)=\{\mathbf u_1\in \mathbb{R}^{n^1_u}_+ \times \mathbb{U}_{d_1}: \mathbf F_1(\mathbf x)\mathbf u_1 \leq \mathbf h+\mathbf G\mathbf x - \mathbf F_0(\mathbf x)\mathbf u_0 \}$. $\hfill\square$ 
\end{prop}
Even if $\tilde f(\bf x)$ does not exist or is not known,  together with Corollary \ref{cor_relaxation_f_mapping}, we have a non-trivial relaxation by simply augmenting $\bf x$ with $\bf u$ and adopting $\mathcal U_1(\mathbf x,\bf u_0)$ as the DDU set. 
\begin{cor}
Let $\bf u'_1$ be a replicate of $\bf u_1$. We have
\begin{eqnarray}
	\label{eq_2RO_U_split}
	w_R\leq w_{(\bf u_0,\bf u_1)}\equiv\min_{\mathbf x\in \mathcal{X},(\mathbf u_0,\bf u'_1)\in \cal{U}(\bf x)}\mathbf c_1\mathbf x+\max_{\bf u_1\in\mathcal U_1(\mathbf x,\bf u_0)} \ \min_{(\mathbf y_c,\mathbf y_d)\in\mathcal Y(\mathbf x,\mathbf u_0,\bf u_1)}\mathbf c_2\mathbf y\leq w^*, 
\end{eqnarray}
and if $w_{(\bf u_0,\bf u_1)}=+\infty$, so is $w^*$. \pushQED{\qed}\qedhere
\end{cor}
Hence, if one variable in DDU set cannot be uniquely determined by $\bf x$, it should not be relocated to $\mathcal X$. Otherwise, the resulting formulation is generally a relaxation to the original $\mathbf{2-Stg \ RO}$. It is also interesting to compare two formulations: the aforementioned one for $w_{(\bf u_d,\bf u_c)}$ that defines a lower bound to $w^*$, and the one for $\mathit{w}\big(\cdot,\mathcal{U}_r(\mathbf x),\cdot\big)$, which considers the continuous relaxation of $\cal U(\bf x)$ and yields an upper bound for $w^*$.

\section{Two-Stage RO with Mixed Integer DDU Set}
\label{sec_MIDDU}

In this section, we consider $\mathbf{2-Stg \ RO}$ with a mixed integer DDU set and a recourse problem in the form of linear program. Our focus is on the derivation of computationally more friendly equivalences,  the development of a C\&CG type of algorithm and its convergence analysis. According to Proposition \ref{prop_DDUto1stg}, we assume  without loss of generality that discrete variables in $\mathcal U(\mathbf x)$ cannot be uniquely determined by the first-stage decision  $\mathbf x$.

\subsection{Enumeration Based Reformulations}
\label{ssec_EnReform}
	
Given the recourse problem is a linear program, we present its dual problem for a given $\bf u$. 
\begin{align}
	\label{eq_recourse_dual}
	\max\Big\{(\mathbf d-\mathbf B_1\mathbf x-\mathbf E\mathbf u)^\intercal\pi: \mathbf B^\intercal\pi\leq \mathbf c_2^\intercal,\ \pi\geq \mathbf 0\Big\}=\max\Big\{(\mathbf d-\mathbf B_1\mathbf x-\mathbf E\mathbf u)^\intercal\pi:\pi\in \Pi\Big\}
\end{align}
where $\Pi=\{\mathbf B^\intercal\pi\leq \mathbf c_2^\intercal,\ \pi\geq \mathbf 0\}$. Let $\mathcal{P}_{\Pi}$ be the set of extreme points and $\mathcal{R}_{\Pi}$ the set of extreme rays of $\Pi$, respectively. Clearly,  an optimal solution of \eqref{eq_recourse_dual} belongs to $\mathcal{P}_{\Pi}$ if it is finite, or it is unbounded following some direction in $\mathcal R_{\Pi}$.  
Next, we generalize a result of~\cite{zeng2022two}, with minor modifications, from polytope $\mathcal U(\mathbf x)$ to the case where it is a mixed integer set.

\begin{thm}
	\label{thm_OptEquiv}[Adapted from \cite{zeng2022two}]
	Formulation $\mathbf{2-Stg \ RO}$ in \eqref{eq_2RO} (and its other  equivalences) is equivalent to a bilevel linear optimization program as in the following.
	\begin{subequations}
		\label{eq_2stgRO_PI_OM}
		\begin{align}
			\mathbf{2\!-\!Stg \ RO(\Omega_{\Pi})}: \   w^*=\min \ & \ \mathbf{c}_1\mathbf x+ \eta   \label{equiv_obj}\\
			\mathrm{s.t.} \ & \ \mathbf x\in \mathcal{X}\\
			&  \Big\{\eta\geq  \mathbf{c_2}\mathbf y, \ \mathbf y\in \mathcal Y(\mathbf x, \mathbf u)\notag \\
			&  \  \ \mathbf{u}\in \arg\max\big\{(-\mathbf E\mathbf u)^\intercal\pi:\mathbf u\in \mathcal{U}(\mathbf x)\big\}\Big\}  \ \forall \pi\in \mathcal{P}_{\Pi} \label{eq_enu2_optimality}\\
			& \Big\{ \mathbf y\in \mathcal Y(\mathbf x, \mathbf v)\notag \\
			&  \ \ \mathbf{v}\in \arg\max\big\{(-\mathbf {Ev})^\intercal\gamma: \mathbf v\in \mathcal{U}(\mathbf x)\big\}\Big\}  \ \forall \gamma\in \mathcal{R}_{\Pi}\label{eq_enu2_feasibility}
		\end{align}
	\end{subequations}
\end{thm}

Unlike the polytope case in \cite{zeng2022two}, with the mixed integer structure of $\mathcal U(\mathbf x)$, the optimal solution set of the optimization problem in \eqref{eq_enu2_optimality} or \eqref{eq_enu2_feasibility} cannot be represented by any known optimality conditions. Nevertheless, once the discrete variables are fixed, such problems reduce to linear programs. It is thus tempting to enumerate discrete variables so that only the continuous portion of $\mathcal U(\mathbf x)$ is under consideration.  To achieve this, we first define an extension of $\mathcal U(\mathbf x)$ in the following, where $\mu_u$ denotes the number of constraints of $\mathcal U(\mathbf x)$.
\begin{align*}
	\tilde{\mathcal U}(\mathbf x)\equiv\Big\{(\mathbf u_c,\mathbf u_d, \tilde {\mathbf u})\in \mathbb{R}^{n_u}_+ \times \mathbb{U}_d\times \mathbb{R}^{\mu_u}_+: \mathbf F_c(\mathbf x)\mathbf u_c+ \mathbf F_d(\mathbf x)\mathbf{u}_d-\tilde{\bf u}\leq\mathbf h+\mathbf{Gx}\Big\}
\end{align*}
Also, let $\tilde{\mathcal U}(\bf x|\bf{u}^0_d)\equiv\tilde{\mathcal U}(\bf x)\cap\{(\mathbf u_c,\mathbf u_d, \tilde {\mathbf u})\in \mathbb{R}^{n_u}_+ \times \mathbb{U}_d\times \mathbb{R}^{\mu_u}_+|\mathbf u_d=\mathbf u_d^0\}$. Due to the existence of $\tilde{\bf u}$,  $\tilde{\mathcal U}(\bf x|\bf{u}^0_d)$ is not empty for any fixed $(\mathbf x,\mathbf u^0_d)\in \mathcal{X}\times\mathbb U_d$.

Without loss of generality, assume that $\mathbb{U}_d=\{\mathbf u^1_d, \dots, \mathbf u^T_d\}$. Given fixed $(\mathbf x,\mathbf u_d^t)$ with $\mathbf x\in \mathcal X$ and $\mathbf u_d^t\in \mathbb{U}_d$, we abstract the two optimization problems in  \eqref{eq_enu2_optimality} and \eqref{eq_enu2_feasibility} into the following linear program defined with respect to $\tilde{\mathcal U}(\bf x|\bf{u}^t_d)$. With $\mathbf 1 \equiv(1,\dots, 1)$, i.e., a column vector with all its entries being ones, we have
\begin{eqnarray}
	\label{eq_LP_parametric_ed}
	\begin{split}
		\mathbf{LP}(\mathbf x, \mathbf u^t_d, \boldsymbol\beta): \ \max\Big\{
		(-\mathbf{Eu}_c)^\intercal\boldsymbol\beta-M(\mathbf 1^\intercal\tilde{\bf u}): \ & \mathbf F_c(\mathbf x)\mathbf u_c - \tilde{\bf u}\leq\mathbf h+\mathbf{Gx}-\mathbf F_d(\mathbf x)\mathbf{u}^t_d, \\
		& \mathbf u_c\geq \mathbf 0, \ \tilde{\bf u}\geq \mathbf 0 \Big\}.
	\end{split}
\end{eqnarray}
Note that $M$ is a sufficiently large number, and coefficient $\boldsymbol\beta$ could be either an extreme point or an extreme ray of $\Pi$. The next result follows trivially. 
\begin{lem}
	\label{lem_FeasDisU}
	For a given $\mathbf x$, $\mathbf u^t_d$ and $\boldsymbol\beta$, if
	\begin{align*}
	\mathcal U(\mathbf x|\mathbf u_d^t)&\equiv\mathcal U(\mathbf x)\cap\{(\mathbf u_c,\mathbf u_d)\in \mathbb{R}^{n_u}_+ \times \mathbb{U}_d|\mathbf u_d=\mathbf u^t_d\}\\
	 &=\left\{\mathbf u_c\in \mathbb{R}^{n_u}_+: \mathbf F_c(\mathbf x)\mathbf u_c\leq\mathbf h+\mathbf{Gx}-\mathbf F_d(\mathbf x)\mathbf u^t_d \right\} \neq \emptyset, 
	\end{align*}
	we have $\tilde{\bf u}=\mathbf 0$ in any optimal solution to \eqref{eq_LP_parametric_ed}. Otherwise, we have $\tilde{\bf u}\geq \mathbf 0$ satisfying $\sum_{j=1}^{\mu_u} \tilde u_j>0$ for  its all optimal solutions. \hfill $\square$
\end{lem}
Noting that \eqref{eq_LP_parametric_ed} is a linear program that is always finitely  optimal,  its optimal solution set can be characterized by its optimality conditions. Similar to the derivation presented in~\cite{zeng2022two}, we next use  $\mathbf{LP}(\mathbf x, \mathbf u_d, \pi)$'s KKT conditions to define the optimal solution set, which is denoted by $\mathcal{OU}(\mathbf x, \mathbf u_d, \pi)$. 
\begin{eqnarray}
	\label{eq_KKT_point}
	\mathcal{OU}(\mathbf x, \mathbf u_d, \pi) =  \left\{\begin{array}{l}
		\mathbf F_c(\mathbf x)\mathbf u_c^{(\bf u_d,\pi)} -\tilde{\bf u}^{(\bf u_d,\pi)}\leq \mathbf h+\mathbf{Gx} - \mathbf F_d(\mathbf x)\mathbf{u}_d \\ 
		\mathbf F_c(\mathbf x)^\intercal\lambda^{(\bf u_d,\pi)}\geq - \mathbf E^\intercal\pi,  \ \lambda^{(\bf u_d,\pi)}\leq M\mathbf 1 \\
		\lambda^{(\bf u_d,\pi)} \circ (\mathbf h+\mathbf{Gx}- \mathbf F_d(\mathbf x)\mathbf{u}_d - \bf F_c(\mathbf x)\mathbf u^{(\bf u_d,\pi)}_c +\tilde{\bf u}^{(\bf u_d,\pi)})=\mathbf{0} \\
		\mathbf u^{(\bf u_d,\pi)}_c \circ (\mathbf F_c(\mathbf x)^\intercal\lambda^{(\bf u_d,\pi)}+ \mathbf E^\intercal\pi) = \mathbf 0,  \  \tilde{\bf u}^{(\bf u_d,\pi)}\circ(M\mathbf 1 - \lambda^{(\bf u_d,\pi)}) = \mathbf 0 \\
		\mathbf u^{(\bf u_d,\pi)}_c\geq \mathbf{0}, \ \tilde{\bf u}^{(\bf u_d,\pi)}\geq \mathbf 0,  \ \lambda^{(\bf u_d,\pi)}\geq \mathbf{0}
	\end{array}\right\}
\end{eqnarray}
Similarly, when $\boldsymbol\beta$ is a ray, i.e., $\gamma$ in $\cal R_{\Pi}$, we use $\mathcal{OV}(\mathbf x, \mathbf u_d, \gamma)$ to denote  the corresponding optimal solution set of \eqref{eq_LP_parametric_ed}, and let $\mathbf v^{(\bf u_d,\gamma)}_c$, $\tilde{\mathbf v}^{(\bf u_d,\gamma)}$, and $\zeta^{(\bf u_d,\gamma)}$ be the aliases of $\mathbf u_c$, $\tilde{\bf u}$, and $\lambda$, respectively, to represent $\mathcal{OV}$. Then, $\mathcal{OV}(\mathbf x, \mathbf u_d, \gamma)$  can be     
simply obtained by employing variables $\mathbf v^{(\bf u_d,\gamma)}_c$, $\tilde{\mathbf v}^{(\bf u_d,\gamma)}$, and $\zeta^{(\bf u_d,\gamma)}$ in~\eqref{eq_KKT_point}.  Note also that primal-dual optimality condition can be used to define 
$\mathcal{OU}(\mathbf x, \mathbf u_d, \pi)$ and $\mathcal{OV}(\mathbf x, \mathbf u_d, \gamma)$ as well. 

Next, we introduce an indicator function $\dot{u} (\tilde{\bf u})$ based on $\tilde{\bf u}$'s value. 
\begin{eqnarray}
	\label{eq_penVar_MIDDU}
\dot{u}(\tilde{\bf u}) = \left\{ \begin{array}{ll}
	0 & \mathrm{if }\ \tilde{\bf u} = \mathbf 0\\
	+\infty & \mathrm{if }\ \tilde{\bf u} \neq \mathbf 0
\end{array} \right.
\end{eqnarray}
In the following, we present a single-level equivalence of \eqref{eq_2stgRO_PI_OM}. Note that $\mathcal{OU}$ or $\mathcal{OV}$ set introduces many variables to represent optimality conditions, while we are only concerned with $(\bf u_c,\tilde{\bf u})$ or $(\bf v_c,\tilde{\bf v})$. Hence, we simply use ``$\cdot$'' to collectively represent non-critical variables.  

\begin{cor}\label{cor_Pi_equivMIPDDU} Formulation $\mathbf{2-Stg \ RO}$ in \eqref{eq_2RO} (and its other  equivalences) is equivalent to the following single-level optimization program.
	\begin{subequations}
		\label{eq_single_MIPDDU}
		\begin{align}
			\min \ &  \mathbf{c}_1\mathbf x+ \eta\label{MP_obj} \\
			\mathrm{s.t.} \ & \mathbf x\in \mathcal{X}\\
			& \Big\{ \eta\geq  \mathbf{c_2}\mathbf y^{(\bf u_d,\pi)}-\dot{u}(\tilde{\bf u}^{(\bf u_d,\pi)}), \  \mathbf y^{(\bf u_d,\pi)}\geq \mathbf 0\label{MP_eta} \\
			&\mathbf{B}_2\mathbf y^{(\bf u_d,\pi)}+\mathbf{E}_c\mathbf{u}_c^{(\bf u_d,\pi)}+
			\dot{u}(\tilde{\bf u}^{(\bf u_d,\pi)})\geq \mathbf{d-\mathbf B}_1\mathbf x-\mathbf{E}_d\mathbf{u}_d\label{MP_DFeas}\\
			& \big(\mathbf{u}^{(\bf u_d,\pi)}_c, \tilde{\bf u}^{(\bf u_d,\pi)}, \cdot\big)\in
			\mathcal{OU}(\mathbf x,\mathbf u_d, \pi)\Big\} 
			\ \forall (\bf u_d,\pi)\in \mathbb{U}_d\times\mathcal{P}_{\Pi}\\
			&\Big\{\mathbf{B}_2\mathbf y^{(\bf u_d,\gamma)}+\mathbf{E}_c\mathbf v^{(\bf u_d,\gamma)}_c + \dot{u}(\tilde{\mathbf v}^{(\bf u_d,\gamma)})\geq \mathbf{d-\mathbf B}_1\mathbf x-\mathbf{E}_d\mathbf{u}_d\\
			& \mathbf y^{(\bf u_d,\gamma)}\geq \mathbf 0, \ \big(\mathbf{v}^{(\bf u_d,\gamma)}_c, \tilde{\mathbf v}^{(\bf u_d,\gamma)}, \cdot\big)\in \mathcal{OV}(\mathbf x, \mathbf u_d, \gamma)\Big\} \ \forall (\bf u_d,\gamma)\in \mathbb{U}_d\times\mathcal{R}_{\Pi} 
		\end{align}
	\end{subequations}
\end{cor}

\begin{cor}
	\label{cor_Pi_partial_KKTOmegaMIPDDU} 
	Let $\widehat{(\bf u_{od},\pi)}\subseteq \mathbb{U}_d\times\mathcal{P}_{\Pi}$ denote a set of $(\bf u_d,\pi)$'s and $\widehat{(\bf u_{fd},\gamma)}\subseteq \mathbb{U}_d\times\mathcal{R}_{\Pi}$ a set of $(\bf u_d,\gamma)$'s. A formulation of \eqref{eq_single_MIPDDU} defined on $\widehat{(\bf u_{od},\pi)}$ and $\widehat{(\bf u_{fd},\gamma)}$ is a relaxation to $\mathbf{2-Stg \ RO}$ in \eqref{eq_2RO}, and its optimal value is smaller than or equal to $w^*$.   \hfill $\square$
\end{cor}

Results in Corollaries \ref{cor_Pi_equivMIPDDU} and \ref{cor_Pi_partial_KKTOmegaMIPDDU} generalize their counterparts developed in \cite{zeng2022two} for the case  where $\cal U(\bf x)$ is a polytope. Nevertheless, the discrete variables in $\mathcal U(\mathbf x)$ significantly increase the complexity of this single-level equivalence.  For every $(\bf u_d,\pi)$ or $(\bf u_d,\gamma)$, we need to introduce a set of recourse variables and associated constraints. Hence, it would be desired to limit the number of discrete variables in $\cal{U}(\mathbf x)$, and convert, whenever possible, any discrete $u_{d,j}$ into a first-stage decision variable.

\subsection{Parametric C\&CG-MIU to Handle Mixed Integer DDU}
The single level equivalence presented in Corollary \ref{cor_Pi_equivMIPDDU} and the relaxation in Corollary \ref{cor_Pi_partial_KKTOmegaMIPDDU} allow us to modify standard parametric C\&CG method \citep{zeng2022two} to handle the challenge arising from the mixed integer DDU. Note that standard parametric C\&CG is an iterative procedure that computes between a master problem and a few subproblems. In the following, we first introduce subproblems built upon the mixed integer DDU set, and then present the overall procedure with the modified master problem.  
	
For a given $\mathbf x^*$, the first subproblem is constructed to check its feasibility. Recall that by definition, $\mathbf x^*$ is feasible if  the recourse problem is feasible for all scenarios in $\mathcal{U}(\mathbf x^*)$. 
\begin{eqnarray}
	\label{eq_SP1}
	\begin{split}
		\mathbf{SP1}: \ \eta_f (\mathbf x^*) = \max_{(\mathbf u_c, \mathbf u_d)\in \mathcal U(\mathbf x^*)}  & \min
		\big\{\mathbf{1}^\intercal\tilde{\mathbf y}: \mathbf{B}_2\mathbf y+\tilde{\mathbf y}\geq \mathbf d-\mathbf B_1\mathbf x^*-\mathbf E_c\mathbf u_c - \mathbf E_d\mathbf u_d \\
		& \ \ \ \ \mathbf y\geq \mathbf 0,\ \tilde{\mathbf y}\geq
		\mathbf 0 \big\}
	\end{split}
\end{eqnarray}
As $\tilde{\mathbf y}\geq \bf0$, it is clear that $\mathbf x^*$ is feasible to $\mathbf{2-Stg \ RO}$ in \eqref{eq_2RO} (and its equivalences) if and only if  $\eta_f(\mathbf x^*)=0$. In the case where $\eta_f(\mathbf x^*)=0$, we compute the second subproblem, which is the original $\max-\min$ substructure of \eqref{eq_2RO} .  
\begin{eqnarray}
	\label{eq_SP2}
	\mathbf{SP2}: \ \eta_o(\mathbf x^*)=\max_{(\bf u_c, \bf u_d)\in \mathcal{U}(\mathbf x^*)} \min\big\{
	\mathbf c_2\mathbf y: \mathbf y\in \mathcal{Y}(\mathbf x^*, \mathbf u_c,\bf u_d)\big\}
\end{eqnarray}
Let $(\bf u^*_o,\pi^*)$ denote its optimal solution, where $\bf u^*_o=(\mathbf u^*_{oc}, \mathbf u^*_{od})$ represents the worst case scenario, and $\pi^*$ the associated optimal dual solution to the recourse problem. Without loss of generality, $\pi^*\in \mathcal P_{\Pi}$. Note that $\eta_o(\bf x^*)$ is the worst case performance of $\mathbf x^*$.

In the case where $\eta_f(\mathbf x^*)>0$, its optimal solution to \eqref{eq_SP1}, denoted by $\bf u^*_f=(\mathbf u^*_{fc}, \mathbf u^*_{fd})$, causes the recourse problem infeasible. Then, we consider in the following the third subproblem, which is the dual of the recourse problem with respect to $\mathbf u^*_f$.
\begin{eqnarray}
	\label{eq_SP3}
	\mathbf{SP3}: \ \max\big\{(\mathbf d- \mathbf B_1 \mathbf x^*- \mathbf {E}_c\mathbf{u}^*_{fc} - \mathbf {E}_d\mathbf{u}^*_{fd})^\intercal\pi: \pi\in \Pi\big\}
\end{eqnarray}
Note that  $\mathbf{SP3}$ is unbounded for $(\mathbf x^*, \mathbf u^*_{fc}, \bf u^*_{fd})$. Solving it by a linear program solver will yield an extreme ray in $\mathcal{R}_{\Pi}$, denoted by $\gamma^*$, through which $\mathbf{SP3}$ becomes infinity. Hence, by convention, we set $\eta_o(\mathbf x^*)$ to $+\infty$ in this case.

With subproblems defined, we next present detailed operations of the modified parametric C\&CG to handle mixed integer DDU, which is hereafter referred to as \textit{parametric C\&CG-MIU}. Let $\widehat{(\bf u_{od},\pi)}$ and $\widehat{(\bf u_{fd},\gamma)}$ store $(\bf u^*_{od},\pi^*)$ and $(\bf u^*_{fd},\gamma^*)$ obtained from solving subproblems in all previous iterations. By Corollary \ref{cor_Pi_partial_KKTOmegaMIPDDU}, computing the master problem, which is built upon $\widehat{(\bf u_{od},\pi)}$ and $\widehat{(\bf u_{fd},\gamma)}$,  yields lower bound $\textit{LB}$ to \eqref{eq_2RO}. Also, a feasible $\mathbf x$, together with $\eta_o(\mathbf x)$, provides upper bound $\textit{UB}$. With $t$ being the iteration counter, this algorithm iteratively refines those bounds until  the optimality tolerance $T
\!O\!L$ is achieved.  
	
	\noindent\underline{\textbf{Parametric C\&CG-MIU}}
\begin{description}
	\item[Step 1] Set $LB = -\infty$, $UB=+\infty$, $t=1$, and set 
	 $\widehat{(\bf u_{od},\pi)}$ and $\widehat{(\bf u_{fd},\gamma)}$ by an initialization strategy. 

	\item[Step 2] Solve the following master problem.
	\begin{equation*}
		\label{eq_master_1}
		\begin{split}
			\mathbf{MP}: \underline w=\min \  &  \mathbf{c}_1\mathbf x+ \eta \\
			\mathrm{s.t.} \  & \mathbf x\in \mathcal{X}\\
			&  \Big\{\eta\geq  \mathbf{c_2}\mathbf y^{(\mathbf u_d,\pi)}-\dot{u}(\tilde{\bf u}^{(\mathbf u_d,\pi})), \ \mathbf y^{(\mathbf u_d, \pi)}\geq \mathbf 0\\
			& \ \ \mathbf{B}_2\mathbf y^{(\mathbf u_d,\pi)}+\mathbf{E}_c\mathbf{u}_c^{(\mathbf u_d,\pi)}+\dot{u}(\tilde{\bf u}^{(\mathbf u_d, \pi)})\geq \mathbf{d-\mathbf B}_1\mathbf x-\mathbf{E}_d\mathbf{u}_d\\
			& \ \ (\mathbf{u}^{(\mathbf u_d,\pi)}_c, \tilde{\bf u}^{(\mathbf u_d,\pi)}, \cdot)\in
			\mathcal{OU}(\mathbf x,\mathbf u_d, \pi) \Big\} \ \forall (\mathbf u_d,\pi)\in \widehat{(\bf u_{od},\pi)}\\
			&\Big\{\mathbf{B}_2\mathbf y^{(\mathbf u_d,\gamma)}+\mathbf{E}_c\mathbf{v}_c^{(\bf u_d,\gamma)}+\dot{u}(\tilde{\mathbf v}^{(\mathbf u_d,\gamma)})\geq \mathbf{d-\mathbf B}_1\mathbf x -\mathbf{E}_d\mathbf u_d\\
			& \ \ \mathbf y^{(\mathbf u_d,\gamma)}\geq \mathbf 0, \  (\mathbf{v}_c^{(\bf u_d\gamma)}, \tilde{\mathbf v}^{(\mathbf u_d,\gamma)}, \cdot)\in \mathcal{OV}(\mathbf x, \mathbf u_d, \gamma) \Big\}  \ \forall (\mathbf u_d,\gamma)\in \widehat{(\bf u_{fd},\gamma)}
		\end{split}
	\end{equation*}

	If it is infeasible, report infeasibility of $\mathbf{2-Stg \ RO}$ in \eqref{eq_2RO} and terminate. Otherwise, derive its optimal solution $(\mathbf x^*, \eta^*, \cdot)$, and its optimal value $\underline w$. Update $LB= \underline w$.
	\item[Step 3] Solve subproblem  $\mathbf{SP1}$ in \eqref{eq_SP1} and derive optimal $(\mathbf u^*_{fc}, \mathbf u^*_{fd})$ and $\eta_f(\mathbf x^*)$.
	\item[Step 4] Cases based on $\eta_f(\mathbf x^*)$
	\begin{description}
		\item[(Case A): $\eta_f(\mathbf x^*)=0$] \textrm{}\\ $(i)$ compute $\mathbf{SP2}$ in \eqref{eq_SP2} to derive $\eta_o(\mathbf x^*)$, $(\mathbf u^*_{oc}, \mathbf u^*_{od})$ and corresponding $\pi^*$; $(ii)$ update $\widehat{(\bf u_{od},\pi)}=\widehat{(\bf u_{od},\pi)}\cup\{(\mathbf u^*_{od}, \pi^*)\}$; $(iii)$ augment master problem $\mathbf{MP}$ accordingly, i.e., create variables $(\mathbf u_c^{(\mathbf u^*_{od},\pi^*)}$, $\tilde{\bf u}^{(\mathbf u^*_{od},\pi^*)}, \cdot)$ and $\mathbf y^{(\mathbf u^*_{od},\pi^*)}$, and add the following constraints (referred to as \textit{the optimality cutting set}) to $\mathbf{MP}$.
		\begin{equation}
			\label{eq_CCG_optimality}
			\begin{split}
				& \eta\geq \mathbf{c_2}\mathbf y^{(\mathbf u^*_{od}, \pi^*)}  -\dot{u}(\tilde{\bf u}^{(\mathbf u^*_{od},\pi^*)}) \\
				& \mathbf{B}_2\mathbf y^{(\mathbf u^*_{od},\pi^*)}+\mathbf{E}_c\mathbf{u}_c^{(\mathbf u^*_{od},\pi^*)}+\dot{u}(\tilde{\bf u}^{(\mathbf u^*_{od},\pi^*)})\geq \mathbf{d-\mathbf B}_1\mathbf x-\mathbf{E}_d\mathbf{u}^*_{od}\\
				& (\mathbf{u}_c^{(\mathbf u^*_{od},\pi^*)}, \tilde{\bf u}^{(\mathbf u^*_{od},\pi^*)}, \cdot)\in
				\mathcal{OU}(\mathbf x, \mathbf u^*_{od}, \pi^*), \ \mathbf y^{(\mathbf u^*_{od},\pi^*)} \geq \mathbf 0
			\end{split}
		\end{equation}
			
		\item[(Case B): $\eta_f(\mathbf x^*) >0$] \textrm{}\\
		$(i)$ compute  $\mathbf{SP3}$ in \eqref{eq_SP3} to derive an extreme ray $\boldsymbol\gamma^*$ of $\Pi$, and set $\eta_o(\mathbf x^*)=+\infty$; $(ii)$ update $\widehat{(\bf u_{fd},\gamma)}=\widehat{(\bf u_{fd},\gamma)}\cup\{(\mathbf u^*_{fd},\gamma^*)\}$; $(iii)$ augment master problem $\mathbf{MP}$ accordingly, i.e., 
		create variables $(\mathbf v_c^{(\mathbf u^*_{fd},\gamma^*)}$, $\tilde{\bf v}^{(\mathbf u^*_{fd},\gamma^*)}, \cdot)$ and $\mathbf y^{(\mathbf u^*_{fd},\gamma^*)}$, and add the following constraints (referred to as \textit{the feasibility cutting set}) to $\mathbf{MP}$.	
		\begin{equation}
			\label{eq_CCG_feasibility}
			\begin{split}
				& \mathbf{B}_2\mathbf y^{(\mathbf v^*_{fd},\gamma^*)}+\mathbf{E}_c\mathbf{v}_c^{(\mathbf u^*_{fd},\gamma^*)}+\dot{u}(\tilde{\bf v}^{(\mathbf u^*_{fd},\gamma^*)})\geq \mathbf{d-\mathbf B}_1\mathbf x- \mathbf{E}_d\mathbf{u}^*_{fd}\\
				& (\mathbf{v}_c^{(\mathbf u^*_{fd},\gamma^*)}, \tilde{\bf v}^{(\mathbf u^*_{fd},\gamma^*)}, \cdot)\in
				\mathcal{OV}(\mathbf x, \mathbf u^*_{fd}, \gamma^*), \ \mathbf y^{(\mathbf u^*_{fd},\gamma^*)} \geq \mathbf 0
			\end{split}
		\end{equation}
				
	\end{description}
	\item[Step 5] Update $UB=\min\{UB, \mathbf{c}_1\mathbf{x}^* +\eta_o(\mathbf x^*)\}.$
	\item[Step 6] If $UB-LB\leq T\!O\!L$, return $\mathbf{x}^*$ and 
	terminate. Otherwise, let $t=t+1$ and go to \textbf{Step~2}.~\hfill $\square$
\end{description}
	
\begin{rem}
	$(i)$ In \textbf{Step 1}, sets $\widehat{(\bf u_{od},\pi)}$ and $\widehat{(\bf u_{fd},\gamma)}$ can be initialized by different strategies. The naive  one is to let them simply be $\emptyset$ and then populate them in the following operations. Another simple initialization strategy is to solve relaxations  in \eqref{eq_MIP_RO} or \eqref{eq_2RO_U_split}, and use the corresponding $\bf u_d$ and $\pi$ (obtained from recomputing the associated recourse problem if needed) as the first components of  $\widehat{(\bf u_{od},\pi)}$. Note also that we can set $LB$ to $w_R$ or $w_{(\bf u_d,\bf u_c)}$. \\ 	
	$(ii)$ Given the fact that $\mathcal{OV} (\bf x, \cdot, \cdot) \subseteq \mathcal U(\bf x)$, the feasibility cutting set in \eqref{eq_CCG_feasibility} can be implemented using the same format as that of the optimality cutting set in \eqref{eq_CCG_optimality} without causing any problem. That is, as noted in  \cite{zeng2013solving,zeng2022two}, they can be unified into a single type of cutting set in the form of \eqref{eq_CCG_optimality}. Hence, unless otherwise stated, we use $\widehat{(\bf u_{d},\pi)}$ to store components from both $\widehat{(\bf u_{od},\pi)}$ and $\widehat{(\bf u_{fd},\gamma)}$, and augment $\mathbf{MP}$ using the unified cutting sets.  \\	
$(iii)$Computationally, function $\dot{u}(\tilde{\bf u})$ can be realized by two approaches. One is to simply use linear function $\dot{u}(\tilde{\bf u})= M(\mathbf{1}^\intercal\tilde{\bf u}),$
	which is similar to the idea used in \eqref{eq_LP_parametric_ed}. Another one can be achieved  by making use of a binary variable $\theta$ and an inequality as in the following. 
	\begin{eqnarray}
		\label{eq_const_relax_functionBin}
		\begin{split}
			\dot{u}(\tilde{\bf u}) & = \Big\{M\theta: \ \theta \leq M(\mathbf{1}^\intercal\tilde{\bf u}), \ \theta \in \{0,1\}\Big\}
		\end{split}  
	\end{eqnarray}
	It can be seen from $\mathbf{MP}$ that $\dot{u}(\tilde{\bf u})$ tends to be $+\infty$, i.e., $\theta=1$. So, we only need to restrict $\theta$ to $0$ whenever $\tilde{\bf u}=\mathbf 0$, which is achieved by the inequality constraint in \eqref{eq_const_relax_functionBin}. Compared to the first approach, the second one is rather numerically more stable. \\
$(iv)$ Assuming $\cal X$ is an integer set, if DDU set $\cal U(\bf x)$ has the TU property, we can, according to Proposition \ref{prop_extremepoint}, simply replace it by its continuous relaxation and employ standard parametric C\&CG to compute. Nevertheless, we note that it generally reduces the number of iterations if integral restrictions can be additionally imposed. Let $\bf u_d$ be discrete variables of the DDU that relax to continuous ones. Specifically, for $\cal{OU}$ and $\cal{OV}$, in addition to constraints for optimality conditions, we impose $\bf u_d\in \mathbb{Z}^{m_u}_+$ to enforce integral optimal solutions.  
\end{rem}

\subsection{Convergence and Complexity}
To analyze parametric C\&CG-MIU on a consistent basis, without loss of generality, we assume that $T\!O\!L=0$, and set $\widehat{(\bf u_{d},\pi)}$ (including both $\widehat{(\bf u_{fd},\pi)}$ and $\widehat{(\bf u_{fd},\gamma)}$) is initialized by the naive strategy, i.e., it is an empty set at the beginning of parametric C\&CG-MIU.
	
\begin{thm}
	\label{thm_mipddu_conv}
	 Parametric C\&CG-MIU will not repeatedly generate any $(\mathbf u^*_{fd},\gamma^*)$ or $(\mathbf u^*_{od},\pi^*)$ unless it terminates. Upon termination, it either reports that $\mathbf{2-Stg \ RO}$ in (\ref{eq_2RO}) is infeasible, or converges to its optimal value with an exact solution.
\end{thm}
\begin{proof}
	Note that, according to Corollaries \ref{cor_Pi_equivMIPDDU} and \ref{cor_Pi_partial_KKTOmegaMIPDDU}, it is clear that if $\mathbf{MP}$ is infeasible, so is $\mathbf{2-Stg \ RO}$ in (\ref{eq_2RO}). 
	
	\noindent Claim 1:  Assume that $\mathbf x^*$, obtained from computing $\mathbf{MP}$, is infeasible to $\mathbf{2-Stg \ RO}$, and $(\mathbf u^*_{fd},\gamma^*)$ is derived after solving $\mathbf{SP1}$ and $\mathbf{SP3}$ in the current iteration. Then, $(\mathbf u^*_{fd},\gamma^*)$ will not be derived in any following iterations.
	\begin{proof}[Proof of Claim 1:]

	Without loss of generality, we assume $\mathbf x^\circ$ is the first-stage solution obtained from solving  $\mathbf{MP}$ in some following iteration. Note that a feasibility cutting set \eqref{eq_CCG_feasibility} defined on $(\mathbf u^*_{fd},\gamma^*)$ has been added to $\mathbf{MP}$. If $\mathcal U(\mathbf x^\circ|\mathbf u^*_{fd})=\emptyset$, then $\mathbf u^*_{fd}$ will not appear in the solution of any subproblem in this iteration, which concludes our claim. If $\mathcal U(\mathbf x^\circ|\mathbf u^*_{fd})\neq\emptyset$,	then by Lemma~\ref{lem_FeasDisU} we have $\dot{u}(\tilde{\mathbf v}^{(\mathbf u^*_{fd},\gamma^*)})=0$ in this feasibility cutting set. Hence, the following problem is guaranteed to be feasible.
	\begin{align*}
		\min&\big\{0: \ \bf y^{(\mathbf u^*_{fd},\gamma^*)}\in \mathcal{Y}(\bf x^\circ,\mathbf{v}_c^{(\mathbf u^*_{fd},\gamma^*)},\bf u^*_{fd}),\ (\mathbf{v}_c^{(\mathbf u^*_{fd},\gamma^*)}, \cdot)\in
		\mathcal{OV}(\mathbf x^\circ, \mathbf u^*_{fd}, \gamma^*)\big\}\\
			&=\min_{(\bf v_c^{(\mathbf u^*_{fd}, \gamma^*)},\cdot)\in \mathcal{OV}(\mathbf x^\circ, \bf u^*_{fd}, \gamma^*)}\min\big\{0:  \bf y^{(\mathbf u^*_{fd},\gamma^*)}\in \mathcal{Y}(\bf x^\circ, \mathbf{v}_c^{(\mathbf u^*_{fd},\bf 0, \gamma^*)}, \bf u^*_{fd})\big\}
	\end{align*}
   By considering the dual of the inner minimization problem, we have
    \begin{align*}
		0&\geq\min_{(\bf v_c^{(\mathbf u^*_{fd},\gamma^*)},\cdot)\in \mathcal{OV}(\mathbf x^\circ, \mathbf u^*_{fd}, \gamma^*)}\max  \big\{(\mathbf d- \mathbf B_1 \mathbf x^\circ- \mathbf {E}_c\mathbf{v}_{c}^{(\mathbf u^*_{fd},\gamma^*)} - \mathbf {E}_d\mathbf{u}^*_{fd})^\intercal\gamma: \gamma\in \mathcal{P}_{\Pi}\cup\mathcal{R}_{\Pi}\big\}\\
		&\geq\min_{(\bf v_c^{(\mathbf u^*_{fd}, \gamma^*)},\cdot)\in \mathcal{OV}(\mathbf x^\circ, \mathbf u^*_{fd}, \gamma^*)}\big\{(\mathbf d- \mathbf B_1 \mathbf x^\circ- \mathbf {E}_c\mathbf{v}_{c}^{(\mathbf u^*_{fd},\gamma^*)} - \mathbf {E}_d\mathbf{u}^*_{fd})^\intercal\gamma^*\big\}\\
		&=\max_{(\bf v_c^{(\mathbf u^*_{fd}, \gamma^*)},\cdot)\in \mathcal{OV}(\mathbf x^\circ, \mathbf u^*_{fd}, \gamma^*)}\big\{(\mathbf d- \mathbf B_1 \mathbf x^\circ- \mathbf {E}_c\mathbf{v}_{c}^{(\mathbf u^*_{fd},\gamma^*)} - \mathbf {E}_d\mathbf{u}^*_{fd})^\intercal\gamma^*\big\},
	\end{align*}
	where the last equality holds by the definition of $\mathcal{OV}$, i.e., 
	 it is the optimal solution set of $\mathbf{LP}(\mathbf x^\circ, \mathbf u^*_{fd}, \gamma^*)$. Hence, $(\mathbf u^*_{fd},\gamma^*)$ will not be derived by $\textbf{SP3}$.
	\end{proof}

Next, we consider the other case where $\mathbf{MP}$ generates feasible $\bf x^*$.
		
\noindent Claim 2: Assume that  $\mathbf x^*$ generated by $\mathbf{MP}$ is feasible, and $(\mathbf u^*_{od},\pi^*)$ is derived as  an optimal solution to $\mathbf{SP2}$ in the current iteration. If $(\mathbf u^*_{od},\pi^*)$ has appeared as an optimal solution to $\mathbf{SP2}$ in some previous iteration, we have $UB=LB$.
\begin{proof}[Proof of Claim 2:]
    It is obvious that $\mathcal U(\mathbf x^*|\mathbf u^*_{od})$ is non-empty, and $\dot{u}(\tilde{\bf u}^{(\mathbf u^*_{od}, \pi^*)})=0$ in the optimality cutting set (\ref{eq_CCG_optimality}) defined on $(\mathbf u^*_{od},\pi^*)$, which has been a part of $\mathbf{MP}$ after the first time $(\mathbf u^*_{od},\pi^*)$ is generated.  Let $\eta^*$ and $\mathbf y^{*(\bf u^*_{od}, \pi^*)}$ denote the optimal values for $\eta$ and $\mathbf y^{(\bf u^*_{od}, \pi^*)}$, respectively, of the current $\mathbf{MP}$. We have
		\begin{align*}
			\eta^*\geq&\mathbf{c_2}\mathbf y^{*(\bf u^*_{od}, \pi^*)}\\
			\geq&\min\Big\{\mathbf{c_2}\mathbf y^{(\bf u^*_{od}, \pi^*)}: \mathbf{B}_2\mathbf y^{(\bf u^*_{od}, \pi^*)}\geq \mathbf{d-\mathbf B}_1\mathbf x^*-\mathbf{E}_c\mathbf{u}_c^{(\bf u^*_{od}, \pi^*)}- \mathbf{E}_d\mathbf{u}^*_{od},\\
			& \ \ \ \ \ \  (\mathbf{u}_c^{(\bf u^*_{od}, \pi^*)},\cdot)\in\mathcal{OU}(\mathbf x^*, \mathbf u^*_{od}, \pi^*), \ \mathbf y^{(\bf u^*_{od}, \pi^*)} \geq \mathbf 0\Big\}\\
		   = &\min_{(\mathbf{u}_c^{(\bf u^*_{od}, \pi^*)}, \cdot)\in\mathcal{OU}(\mathbf x^*, \mathbf u^*_{od}, \pi^*)} \max_{\pi\in \cal P_{\Pi}}\Big\{(\mathbf{d-\mathbf B}_1\mathbf x^*-\mathbf{E}_c\mathbf{u}_c^{(\bf u^*_{od}, \pi^*)}- \mathbf{E}_d\mathbf{u}^*_{od})^\intercal\pi\Big \}\\
		   \geq &\min_{(\mathbf{u}_c^{(\bf u^*_{od}, \pi^*)}, \cdot)\in\mathcal{OU}(\mathbf x^*, \mathbf u^*_{od}, \pi^*)} (\mathbf{d-\mathbf B}_1\mathbf x^*-\mathbf{E}_c\mathbf{u}_c^{(\bf u^*_{od}, \pi^*)}- \mathbf{E}_d\mathbf{u}^*_{od})^\intercal\pi^*\\
		   = &\max_{(\mathbf{u}_c^{(\bf u^*_{od}, \pi^*)}, \cdot)\in\mathcal{OU}(\mathbf x^*, \mathbf u^*_{od}, \pi^*)} (\mathbf{d-\mathbf B}_1\mathbf x^*-\mathbf{E}_c\mathbf{u}_c^{(\bf u^*_{od}, \pi^*)}- \mathbf{E}_d\mathbf{u}^*_{od})^\intercal\pi^*\\
			=&\max\Big\{(\mathbf{d-\mathbf B}_1\mathbf x^*-\mathbf{E}_c\mathbf{u}_c-\mathbf{E}_d\mathbf{u}^*_{od})^\intercal\pi^*: \ \mathbf F_c(\mathbf x^*)\mathbf u_c\leq\mathbf h+\mathbf{Gx}-\mathbf F_d(\mathbf x^*)\mathbf{u}^*_{od},\ \mathbf u_c\geq \mathbf 0\Big\}\\
			= & \max_{\bf u_d=\bf u^*_{od}, \bf u_c\in \cal U(\bf x^*|\bf u^*_{od})}\max_{\pi\in \Pi}\Big\{(\mathbf{d-\mathbf B}_1\mathbf x^*-\mathbf{E}_c\mathbf{u}_c-\mathbf{E}_d\mathbf{u}^*_{od})^\intercal\pi\Big\}\\
			= & \eta_o(\mathbf x^*)
		\end{align*}
	    The first equality follows the duality of recourse problem, the  second and third equalities hold by the definition of  $\mathcal{OU}(\mathbf x^*, \mathbf u^*_{od}, \pi^*)$, and the last two equalities follow from the claim statement that $(\bf u^*_{od}, \pi^*)$ is an optimal solution of $\mathbf{SP2}$ (with the dual of recourse). Hence, we have $LB=\bf c_1\bf x^*+\eta^*\geq \bf c_1\bf x^*+\eta_o(\bf x^*)\geq UB$, which leads to the desired conclusion. 
\end{proof}
Given that both sets $\mathbb U_d$ and $\mathcal P_{\Pi}\cup\mathcal R_{\Pi}$ are fixed and finite, it is straightforward to conclude that parametric C\&CG-MIU algorithm will terminate by either reporting the infeasibility of $\mathbf{2-Stg \ RO}$ or  converging to an exact solution after a finite number of iterations.
\end{proof}
	
According to the proof of Theorem \ref{thm_mipddu_conv}, we can easily bound the number of iterations.
	
\begin{cor}
	\label{cor_mipddu_comp1}
	The number of iterations of parametric C\&CG-MIU  before termination is bounded by $|\mathbb U_d|\times(|\mathcal P_\Pi|+|\mathcal R_\Pi|)$. Hence, the algorithm is of $O(|\mathbb U_d|\times\binom{n_y+\mu_y}{\mu_y})$ iteration complexity, where $\mu_y$ denotes the number of rows of matrix $\mathbf B_2$ in $\mathcal Y(\mathbf x,\mathbf u)$. $\hfill\square$
\end{cor}
	
The next result follows easily from structures of feasibility and optimality cutting sets and derivations of lower and upper bounds.  
	
\begin{prop}
	\label{prop_DupInf}
	$(i)$ If solution $\mathbf x^*$ of $\mathbf{MP}$ is infeasible according to $\mathbf{SP1}$ and $\mathbf{SP3}$ and a corresponding feasibility cutting set is included in $\mathbf{MP}$, then it will not be generated by $\mathbf{MP}$ in the following iterations. \\
	$(ii)$ If $\mathbf x^*$ is an optimal solution to $\mathbf{MP}$ in both iterations $t_1$ and $t_2$ with $t_1<t_2$, then the parametric C\&CG-MIU algorithm terminates and $\mathbf x^*$ is optimal to $\mathbf{2-Stg \ RO}$. $\hfill\square$
\end{prop}

Consequently, we have the following bound on the iteration complexity.

\begin{cor}
	\label{cor_mipddu_comp2}
	If set $\mathcal X$ is discrete and of a finite number of elements, then the number of iterations of parametric C\&CG-MIU before termination is bounded by $|\mathcal X|$, i.e., the algorithm is of $O(|\mathcal X|)$ iteration complexity.
\end{cor}

In \cite{zeng2022two}, a deeper understanding on the iteration complexity of parametric C\&CG is  developed, which is based on the concept of basis of polytope $\mathcal U(\bf x)$. It is particularly useful when decision dependence appears in RHS of $\mathcal U(\bf x)$, where the iteration complexity generalizes that of basic C\&CG developed for RO with DIU. Next, we present a result of \cite{zeng2022two} modified according to the structure of mixed integer DDU to support our analysis. 

\begin{lem}
	\label{lem_feasible_then_optimal}
	[Adapted from Lemma 22 in \cite{zeng2022two}] 	Assume that $\mathcal{U}(\bf x)$ is with RHS decision-dependence. Consider $\mathbf{LP}(\mathbf x^0, \bf u^*_d, \boldsymbol\beta)$ for fixed $(\bf u^*_d,\boldsymbol\beta)$, and basis $\mathbb{B}^0$ with $\tilde{\bf u}\notin \mathbb{B}^0$ that is an optimal basis, i.e., its basic solution (BS) with respect to $\tilde{\mathcal{U}}(\mathbf x^0|\mathbf u^*_{d})$ is feasible and optimal. 
	If $\mathbb{B}^0$'s BS with respect to $\tilde{\mathcal{U}}(\mathbf x^1|\bf u^*_d)$	is feasible, i.e., a basic feasible solution, it is also optimal to $\mathbf{LP}(\mathbf x^1, \bf u^*_d, \boldsymbol\beta)$.  Moreover, if $\mathbb{B}^0$ yields the unique optimal solution to $\mathbf{LP}(\mathbf x^0, \bf u^*_d, \boldsymbol\beta)$, it also yields the unique one to $\mathbf{LP}(\mathbf x^1, \bf u^*_d, \boldsymbol\beta)$. $\hfill\square$
\end{lem}

Following this lemma, we study the consequence of repeated bases. Without loss of generality, we assume that the unified cutting sets in the form of (\ref{eq_CCG_optimality}) are supplied to $\mathbf{MP}$ to simplify our exposition.

\begin{thm}
	\label{thm_complex_2^n} 
	Let $\hat{\bb U}_d^t$ be the set of all discrete $\mathbf u_{od}$ and $\bf u_{fd}$ appeared in $\mathbf{MP}$ at the end of iteration $t$, and for each $\mathbf u_d\in\hat{\bb U}_d^t$, let $\mathfrak B^t_{\mathbf u_d}$ be the collection of bases obtained from the sets $\mathcal{OU}(\mathbf x,\mathbf u_d,\cdot)$ and $\mathcal{OV}(\mathbf x,\mathbf u_d,\cdot)$ at the end of iteration $t$. For any iterations $t_1<t_2$, if $\hat{\bb U}_d^{t_1}=\hat{\bb U}_d^{t_2}$, and $\forall\mathbf u_d\in\hat {\bb U}_d^{t_1}$, $\mathfrak B^{t_1}_{\mathbf u_d}=\mathfrak B^{t_2}_{\mathbf u_d}$, parametric C\&CG-MIU  terminates at $t_2$, and $\mathbf x^{t_1}$, an optimal solution to $\mathbf{MP}$ in iteration $t_1$, is optimal to $\mathbf{2-Stg \ RO}$. 
\end{thm}
\begin{proof}
	Let $(\mathbf x^t, \eta^t, \cdot)$ denote the optimal solution of $\mathbf{MP}$ in iteration $t$. Since $\hat {\bb U}_d^{t_1}=\hat {\bb U}_d^{t_2}$, there exists $\mathbf u_d^*\in\hat{\bb U}_d^{t_1}$ such that $\mathcal U(\mathbf x^{t_2-1}|\bf u^*_d)\neq\emptyset$. Note otherwise that a new $\bf u_d$ would be generated from subproblems in iteration $t_2-1$, which leads to $\hat {\bb U}_d^{t_1}\neq \hat {\bb U}_d^{t_2}$. In addition, starting from iteration $t_2$, for all $\mathbf x\in\mathcal X$, there exists $\mathbf u_d^*\in\hat {\bb U}_d^{t_1}$ such that $\mathcal U(\mathbf x|\bf u^*_{d})\neq\emptyset$. Otherwise, such $\mathbf x$ renders the objective value of $\mathbf{MP}$ to $-\infty$ and $\hat {\bb U}_d^{t_1}=\hat {\bb U}_d^{t_2}$ does not hold. 
	
	As $\mathbf{MP}$ is finite in iteration $t_2$, it is sufficient to examine any $\mathbf u^*_d\in\hat {\bb U}_d^{t_2}$ whose associated $\tilde{\bf u}^{\bf u^*_d,\pi}=0$ for some $\pi$.
	Hence,  suppose $(\mathbf u^*_d, \pi^*)\in \widehat{(\bf u_{d},\pi)}$ whose associated $\tilde{\bf u}^{\bf u^*_d,\pi^*}=0$, and $\mathbb B_{\mathbf u_d^*}\in\mathfrak B^{t_1}_{\mathbf u_d^*}=\mathfrak B^{t_2}_{\mathbf u_d^*}$ yields an optimal solution belonging to set $\mathcal{OU}(\mathbf x^{t_2},\mathbf u_d^*,\pi^*)$. By Lemma~\ref{lem_feasible_then_optimal}, given that $\mathbb B_{\mathbf u_d^*}$ is feasible to $\mathcal U(\mathbf x^{t_1}|\mathbf u^*_{d})$, 
	 it defines an optimal solution in $\mathcal{OU}(\mathbf x^{t_1},\mathbf u_d^*,\pi^*)$. It follows that $(\mathbf x^{t_1},\eta^{t_1}, \cdot)$ is feasible (and optimal) to $\mathbf{MP}$ in $t_2$,  which, according to Proposition~\ref{prop_DupInf}.$(ii)$, terminates the algorithm with $\bf x^{t_1}$ being an optimal solution.
\end{proof}

Because all possible combination of bases is finite, the next result follows.
\begin{cor}
	The number of iterations of  parametric C\&CG-MIU before termination is bounded by $|\mathbb U_d|\times2^{\binom{n_u+\mu_u}{\mu_u}}$, i.e., the algorithm is of $O(|\mathbb U_d|\times2^{\binom{n_u+\mu_u}{\mu_u}})$ iteration complexity, where $\mu_u$ is the number of constraints in set $\mathcal U(\mathbf x)$.
\end{cor}

Result in Theorem \ref{thm_complex_2^n} actually can be significantly improved if $\mathbf{LP}(\mathbf x, \mathbf u_d, \boldsymbol\beta)$ always has a unique optimal solution for any $(\bf u_d,\beta)$ with $\bf u_d\in\bb{U}_d$ and $\beta\in \cal{P}_{\Pi}\cup\cal{R}_{
\Pi}$, which is referred to as the unique optimal solution (or ``uniqueness'' for short) property. 

\begin{thm}
	\label{thm_complex_MIU_basis}
	Assume that $\mathbf{LP}(\mathbf x, \mathbf u_d, \boldsymbol\beta)$ has the uniqueness property, and let $\mathbb B_{\mathbf u_d^t}$ be the basis of $\mathcal U(\mathbf x|\bf u^t_d)$	associated with optimal solution $(\mathbf u_c^t,\mathbf u_d^t)$ of $\mathbf{SP1}$ or $\mathbf{SP2}$ in iteration $t$. If 2-tuple $(\mathbf u_d^t,\mathbb B_{\mathbf u_d^t})$ has appeared in some previous iteration, parametric C\&CG-MIU terminates.
\end{thm}

\begin{proof}
	Assume the repeat happens in iterations $t_1$ and $t_2$ with $t_1<t_2$, and let $\mathbf x^{t_1}$, $\mathbf x^{t_2}$ be the first stage solutions of $\mathbf {MP}$, $(\mathbf u_c^{t_1},\mathbf u_d^{t})$, $(\mathbf u_c^{t_2},\mathbf u_d^{t})$ be the worst case scenarios from solving either $\mathbf{SP1}$ or $\mathbf{SP2}$ subproblem, and $\pi^{t_1}$ and $\pi^{t_2}$ be the corresponding extreme points or rays of $\Pi$ obtained in those iterations, respectively. Then $\mathcal U(\mathbf x^{t_1}|\mathbf u_d^t)\neq\emptyset$, $\mathcal U(\mathbf x^{t_2}|\mathbf u_d^t)\neq\emptyset$, and $\mathbb B_{\mathbf u_d^{t}}$ defines a basic feasible solution to both of them.
	
	Since $\mathbb{B}_{\mathbf u_d^{t}}$ defines the optimal solution $(\mathbf u_c^{t_1},\mathbf 0)$ of $\mathbf{LP}(\mathbf x^{t_1},\mathbf{u}^t_d, \pi^{t_1})$, by Lemma \ref{lem_feasible_then_optimal}, it also defines the optimal solution of $\mathbf{LP}(\mathbf x^{t_2}, \mathbf u^t_d, \pi^{t_1})$. 
	It is $(\mathbf u_c^{t_2},\mathbf 0)$, according to the assumption that $\mathbb B_{\mathbf u_d^{t}}$ defines the optimal solution of $\mathbf{LP}(\mathbf x^{t_2}, \mathbf u^t_d, \pi^{t_2})$. Note that a cutting set with $\mathcal{OU}(\mathbf x, \mathbf u_d^t, \pi^{t_1})$ or $\mathcal{OV}(\mathbf x, \mathbf u_d^t, \pi^{t_1})$ has been added to $\mathbf{MP}$ in iteration $t_1$. Hence we have $\mathcal{OU}(\mathbf x^{t_2}, \mathbf u_d^t, \pi^{t_1})=\{(\mathbf u_c^{t_2},\mathbf 0,\cdot)\}$ or $\mathcal{OV}(\mathbf x^{t_2}, \mathbf u_d^t, \pi^{t_1})=\{(\mathbf u_c^{t_2},\mathbf 0,\cdot)\}$. If the unified cutting sets are adopted, we have
	\begin{align*}
		LB=\mathbf c_1\mathbf x^{t_2}+\eta\geq\mathbf c_1\mathbf x^{t_2}+\min\{\mathbf c_2\mathbf y:\mathbf y\in\mathcal Y(\mathbf x^{t_2},\mathbf u_c^{t_2},\mathbf u_d^{t})\}=\mathbf c_1\mathbf x^{t_2}+\eta_o(\mathbf x^{t_2})=UB. \quad\pushQED{\qed}\qedhere
	\end{align*}
\end{proof}
Actually, as shown in \cite{zeng2022two}, it is without loss of generality to modify any  $\bsm\beta$ slightly to ensure a unique optimal solution in the execution of parametric C\&CG-MIU.   
\begin{cor}
	If $\mathbf{LP}(\mathbf x, \mathbf u_d, \boldsymbol\beta)$ has the uniqueness property, then the number of iterations required for parametric C\&CG-MIU before termination is bounded by $|\mathbb U_d|\times\binom{n_u+\mu_u}{\mu_u}$. Therefore, the algorithm is of $O\left(\min\Big\{|\mathbb U_d|\times\binom{n_y+\mu_y}{\mu_y},|\mathcal X|,|\mathbb U_d|\times\binom{n_u+\mu_u}{\mu_u}\Big\}\right)$ iteration complexity.
\end{cor}

\section{Two-Stage RO DDU with Mixed Integer Recourse}
\label{sect_MIPrecourse}
In  \cite{zhao2012exact}, basic C\&CG is extended in a nested fashion to exactly compute DIU-based  $\mathbf{2-Stg \ RO}$ with an MIP recourse, which then has been adopted to solve many practical problems (e.g., \cite{danandeh2014job,fang2019adaptive,garcia2021robust}).   
To handle MIP recourse in the context of DDU-based $\mathbf{2-Stg \ RO}$, we extend parametric C\&CG in a similar approach. To help with the exposition of rather complex algorithm structure,  we concentrate on polytope DDU set in this section to provide detailed descriptions and theoretical analysis. We then present  in the next section additional modifications on algorithm operations for the case of mixed integer DDU.
Note that a new inner subroutine is introduced, which eliminates the assumption of relatively complete recourse property made in \cite{zhao2012exact} for the MIP recourse problem.    

\subsection{Equivalent Reformulations by Enumerations}
As we focus on polytope DDU set where $\bf u$ is continuous, $\bf u$ and $\bf u_c$ are interchangeable in this section, and so do $\bf E_c$ and $\bf E$.   

\subsubsection{The First Equivalent Reformulation}
Consider the mixed integer recourse problem of $\mathbf{2-Stg \ RO}$ in \eqref{eq_2RO}. Similar to \cite{zhao2012exact}, we treat discrete and continuous recourse variables separately as in the following.  Let $\mathcal Y(\mathbf x, \mathbf u, \mathbf y_d)\equiv\{\mathbf B_{2,c}\mathbf y_c \geq \mathbf d-\mathbf B_1\mathbf x-\mathbf E\mathbf u-\mathbf B_{2,d}\mathbf y_d, \ \mathbf y_c\geq \mathbf 0\}$, and, by a slight abuse of notation, $\Pi=\{\mathbf B_{2,c}^\intercal\pi\leq \mathbf{c}^{\intercal}_{2,c},\ \pi\geq\mathbf 0\}$, which is not empty according to Assumption $(\textit{A3})$. We have
\begin{subequations}
 \begin{align}
	&\min_{\mathbf y_c,\mathbf y_d} \Big\{\mathbf c_{2,c}\mathbf y_c+\mathbf c_{2,d}\mathbf y_d: \mathbf y_d\in \mathbb Y_d, \mathbf y_c\in \mathcal Y(\mathbf x, \mathbf u, \mathbf y_d)\Big\}\\
	=& \min_{\mathbf y_d\in \mathbb Y_d} \Big\{\mathbf c_{2,d}\mathbf y_d+\min_{\mathbf y_c\in \mathcal Y(\mathbf x, \mathbf u, \mathbf y_d)}\{\mathbf c_{2,c}\mathbf y_c\} \Big\} \label{eq_max_min_min}\\
	=& \min_{\mathbf y_d\in \mathbb Y_d} \Big\{\mathbf c_{2,d}\mathbf y_d+\max_{\pi\in \Pi}\{(\mathbf d-\mathbf B_1\mathbf x-\mathbf E\mathbf u-\mathbf B_{2,d}\mathbf y_d)^\intercal\pi\} \Big\}. \label{eq_recourse-min-max}
\end{align}
\end{subequations}
 \begin{rem}
 	\label{rem_MIP_recourse_extreme_pointray}
Note again that $\Pi$ is a fixed polyhedron independent of any decision variable.  Following our convention, we let $\mathcal{P}_{\Pi}$ and $\mathcal{R}_{\Pi}$ being the sets of extreme points and extreme rays of $\Pi$. It is straightforward to see that, for any given $(\mathbf x^0, \mathbf u^0, \mathbf y_d^0)$, an optimal solution to the inner maximization problem of \eqref{eq_recourse-min-max} is an extreme point or ray of $\Pi$. 
 \end{rem}

 Assume $\mathbb Y_d=\{\bf y^1_d,\dots, \bf y^{|\mathbb Y_d|}_d\}$.  With \eqref{eq_recourse-min-max} and for a fixed $\mathbf x$, we now consider the max-min substructure of \eqref{eq_2RO}, which can be reformulated equivalently as the following. 
 \begin{subequations}
 \label{eq_max-min_substructure0}
  \begin{align}
   \eta_o (\mathbf x) =& \max_{\mathbf u\in \mathcal{U}(\mathbf x)}\min_{(\mathbf y_c,\mathbf y_d)\in \mathcal{Y}(\mathbf x,\mathbf u)} (\mathbf c_{2,c}\mathbf y_c+\mathbf c_{2,d}\mathbf y_d) \label{eq_max-min_substructure}\\
 	=& \max_{\mathbf u\in \mathcal{U}(\mathbf x)} \min_{\mathbf y_d\in \mathbb Y_d} \Big\{\mathbf c_{2,d}\mathbf y_d+\max_{\pi\in \Pi}\{(\mathbf d-\mathbf B_1\mathbf x-\mathbf E\mathbf u-\mathbf B_{2,d}\mathbf y_d)^\intercal\pi\} \Big\} \label{eq_max-min_DAD1}\\
 	= &\max \Big\{\hat\eta: \mathbf u\in \mathcal{U}(\mathbf x), \notag \\	&  \ \ \ \ \ \ \ \ \ \ \ \ \hat\eta\leq \mathbf{c}_{2,d}\mathbf y_d + \max_{\pi\in \Pi}\{(\mathbf d-\mathbf B_1\mathbf x-\mathbf E\mathbf u-\mathbf B_{2,d}\mathbf y_d)^\intercal\pi\}, \ \forall \mathbf y_d \in \mathbb{Y}_d \Big\} \label{eq_max-min_DAD2}\\
 	=& \max \Big\{\hat\eta:  \mathbf u\in \mathcal{U}(\mathbf x),\  \pi^t\in \cal{P}_{\Pi}\cup\cal{R}_{\Pi}, \ t=1,\dots, |\mathbb{Y}_d| \notag\\
 	&  \ \ \ \ \ \ \ \ \ \ \ \ \ \hat\eta\leq \mathbf{c}_{2,d}\mathbf y^t_d + (\mathbf d-\mathbf B_1\mathbf x-\mathbf E\mathbf u-\mathbf B_{2,d}\mathbf y^t_d)^\intercal\pi^t, \ t=1,\dots, |\mathbb{Y}_d|  \Big\}\label{eq_max-min_single_enumer}	
 \end{align}
\end{subequations}
Note that even if the minimization problem in \eqref{eq_max-min_substructure} is infeasible for some $\bf u \in \mathcal U(\bf x)$, all maximization problems in \eqref{eq_max-min_DAD2} are unbounded and hence  the  optimal value of \eqref{eq_max-min_single_enumer} goes to $\infty$, indicating the whole derivation still holds. We also mention that \eqref{eq_max-min_DAD1} has a two-stage RO's min-max-min structure that is suitable for some C\&CG method.  To facilitate our exposition, we 
introduce 
$\bsm{\pi}$, a $|\mathbb Y_d|$-tuple of extreme points or rays defined as $$\bsm{\pi}\equiv(\pi^1,\dots, \pi^{|\mathbb Y_d|}) \in  (\mathcal{P}_{\Pi}\bigcup\mathcal{R}_{\Pi})^{|\mathbb Y_d|} \equiv (\mathcal{P}_{\Pi}\bigcup\mathcal{R}_{\Pi})\times \dots\times (\mathcal{P}_{\Pi}\bigcup\mathcal{R}_{\Pi}).$$

And for a given $\bsm{\pi}$,  we define  
\begin{align}
\label{eq_eta_inner}
\begin{split}
	\hat \eta(\mathbf x, \bsm{\pi})  = \max\Big\{\hat \eta: & \mathbf u\in \mathcal{U}(\mathbf x), \\
 &\hat \eta\leq \bf{c}_{2,d}\bf y^t_d + (\mathbf d-\mathbf B_1\mathbf x-\mathbf E\mathbf u-\mathbf B_{2,d}\mathbf y^t_d)^\intercal\pi^{t}, \ t=1,\dots, |\mathbb{Y}_d|\Big\},
\end{split}
 \end{align}
which is a linear program. Then, according to the derivation in \eqref{eq_max-min_substructure0}, we have 
\begin{align}
\label{eq_max-min-to-max_allpi}
	\eta_o(\bf x) = \max_{\mathbf u\in \mathcal{U}(\mathbf x)}\min_{(\mathbf y_c,\mathbf y_d)\in \mathcal{Y}(\mathbf x,\mathbf u)} (\mathbf c_{2,c}\mathbf y_c+\mathbf c_{2,d}\mathbf y_d) =\max_{\bsm{\pi}\in (\mathcal{P}_{\Pi}\bigcup\mathcal{R}_{\Pi})^{|\mathbb Y_d|} }\{\hat \eta(\mathbf x, \bsm{\pi)}\}.
\end{align}

Next, we build a bilevel reformulation that is equivalent to the original
 $\mathbf{2-Stg \ RO}$. 
\begin{thm}
	When it has a mixed integer recourse problem, $\mathbf{2-Stg \ RO}$ in \eqref{eq_2RO} (and its other  equivalences) is equivalent to a bilevel linear optimization program as in the following.
	\begin{subequations}
		\label{eq_2stgRO_PI_MIP}
		\begin{align}
		    w^*=\min \ & \ \mathbf{c}_1\mathbf x+ \eta   \\
			\mathrm{s.t.} \ & \ \mathbf x\in \mathcal{X}\\
			&  \bigg\{ \eta\geq \hat \eta(\mathbf x, \bsm{\pi}) \bigg\}  \ \forall \bsm{\pi}\in (\mathcal{P}_{\Pi}\bigcup\mathcal{R}_{\Pi})^{|\mathbb Y_d|} \label{eq_RMIP_Pi_expo}
		\end{align}
	\end{subequations}
\end{thm}

\begin{rem}
	\label{rem_MIP_recourse_ER1}
	Obviously \eqref{eq_2stgRO_PI_MIP} is a very large-scale bilevel optimization model that is impractical to compute directly.  Constraints in \eqref{eq_RMIP_Pi_expo} clearly demonstrate a structure that is friendly to develop cutting set-based algorithms. Nevertheless, there are two fundamental technical challenges associated with $\hat\eta(\mathbf x, \bsm{\pi})$. One  is that  \eqref{eq_RMIP_Pi_expo}, which completely depend on $\bsm{ \pi}$,  are duality based cutting sets that are likely to be weak.  
	The other one is that not only the number of linear program $\hat\eta(\mathbf x, \bsm{ \pi})$'s are exponential, but also the size of such linear program is enormous, noting $|\mathbb Y_d|$ is exponential with respect to the dimension of $\mathbf y_d$. Actually, the latter issue is more critical as converting $\hat\eta(\mathbf x, \bsm{ \pi})$ into (linearized) constraints by its optimality conditions is practically infeasible. Hence, we believe that more effective  strategies should be investigated to develop  exact and efficient algorithms.  
\end{rem}

\subsubsection{The Second Equivalent Reformulation}

Rather than listing all $\bf y_d$'s in $\mathbb Y_d$ to constructing $\hat \eta(\bf x,\bsm{ \pi})$, we consider a simpler one that only involves a subset of $(\bf y_d, \pi)$ pairs. Let $\widehat{(\bf y_d,\pi)}\equiv\big\{(\bf y^1_d, \pi^1), \dots, (\bf y^{|\bb{\hat Y}_d|}, \pi^{|\bb{\hat Y}_d|})\big\}$ denote such a subset with $\bb{\hat Y}_d\equiv \{\bf y^1_d,\dots,\bf y^{|\bb{\hat Y}_d|}\}$ being the associated set of $\bf y_d$'s and  $\bsm{\hat \pi}\equiv\{\pi^1,\dots,\pi^{|\bb{\hat Y}_d|}\}$. We actually restrict $\widehat{(\bf y_d,\pi)}$ such that $\bf y_d$'s in $\bb{\hat Y}_d$  are distinct while some $\pi$'s from $\bsm{\hat \pi}$ could be same. Next, we define the following function.  
\begin{align}
\label{eq_eta_inner_subset}
\begin{split}
	 	\hat \eta\big(\mathbf x, \widehat{(\bf y_d,\pi)}\big)  = \max\Big\{\hat\eta: & \mathbf u\in \mathcal{U}(\mathbf x),\\ 
		\hat \eta\leq & \mathbf{c}_{2,d}\mathbf y_d + (\mathbf d-\mathbf B_1\mathbf x-\mathbf E\mathbf u-\mathbf B_{2,d}\mathbf y_d)^\intercal\pi, \ \forall (\bf y_d,\pi)\in \widehat{(\bf y_d,\pi)}\Big\}
 \end{split}
\end{align}
 Compared to \eqref{eq_eta_inner}, the size of \eqref{eq_eta_inner_subset} could be much smaller if $\widehat{(\bf y_d,\pi)}$ is of a small cardinality.   In the following, we make use of \eqref{eq_eta_inner_subset} to provide a strong lower bound to $\mathbf{2-Stg \ RO}$.  Indeed, we mention that the maximization  problem in \eqref{eq_eta_inner_subset} essentially plays the same role as \eqref{eq_LP_parametric_ed}, and hence we denote it by $\mathbf{LP}\big(\mathbf x, \widehat{(\bf y_d,\pi)}\big)$ in development of the new reformulation.

Basically, we can take advantage of the projection of the optimal solution set of  \eqref{eq_eta_inner_subset}, which is a linear program, onto the space hosting $\cal U(\bf x)$. Similar to \eqref{eq_2stgRO_PI_OM}, we supplement it with a replicate of recourse problem. Next, we present the complete reformulation of \eqref{eq_2RO} obtained using this idea, which lays the foundation for the development of a C\&CG type of algorithm.  
Let $\mathcal{OU}\big(\mathbf x, \widehat{(\bf y_d,\pi)}\big)$ denote this optimal solution set. We can employ the KKT conditions to characterize it as in the following, noting that the strong duality can help too. 
	\begin{align*}
	\mathcal{OU}\big(\mathbf x, \widehat{(\bf y_d,\pi)}\big) = \left\{\begin{array}{l}
		\hat \eta\leq \mathbf{c}_{2,d}^\intercal\mathbf y^t_d + (\mathbf d-\mathbf B_1\mathbf x-\mathbf E\mathbf u-\mathbf B_{2,d}\mathbf y^t_d)^\intercal\pi^{t}, t=1,\dots, |\hat{\mathbb{Y}}_d|\\
		\mathbf F(\bf x)\mathbf u\leq \mathbf h+\mathbf{Gx}\\ 
		\mathbf F^\intercal(\bf x)\boldsymbol\lambda+\sum_{t=1}^{|\hat{\mathbb{Y}}_d|}\mathbf E^\intercal\pi^t\mu^t\geq\mathbf 0 \\
		\sum_{t=1}^{|\hat{\mathbb{Y}}_d|}\mu^t=1\\
		\mu^t\circ(\mathbf{c}_{2,d}^\intercal\mathbf y^t_d + (\mathbf d-\mathbf B_1\mathbf x-\mathbf E\mathbf u-\mathbf B_{2,d}\mathbf y^t_d)^\intercal\pi^t-\hat \eta)=0, t=1,\dots,|\hat{\mathbb{Y}}_d|\\
		\boldsymbol\lambda \circ (\mathbf h+\mathbf{Gx}- \mathbf F(\bf x)\mathbf{u})=\mathbf{0} \\
		\mathbf u\circ (\mathbf F^\intercal(\bf x)\boldsymbol\lambda+ \sum_{t=1}^{|\hat{\mathbb{Y}}_d|}\mathbf E^\intercal\pi^t\mu^t) = \mathbf 0 \\
		\mathbf u\geq \mathbf{0},\  \boldsymbol\lambda\geq \mathbf{0},\ \mu^t\geq0, \ t=1,\dots, |\hat{\mathbb Y}_d|
	\end{array}\right\}
\end{align*}
 
Different from \eqref{eq_max-min_substructure0}, we next provide an $\mathcal{OU}$-based reformulation for the max-min substructure of $\mathbf{2-Stg \ RO}$ with a mixed integer recourse problem.    
\begin{prop}
\label{prop_OU_max-min_reform}
Let $2^{\mathbb{Y}_d}$ denote the power set of $\mathbb{Y}_d$. For a fixed $\mathbf x$, we let
$$\hat{\mathcal{U}}(\mathbf x)=\bigg\{\mathbf u: (\mathbf u,\cdot)\in \mathcal{OU}\big(\mathbf x, \widehat{(\bf y_d,\pi)}\big), \textrm{with} \ \widehat{(\bf y_d,\pi)} \in  2^{\mathbb Y_d}\times (\mathcal{P}_{\Pi}\bigcup\mathcal{R}_{\Pi})^{|\hat{\mathbb Y}_d|}\bigg\}.$$
 Then, we have 
 $$\max_{\mathbf u\in \mathcal{U}(\mathbf x)}\min_{(\mathbf y_c,\mathbf y_d)\in \mathcal{Y}(\mathbf x,\mathbf u)} (\mathbf c_{2,c}\mathbf y_c+\mathbf c_{2,d}\mathbf y_d) = \max_{\mathbf u\in \hat{\mathcal{U}}(\mathbf x)} \min_{(\mathbf y_c,\mathbf y_d)\in \mathcal{Y}(\mathbf x,\mathbf u)} (\mathbf c_{2,c}\mathbf y_c+\mathbf c_{2,d}\mathbf y_d).$$
\end{prop}
\begin{proof}
	According to the definition of set $\mathcal{OU}$, we have $\hat{\mathcal{U}}(\mathbf x)\subseteq \mathcal{U}(\mathbf x)$.   Hence, it is straightforward to conclude that 
	 $$\max_{\mathbf u\in \mathcal{U}(\mathbf x)}\min_{(\mathbf y_c,\mathbf y_d)\in \mathcal{Y}(\mathbf x,\mathbf u)} (\mathbf c_{2,c}\mathbf y_c+\mathbf c_{2,d}\mathbf y_d) \geq \max_{\mathbf u\in \hat{\mathcal{U}}(\mathbf x)} \min_{(\mathbf y_c,\mathbf y_d)\in \mathcal{Y}(\mathbf x,\mathbf u)} (\mathbf c_{2,c}\mathbf y_c+\mathbf c_{2,d}\mathbf y_d).$$

	Then, it would be sufficient to show the left-hand-side is less than or equal to the right-hand-side. Consider the complete set $\bb Y_d$ and let 
	$$\hat{\cal{U}}_{\bb{Y}_d}(\mathbf x)=\Big\{\mathbf u: (\mathbf u, \lambda)\in \mathcal{OU}\big(\mathbf x, \widehat{(\bf y_d,\pi)}\big), \textrm{with} \ \widehat{(\bf y_d,\pi)} \in  \mathbb Y_d\times (\mathcal{P}_{\Pi}\bigcup\mathcal{R}_{\Pi})^{|{\mathbb Y}_d|}\Big\}.$$
	 Indeed, following \eqref{eq_eta_inner}-\eqref{eq_eta_inner_subset} and the definition of set $\mathcal{OU}$, we have
	 \begin{align*} 
	 	\label{eq_max-min-to-max}
	 	\begin{split}
	 	&\max_{\mathbf u\in \mathcal{U}(\mathbf x)}\min_{(\mathbf y_c,\mathbf y_d)\in \mathcal{Y}(\mathbf x,\mathbf u)} (\mathbf c_{2,c}\mathbf y_c+\mathbf c_{2,d}\mathbf y_d)\\
	 	=&\max_{\mathbf u\in \mathcal{U}(\mathbf x)}\Big\{\hat\eta:\hat\eta\leq\mathbf{c}_{2,d}\mathbf y_d + \min \left\{\mathbf c_{2,c}\mathbf y_c: \mathbf y_c\in\mathcal Y(\mathbf x, \mathbf u, \mathbf y_d)\right\}, \ \forall \bf y_d\in \mathbb{Y}_d\Big\}\\
	 	=&\max_{\bsm{\pi}\in (\mathcal{P}_{\Pi}\bigcup\mathcal{R}_{\Pi})^{|{\bb Y}_d|} }\Big\{\hat \eta: \mathbf u\in \mathcal{U}(\mathbf x),\hat \eta\leq \mathbf{c}_{2,d}\mathbf y_d^t + (\mathbf d-\mathbf B_1\mathbf x-\mathbf E\mathbf u-\mathbf B_{2,d}\mathbf y^t_d)^\intercal\pi^{t}, t=1,\dots, |\bb{Y}_d|\Big\}\\
	 	=&\max_{\mathbf u\in \hat{\mathcal{U}}_{\bb Y_d}(\bf x)}\max_{\bsm{\pi}\in (\mathcal{P}_{\Pi}\bigcup\mathcal{R}_{\Pi})^{|\bb Y_d|} }\Big\{\hat \eta: \hat \eta\leq \mathbf{c}_{2,d}\mathbf y_d^t + (\mathbf d-\mathbf B_1\mathbf x-\mathbf E\mathbf u -\mathbf B_{2,d}\mathbf y^t_d)^\intercal\pi^{t}, t=1,\dots, |\bb{Y}_d|\Big\}\\
	 	\leq&\max_{\mathbf u\in \hat{\mathcal{U}}_{\bb Y_d}(\bf x)}\Big\{\hat \eta: \hat \eta\leq \mathbf{c}_{2,d}\mathbf y_d^t + \min \left\{\mathbf c_{2,c}\mathbf y_c: \mathbf y_c\in\mathcal Y(\mathbf x, \mathbf u, \mathbf y^t_d)\right\},t=1,\dots, |\mathbb{Y}_d|\Big\}\\
	 	=&\max_{\mathbf u\in \hat{\mathcal{U}}_{\bb Y_d}(\bf x)}\min_{(\mathbf y_c,\mathbf y_d)\in \mathcal{Y}(\mathbf x,\mathbf u)} (\mathbf c_{2,c}\mathbf y_c+\mathbf c_{2,d}\mathbf y_d)
	 	\end{split}
	 \end{align*}
  As $\min \left\{\mathbf c_{2,c}\mathbf y_c: \mathbf y_c\in\mathcal Y(\mathbf x, \mathbf u, \mathbf y^t_d)\right\}=\infty$ if $\mathcal Y(\mathbf x, \mathbf u, \mathbf y^t_d)=\emptyset$,  it allows us to consider $(\mathbf y_c,\mathbf y_d)\in \mathcal{Y}(\mathbf x,\mathbf u)$ only. Hence, the last equality follows. Moreover, given the fact that 
    $\mathbb{Y}_d\in 2^{\mathbb{Y}_d}$, we have $\hat{\cal{U}}_{\bb Y_d}\subseteq \hat{\cal U}(\bf x)$, and 
	 $$\max_{\mathbf u\in \mathcal{U}(\mathbf x)}\min_{(\mathbf y_c,\mathbf y_d)\in \mathcal{Y}(\mathbf x,\mathbf u)} (\mathbf c_{2,c}\mathbf y_c+\mathbf c_{2,d}\mathbf y_d) \leq \max_{\mathbf u\in \hat{\mathcal{U}}(\mathbf x)} \min_{(\mathbf y_c,\mathbf y_d)\in \mathcal{Y}(\mathbf x,\mathbf u)} (\mathbf c_{2,c}\mathbf y_c+\mathbf c_{2,d}\mathbf y_d),$$
	 which leads to the desired result.	
\end{proof}

It is worth highlighting that $\hat{\mathcal{U}}(\mathbf x)$ is the union of a finite number of projections of $\mathcal{OU}$ sets, given that $2^{\mathbb{Y}_d}$ and  $(\mathcal{P}_{\Pi}\bigcup\mathcal{R}_{\Pi})^{|\hat{\mathbb Y}|_d}$ are finite. As a consequence of Proposition \ref{prop_OU_max-min_reform}, we present another reformulation using $\mathcal{OU}$ sets. 

\begin{thm}
	 When it has an MIP recourse, formulation $\mathbf{2-Stg \ RO}$ in \eqref{eq_2RO} (and its other  equivalences) is equivalent to the following bilevel linear optimization program.
	\begin{subequations}
		\label{eq_2stgRO_Omega_MIP}
		\begin{align}
			w^*=\min \ & \ \mathbf{c}_1\mathbf x+ \eta   \\
			\mathrm{s.t.} \ & \ \mathbf x\in \mathcal{X}\\
			&  \bigg\{\eta\geq \mathbf{c}_{2,c}\mathbf y^{\widehat{(\bf y_d,\pi)}}_c+\mathbf{c}_{2,d}\mathbf y^{\widehat{(\bf y_d,\pi)}}_d\\
			& \ \ \mathbf{B}_{2,c}\mathbf y^{\widehat{(\bf y_d,\pi)}}_c+\mathbf{B}_{2,d}\mathbf y^{\widehat{(\bf y_d,\pi)}}_d+\mathbf{E}_c\mathbf{u}^{\widehat{(\bf y_d,\pi)}}_c\geq \mathbf d-\mathbf B_1\mathbf x\\
			& \ \ \mathbf{y}^{\widehat{(\bf y_d,\pi)}}_{c}\geq \mathbf{0},\ \mathbf{y}^{\widehat{(\bf y_d,\pi)}}_d\in \mathbb{Y}_d\\
			& \ \ (\mathbf u^{\widehat{(\bf y_d,\pi)}}, \cdot)\in \mathcal{OU}\big(\mathbf x, \widehat{(\bf y_d,\pi)}\big) \bigg\}  \ \forall \widehat{(\bf y_d,\pi)}\in 2^{\mathbb Y_d}\times (\mathcal{P}_{\Pi}\bigcup\mathcal{R}_{\Pi})^{|\hat{\mathbb Y}_d|}\label{eq_RMIP_subset_OU}
		\end{align}
	\end{subequations}
\end{thm}

\begin{cor}
	\label{cor_Pi_partial_KKTOmegaMIDDURE} 
	Let $\widehat{(\hat{\mathbb Y}_d, \hat{ \bsm\pi})}$  denote a set of $\widehat{(\bf y_d,\pi)}$'s. Formulation \eqref{eq_2stgRO_Omega_MIP} defined on $\widehat{(\hat{\mathbb Y}_d, \hat{ \bsm\pi})}$ is a relaxation to $\mathbf{2-Stg \ RO}$ in \eqref{eq_2RO}, and its optimal value is less than or equal to $w^*$.   \hfill $\square$
\end{cor}

\begin{rem}
 Note that the aforementioned derivation leads to the unified cutting sets, regardless the feasibility of $\bf x$. As shown in the next subsection, a modified $\max-\min$ problem provides an easier approach to detect $\bf x$'s infeasibilty and helps us generate cutting sets to enforce feasibilty. We mention that although \eqref{eq_2stgRO_Omega_MIP} is a rather complex equivalence with two types of enumerations, it allows us to compute an exact solution by implementing (parametric) C\&CG method in a nested fashion.  
\end{rem}

\subsection{Handling MIP Recourse by A Nested Parametric C\&CG}
The presented nested parametric C\&CG method generalizes and strengthens the nested basic C\&CG presented in \cite{zhao2012exact} that is developed to handle $\mathbf{2-Stg \ RO}$ with MIP recourse in the context of DIU. Different from \cite{zhao2012exact}, it  includes two inner C\&CG subroutines, one for feasibility and one for optimality, respectively. Indeed, the first one for feasibility eliminates  the dependence on the relatively complete recourse property, which has been assumed to ensure the algorithm's applicability. Another new feature is that both subroutines have a two-phase structure to ensure that $\mathcal{OU}$ set captures worst-case scenarios of $\mathcal U(\bf x)$  precisely.  The outer parametric C\&CG procedure serves as the main framework on top of those inner subroutines. We next provide their detailed descriptions.

\subsubsection {Inner C\&CG Subroutine for Feasibility}
Assuming that a fixed $\mathbf x^*$ is given, we employ the following $\max-\min$ problem to investigate 
the feasibility issue of $\bf x^*$. It resembles the $\max-\min$ portion of $\mathbf{2-Stg \ RO}$ but has a new variable $\tilde{ y}$ and a modified objective function, 
	\begin{align}
    \label{eq_max-min_substructure_feasibility}
		\eta_f(\mathbf x^*) = \max_{\mathbf u\in \mathcal{U}(\mathbf x^*)} \ \ \min_{(\mathbf y_c,\mathbf y_d,\tilde{ y})\in \tilde{\mathcal{Y}}(\mathbf x^*,\mathbf u)} \tilde {y} 
	\end{align}
where  $\tilde{\mathcal Y}(\mathbf x^*, \mathbf u)=\{\mathbf B_{2,c}\mathbf y_c + \mathbf B_{2,d}\mathbf y_d+ \mathbf 1\tilde{y}\geq \mathbf d-\mathbf B_1\mathbf x^*-\mathbf E\mathbf u, \ \mathbf y_c\geq \mathbf 0, \ \mathbf y_d\in \mathbb{Y}_d,\ \tilde{y}\geq 0\}$.
\begin{lem}
	A given $\mathbf x^*$ is feasible with respect to $\mathcal{U}(\mathbf x^*)$ if and only if $\eta_f(\mathbf x^*)=0$. If $\eta_f(\mathbf x^*) > 0$, $\mathbf x^*$ is infeasible to an optimal $\mathbf u$ of  \eqref{eq_max-min_substructure_feasibility}. \hfill $\square$
\end{lem}

Because $\mathcal U(\mathbf x^*)$ is bounded, it is clear that $\eta_f(\mathbf x^*)$ is of a finite value. Also, there is no feasibility issue associated with this model. Indeed, by enumerating choices of discrete variables and through duality, \eqref{eq_max-min_substructure_feasibility} can be reformulated with a max-min-max structure as that in \eqref{eq_max-min_DAD1}, which is solvable by basic C\&CG. We next exploit this structure to develop an inner C\&CG subroutine. In addition to $\mathbf x^*$, we further assume that some $\mathbf u^*\in \mathcal{U}(\mathbf x^*)$ is given to define inner subproblem $\mathbf{ISP_f}$, which is actually the minimization problem of \eqref{eq_max-min_substructure_feasibility}.
\begin{eqnarray}
	\label{eq_SP-if}
		\mathbf{ISP_f}: r_f(\mathbf x^*,\mathbf u^*) =   \min \Big\{\tilde{ y}: (\mathbf y_c,  \mathbf y_d, \tilde{y})\in \tilde{\mathcal Y}(\mathbf x^*, \mathbf u^*)\Big\}
\end{eqnarray}
For this MIP program, we directly compute its optimal value $r_f(\mathbf x^*,\mathbf u^*)$ and optimal $(\bf y^*_c,\mathbf y^*_d, \tilde y^*)$. Then, we define the inner master problem, which can be seen as in  the form \eqref{eq_max-min_DAD2} with respect to the objective function of \eqref{eq_max-min_substructure_feasibility} and $\hat{\mathbb{Y}}_d =\{\bf y^1_d, \dots, \bf y^{|\hat{\mathbb{Y}}_d |}_d\}$, a subset of $\mathbb{Y}_d$.  
\begin{subequations}
	\label{eq_MP-if}
	\begin{align}
		\mathbf{IMP_f}: \ \hat\eta_f(\mathbf x^*, \hat{\mathbb{Y}}_d)= \max \ & \hat\eta_f \\
		\mathrm{s.t.}\ & \mathbf u  \in \mathcal{U}(\mathbf x^*)\\
		& \Big\{\hat \eta_f\leq \min\{\tilde {y}^t: \mathbf B_{2,c}\mathbf y^t_c + \bf1\tilde {y}^t\geq \mathbf d-\mathbf B_1\mathbf x^*-\mathbf E\mathbf u-\mathbf B_{2,d}\mathbf y^t_d  \notag\\
		&  \ \ \ \ \ \ \ \mathbf y^t_c\geq \mathbf 0,\ \tilde{y}^t\geq {0} \}\Big\} \  t=1,\dots, |\hat{\bb Y}_d | \label{eq_feasibility_master}
	\end{align}	
\end{subequations}
\begin{cor}
\label{cor_feasibility_sub}
If $\hat\eta_f(\mathbf x^*, \hat{\mathbb{Y}}_d)=0$ for some $\hat{\mathbb Y}_d$, we have $\eta_f(\bf x^*)=0$, i.e., $\bf x^*$ is feasible for $\cal U(\bf x^*)$. If $r_f(\mathbf x^*,\mathbf u^*)>0$ for some $\bf u^*$, we have $\eta_f(\bf x^*)>0$, i.e., $\bf x^*$ is infeasible for $\cal U(\bf x^*)$. 	  
\end{cor}
\begin{rem}
	\label{rmk_inner_feasibile}
	$(i)$ The \textit{relatively complete recourse property} is often assumed in the literature of multistage MIPs to circumvent the challenge of detecting and resolving the infeasibility issue of an MIP recourse for a first-stage solution. Specific to two-stage RO with MIP recourse~\cite{zhao2012exact}, it says the linear program portion of the recourse problem defined on $\bf y_c$ is feasible for any possible 3-tuple $(\mathbf x, \mathbf u, \mathbf y_d)\in \mathcal{X}\times\mathcal{U}(\bf x)\times\mathbb{Y}_d$. 
	Now, with the preceding subproblems and the following subroutine for feasibility, the infeasibility of $\bf x$ can be detected. Then, it can be addressed rigorously by implementing appropriate cutting sets in the outer C\&CG procedure. Hence, this relatively complete recourse assumption is not needed anymore. Since DDU generalizes DIU, it is clear that this subroutine can also be supplied to the nested basic  C\&CG to expand its solution capacity for DIU-based RO.\\
	$(ii)$ Constraints in \eqref{eq_feasibility_master} can be reformulated and simplified using  duality as 	
	\begin{equation}
	\label{eq_ISF_dual} 
	\big\{\hat\eta_f \leq  (\mathbf d-\mathbf B_1\mathbf x^*-\mathbf E\mathbf u-\mathbf B_{2,d}\mathbf y^t_d)^\intercal\pi^t, \  \mathbf{B}_{2,c}^\intercal \pi^t \leq \mathbf 0, \ \bf1^\intercal\pi^t\leq1,\ \pi^t\geq\mathbf{0}\big\}, \ t=1,\dots, |\hat{\mathbb{Y}}_d|.
	\end{equation}
	Then, problem $\mathbf{IMP_f}$ is converted into a bilinear MIP since it contains the product between $\mathbf u$ and $\pi^t$. It is currently solvable for some professional MIP solvers, while the computational performance might not be satisfactory. Actually, the relatively complete recourse property naturally holds for the LP portion of minimization problem of \eqref{eq_max-min_substructure_feasibility}, i.e., feasible for $\bf x^*$ and $(\bf u, \bf y_d)\in \cal U(\bf x^*)\times \bb Y_d$. Hence, we can replace the minimization problem in \eqref{eq_feasibility_master} by its  (linearized) KKT conditions,  which often yields a computationally more friendly formulation. Note that for either approach, computing $\mathbf{IMP_f}$ generates optimal values for $\pi^t$ for all $t$. \\
	$(iii)$It is worth mentioning that non-zero extreme points of the feasible set of \eqref{eq_ISF_dual} are normalized extreme rays of $\Pi$. This insight directly links $\widehat{(\bf y_d,\pi)}$ derived in the next algorithm to those appeared in \eqref{eq_2stgRO_Omega_MIP}, and therefore helps to define corresponding cutting sets.    
\end{rem}
 
Problems $\mathbf{ISP}_f$ and $\mathbf{IMP}_f$ are solved in a dynamic and iterative fashion to detect the feasibility issue of $\bf x^*$. Nevertheless, we might not obtain sufficient $\bf y_d$'s (and associated $\pi$'s) that precisely capture the set of worst-case scenarios within $\cal{U}(\bf x^*)$. 
To address such a situation, we introduce the following correction problem $\mathbf{ICP_f}$ to retrieve missing $\bf y_d$'s. 
\begin{eqnarray}
	\label{eq_CP-if}
	\begin{split}
		\mathbf{ICP_f}: c_f\big(\bf x^*, \widehat{(\bf y_d,\pi)}\big) =   \min \Big\{\tilde{ y}:(\mathbf y_c,\mathbf y_d,\tilde{ y})\in \tilde{\mathcal{Y}}(\mathbf x^*,\mathbf u), \ (\mathbf u, \cdot)\in \mathcal{OU}\big(\mathbf x^*, \widehat{(\bf y_d,\pi)}\big)\Big\}
	\end{split}
\end{eqnarray}

\begin{rem}
	By Remark \ref{rmk_inner_feasibile}(ii), solving \eqref{eq_MP-if} produces $\pi^t$ for each $\bf y^t_d$. If set $\widehat{(\bf y_d,\pi)}$ in \eqref{eq_CP-if} is same as the one obtained from computing \eqref{eq_MP-if}, it is clear that $c_f\big(\bf x^*, \widehat{(\bf y_d,\pi)}\big)$ is less than or equal to $\hat\eta_f(\mathbf x^*, \hat{\mathbb{Y}}_d).$ Indeed, if $c_f\big(\bf x^*, \widehat{(\bf y_d,\pi)}\big) =0$ and $\hat\eta_f(\mathbf x^*, \hat{\mathbb{Y}}_d)>0$, it indicates that we are able to, by leveraging 	some $\bf y_d\notin \bb{\hat Y}_d$, convert a $\bf u\in \mathcal{OU}\big(\mathbf x^*, \widehat{(\bf y_d,\pi)}\big)$ to a feasible scenario for $\bf x^*$. Note that such a scenario has been deemed an infeasible one for $\bf x^*$ if the discrete recourse decision is restricted to $\bb{\hat Y}_d$. This situation clearly requires us to expand $\bb{\hat Y}_d$ to eliminate the discrepancy. 
\end{rem}

In the following, we present all detailed steps for the inner C\&CG subroutine for feasibility. Note that $\mathbf{ICP_f}$ is called 
after we complete the iterative procedure between  $\mathbf{ISP}_f$ and $\mathbf{IMP}_f$, which renders this subroutine to have two phases. The explicit expression for cutting set in $\mathbf{Step \ 5}$ is ignored to simplify our exposition as it takes the form of \eqref{eq_feasibility_master}.

\noindent\underline{\textbf{Inner C\&CG Subroutine for Feasibility - (ISF-
		C\&CG)}}\\
\textbf{\underline{Phase I}} \vspace{-5pt}
\begin{description}
	\item[Step 1] Set $\hat{\mathbb Y}_d$ by an initialization strategy. 
	\item[Step 2] Solve inner master problem $\mathbf{IMP_f}$ defined with respect to $\hat{\mathbb Y}_d$, and derive optimal $\hat\eta_f(\mathbf x^*,\bb{\hat Y}_d)$, and optimal solution $(\mathbf u^*, \pi^{1*},\dots,\pi^{|\hat{\bb{Y}}_d|*}, \cdot)$. 
	\item[Step 3] If $\hat\eta_f(\mathbf x^*,\bb{\hat Y}_d)=0$, set $\eta_f(\mathbf x^*) =0$ and terminate the subroutine. 
	\item[Step 4] For given $(\mathbf x^*, \mathbf u^*)$, solve inner subproblem $\mathbf{ISP_f}$, derive optimal value $r_f(\mathbf x^*,\mathbf u^*)$ and solution $(\mathbf y_c^*,\mathbf y_d^*,\tilde y)$. 
	\item[Step 5] If $r_f(\mathbf x^*,\mathbf u^*)>0$, define $\widehat{(\bf y_d,\pi)}\equiv\big\{(\bf y^1_d, \pi^{1*}),\dots, (\bf y^{|\hat{\bb Y}_d|}_d, \pi^{|\hat{\bb Y}_d|*})\big\}$ and go to \textbf{Phase~II}. Otherwise, update $\hat{\mathbb{Y}}_d = \hat{\mathbb{Y}}_d \bigcup\{\mathbf y^*_d\}$ and go to $\mathbf{Step \ 2}$.
\end{description}
\textbf{\underline{Phase II}} \vspace{-5pt}
\begin{description}
	\item[Step 6] Solve correction problem (for feasibility) $\mathbf{ICP_f}$, derive its optimal value $c_f(\mathbf x^*,\widehat{(\bf y_d,\pi)})$, optimal solution $(\bf u^*, \bf y_d^*,\cdot)$, and corresponding $\pi^*$ obtained from computing $$\max\big\{(\mathbf d-\mathbf B_1\mathbf x^*-\mathbf E\mathbf u^*-\mathbf B_{2,d}\mathbf y^*_d)^\intercal\pi:  \mathbf{B}_{2,c}^\intercal \pi \leq \mathbf 0, \ \mathbf{1}^\intercal\pi\leq 1,\ \pi\geq\bf0\big\}.$$
	\item[Step 7] If $c_f\big(\mathbf x^*,\widehat{(\bf y_d,\pi)}\big)=0$, update $\widehat{(\bf y_d,\pi)}=\widehat{(\bf y_d,\pi)}\bigcup\{(\mathbf y^*_d,\pi^*)\}$ and go to $\mathbf{Step \ 6}$. Otherwise, set $\eta_f(\mathbf x^*) = r_f(\mathbf x^*,\mathbf u^*)$, return $\widehat{(\bf y_d,\pi)}$, and terminate.
\end{description}

\begin{rem}
	$(i)$ Indeed, correction problem $\mathbf{ICP_f}$ (and $\mathbf{ICP_o}$ introduced for the next subroutine) can be used reversely to eliminate unnecessary one from $\widehat{(\bf y_d,\pi)}$, i.e., whose removal does not change $c_f\big(\bf x^*, \widehat{(\bf y_d,\pi)}\big)$. By examining each component of $\widehat{(\bf y_d,\pi)}$ and deleting unnecessary ones, we can reduce its size and hence the complexity of set $\cal{OU}$.\\        
	$(ii)$ It is not necessary to use a single-dimension variable $\tilde y$ to examine feasibility of $\bf x^*$. It is introduced to help us understand the transition between MIP and LP recourse, and to simplify our convergence analysis for the main solution algorithm. Indeed, similar to \eqref{eq_SP1}, a multi-dimensional $\tilde{\bf y}$ is probably more effective in the feasibility examination, especially when the recourse problem has equality constraints. Certainly some simple changes should be made in defining sub-, master and correction problems for \textbf{ISF-C\&CG} if $\tilde{\bf y}$ is employed.
\end{rem}

Again, different initialization strategies can be applied to set $\hat{\mathbb Y}_d$ in \textbf{Step 1}. In addition to the naive strategy, another one is to set $\hat{\mathbb Y}_d$ by inheriting some or all $\mathbf y_d$'s contained in $\hat{\bb Y}_d$ from the previous $\bf{ISF-C\&CG}$ execution, referred to as the \textit{``inheritance initialization''}. We note that it may help us reduce the number of iterations involved, but computationally is not as stable as the naive one. 
After carrying out $\bf{ISF-C\&CG}$, if $\eta_f(\mathbf x^*)=0$ we will execute next the inner C\&CG subroutine for optimality, which is described in the following. For the other case where $\eta_f(\mathbf x^*)>0$, we will generate one cutting set to strengthen the master problem of the main outer procedure. 

\subsubsection{Inner C\&CG Subroutine for Optimality}
Assuming examination by the previous inner subroutine has been passed, we develop the subroutine for optimality to compute  $\eta_o(\bf x^*)$, the worst case performance of $\mathbf x^*$ as in \eqref{eq_max-min_substructure0}. The inner subproblem for this subroutine is simply the recourse problem for given $(\mathbf x^*,\bf u^*)$.
\begin{eqnarray}
	\label{eq_SP-o}
	\begin{split}
		\mathbf{ISP_o}: \ \ r_o(\mathbf x^*,\mathbf u^*) = \min_{(\mathbf y_c,\mathbf y_d)\in \mathcal{Y}(\mathbf x^*,\mathbf u^*)} (\mathbf c_{2,c}\mathbf y_c+\mathbf c_{2,d}\mathbf y_d)
	\end{split}
\end{eqnarray}
Similar to \eqref{eq_MP-if}, inner master problem $\mathbf{IMP_o}$ is defined with respect to subset $\hat{\mathbb Y}_d$. 
\begin{subequations}
	\label{eq_MP-io}
	\begin{align}
		\mathbf{IMP_o}: \ \hat\eta_o(\mathbf x^*)= \max \ & \hat\eta_o \\
		\mathrm{s.t.}\ &\mathbf u  \in \mathcal{U}(\mathbf x^*)\\
		&\Big\{\hat\eta_o  \leq \mathbf{c}_{2,d}\mathbf y^t_d + \min\{\bf{c}_{2,c}\bf y^t_c: \mathbf{B}_{2,c}\bf y^t_c\geq \mathbf d-\mathbf B_1\mathbf x-\mathbf E\mathbf u-\mathbf B_{2,d}\mathbf y^t_d \notag\\ 
		&  \ \ \  \ \bf y_{2,c}\geq \bf 0\}\Big\}\ t=1,\dots, |\hat{\mathbb{Y}}_d|
       \label{eq_MP-io-const}
	\end{align}	
\end{subequations}

\begin{rem}
\label{rem_inner_optimality}
 By making use of duality, constraints in \eqref{eq_MP-io-const} can be replaced
 by 
 $$\Big\{\hat\eta_o  \leq \mathbf{c}_{2,d}\mathbf y^t_d + (\mathbf d-\mathbf B_1\mathbf x-\mathbf E\mathbf u-\mathbf B_{2,d}\mathbf y^t_d)^\intercal\pi^t, \ \bf B^\intercal_{2,c}\pi^t\leq \bf c^\intercal_{2,c}, \pi^t\geq \bf 0 \Big\} \ t=1,\dots, |\hat{\mathbb{Y}}_d|.$$  
 This reformulation holds even if some of lower-level minimization problem is infeasible, where the dual problem is unbounded.  
 Then, $\mathbf{IMP_o}$ is converted into a mixed integer bilinear program, computable by some professional solvers. Also, if the relatively complete recourse property holds, every minimization problem can be replaced by its (linearized) KKT conditions,  converting $\mathbf{IMP_o}$ into an MIP that probably is computationally more friendly.
\end{rem}

Similar to the feasibility subroutine, the correction problem for optimality, i.e., $\mathbf{ICP_o}$, is introduced next to help us precisely capture the set of worst-case scenarios within in $\cal U(\bf x^*)$.
\begin{eqnarray*}
	\label{eq_CP-io}
	\begin{split}
		\mathbf{ICP_o}: c_o\big(\bf x^*,\widehat{(\bf y_d,\pi)}\big) = \min\!\Big\{\bf{c}_{2,c}\bf y_c+\mathbf{c}_{2,d}\mathbf y_d: (\mathbf y_c,\mathbf y_d)\in \mathcal{Y}(\mathbf x^*,\mathbf u), \ (\mathbf u, \cdot)\in \mathcal{OU}\big(\mathbf x^*, \widehat{(\bf y_d,\pi)}\big)\Big\}\\
	\end{split}
\end{eqnarray*}

In the following, we present the inner C\&CG subroutine for optimality.

\noindent\textbf{\underline{Inner C\&CG Subroutine for Optimality - (ISO-C\&CG)}}\\
\textbf{\underline{Phase I}} \vspace{-5pt}
\begin{description}
	\item[Step 1] Set $LB = -\infty$, $UB=+\infty$, and set $\hat{\mathbb Y}_d$ by an initialization strategy. 
	\item[Step 2] Solve inner master problem $\mathbf{IMP_o}$ defined with respect to $\hat{\mathbb Y}_d$, and derive optimal value $\hat\eta_o(\mathbf x^*)$, optimal solution $(\bf u^*, \pi^{1*}, \dots, \pi^{|\hat{\mathbb Y}_d|*}, \cdot)$. Update $UB=\hat\eta_o(\mathbf x^*)$. 
	\item[Step 3] For given $(\mathbf x^*, \mathbf u^*)$, solve inner subproblem $\mathbf{ISP_o}$, and derive optimal value $r_o(\mathbf x^*,\mathbf u^*)$, and optimal solution $(\mathbf y^*_c,\mathbf y_d^*)$. Update 
	$LB=\max\{LB, r_o(\mathbf x^*,\mathbf u^*)\}$.
	\item [Step 4] If $UB-LB\leq T\!O\!L$, create $\widehat{(\bf y_d,\pi)}\equiv\big\{(\bf y^1_d, \pi^{1*}),\dots, (\bf y^{|\hat{\bb Y}_d|}_d, \pi^{|\hat{\bb Y}_d|*})\big\}$ and go to \textbf{Phase~II}. Otherwise, update $\hat{\mathbb{Y}}_d = \hat{\mathbb{Y}}_d \bigcup\{\mathbf y^*_d\}$ and go to $\mathbf{Step \ 2}$.
\end{description}
\textbf{\underline{Phase II}} \vspace{-5pt}
\begin{description}
	\item[Step 5] Solve correction problem (for optimality) $\mathbf{ICP_o}$, derive its optimal value $c_o\big(\mathbf x^*,\widehat{(\bf y_d,\pi)}\big)$, optimal solution $(\bf u^*, \bf y_d^*,\cdot)$, and corresponding $\pi^*$ obtained from computing  
	$$\max\big\{(\mathbf d-\mathbf B_1\mathbf x^*-\mathbf E\mathbf u^*-\mathbf B_{2,d}\mathbf y^*_d)^\intercal\pi:  \mathbf{B}_{2,c}^\intercal \pi \leq \mathbf c^{\intercal}_{2,c}, \ \pi\geq \mathbf 0\big\}.$$
	\item[Step 6] If $c_o(\mathbf x^*,\widehat{(\bf y_d,\pi)})<LB$, update $\widehat{(\bf y_d,\pi)}=\widehat{(\bf y_d,\pi)}\bigcup\{(\mathbf y^*_d,\pi^*)\}$ and go to $\mathbf{Step \ 5}$. Otherwise, set $\eta_o(\mathbf x^*)=UB$, return $\widehat{(\bf y_d,\pi)}$, and terminate.
\end{description}

Similar to \textbf{ISF-C\&CG}, different initialization strategies can be applied to set $\hat{\mathbb Y}_d$. A useful one is to populate it with some $\bf y_d$'s generated in \textbf{ISF-C\&CG} upon exit, which are deemed effective in dealing with $\bf u$. For inner C\&CG subroutines of both feasibility and optimality,  every call of subproblem or correction problem (except the last call) generates a new $\bf y_d\in \bb{Y}_d$. So, the next iteration complexity result simply follows.

\begin{prop}
	\label{prop_complexity_inner_opt} For a given $\mathbf x^*$, the number of iterations, including those from both $\textbf {Phase I}$ and $\textbf{Phase II}$, for either inner C\&CG subroutine is of $O(|\mathbb Y_d|)$.\hfill $\square$
\end{prop}

\subsubsection{The Main Outer C\&CG Procedure}

The primary component of the outer  C\&CG procedure is its master problem. We present it with the unified cutting sets in the following. Assume that $\widehat{(\hat{\mathbb Y}_d,\hat{\bsm\pi})}$ is the set of $\widehat{(\bf y_d,\pi)}$'s obtained from inner C\&CG subroutines (including both feasibility and optimality ones) in all previous iterations, as well as from initialization. We have     
	\begin{subequations}
	\label{eq_2stgRO_Master_Omega_MIP}
	\begin{align}
	\mathbf{OMP:} \  \underline w=\min \ & \ \mathbf{c}_1\mathbf x+ \eta   \\
		\mathrm{s.t.} \ & \ \mathbf x\in \mathcal{X}\\
		&  \Big\{\eta\geq \mathbf{c}_{2,c}\mathbf y^{\widehat{(\bf y_d,\pi)}}_c+\mathbf{c}_{2,d}\mathbf y^{\widehat{(\bf y_d,\pi)}}_d\\
		& \ \ \mathbf{B}_{2,c}\mathbf y^{\widehat{(\bf y_d,\pi)}}_c+\mathbf{B}_{2,d}\mathbf y^{\widehat{(\bf y_d,\pi)}}_d\geq \mathbf d-\mathbf B_1\mathbf x-\mathbf{E}\mathbf{u}\\
		& \ \ \mathbf{y}^{\widehat{(\bf y_d,\pi)}}_{c}\geq \mathbf{0},\ \mathbf{y}^{\widehat{(\bf y_d,\pi)}}_d\in \mathbb{Y}_d\\
		& \ \ (\mathbf u, \cdot)\in \mathcal{OU}\big(\mathbf x, \widehat{(\bf y_d,\pi)}\big) \Big\}\ \forall \widehat{(\bf y_d,\pi)}\in \widehat{(\hat{\mathbb Y}_d,\hat{\bsm\pi})}. \label{eq_RMIP_MP_subset_OU}
	\end{align}
\end{subequations}
Next, we describe the main outer procedure. \\
\noindent\textbf{\underline{Outer C\&CG Procedure}}
\begin{description}
	\item[Step 1] Set $LB = -\infty$, $UB=+\infty$, and set $\widehat{(\hat{\mathbb Y}_d,\hat{\bsm\pi})}$ by an initialization strategy. 
	\item[Step 2] Solve outer master problem $\mathbf{OMP}$ defined with respect to $\widehat{(\hat{\mathbb Y}_d,\hat{\bsm\pi})}$. If infeasible, report the infeasibility of $\mathbf{2-Stg \ RO}$ and terminate. Otherwise, derive optimal value $\underline w$,  optimal solution $(\mathbf x^*,\cdot)$ and update $LB=\underline w$. 
	\item[Step 3] For given $\mathbf x^*$, call the inner subroutine for feasibility \textbf{ISF-C\&CG}, and obtain $\eta_f(\mathbf x^*)$, as well as $\widehat{(\bf y_d,\pi)}$ if $\eta_f(\mathbf x^*)>0$.   
	\item[Step 4] Cases based on $\eta_f(\mathbf x^*)$.
	\begin{description}
		\item[(Case A): $\eta_f(\mathbf x^*)>0$] \textrm{}
		 \\ Update $\widehat{(\hat{\mathbb Y}_d,\hat{\bsm\pi})}=\widehat{(\hat{\mathbb Y}_d, \hat{ \bsm\pi})}\bigcup\big\{\widehat{(\bf y_d,\pi)}\big\}$. 
		 \item[(Case B): $\eta_f(\mathbf x^*)=0$]  \textrm{}\\
		 Call the inner subroutine for optimality \textbf{ISO-C\&CG}, obtain  $\eta_o(\mathbf x^*)$ and $\widehat{(\bf y_d,\pi)}$. Update $\widehat{(\hat{\mathbb Y}_d,\hat{\bsm\pi})}=\widehat{(\hat{\mathbb Y}_d, \hat{ \bsm\pi})}\bigcup\big\{\widehat{(\bf y_d,\pi)}\big\}$ and $UB=\min\{UB,\ \hat\eta_o(\mathbf x^*)\}$. 
	\end{description} 
	\item [Step 5] If $UB-LB\leq T\!O\!L$, return $\bf x^*$ and terminate. Otherwise, go to \textbf{Step 2}. 
\end{description}

\begin{rem}
Again, various initialization strategies can be applied to set $\widehat{(\hat{\mathbb Y}_d,\hat{\bsm\pi})}$. For example, the naive one is to set it as an  empty set. Also, as showed in the next section, the linear program relaxation of the recourse problem can be used to generate effective $\widehat{(\bf y_d,\pi)}$'s, which also helps to yield valid lower and upper bounds to develop an approximation scheme.  
\end{rem}

\subsection{Convergence and Complexity}
Compared to other parametric C\&CG variants, the analyses of convergence and complexity issues for the aforementioned nested parametric C\&CG are rather different. One main reason is the sophisticated structure of $\bf{LP}\big(\bf x,\widehat{(\bf y_d,\pi)}\big)$ in \eqref{eq_eta_inner_subset}. Note that its optimal solution generally is not an extreme point or not defined by a basis of $\cal U(\bf x)$. Actually, its optimal solution is determined by, in addition to $\cal{U}(\bf x)$, all $(\bf y_d, \pi)$'s in set $\widehat{(\bf y_d,\pi)}$.  
Hence, our convergence and complexity analyses mainly depend on $\widehat{(\bf y_d,\pi)}$.  Again, we assume that $T\!O\!L=0$ for all involved procedures. 
	  
\begin{thm}
\label{thm_miprecourse_conv}
The aforementioned nested parametric C\&CG  will not repeatedly generate any $\widehat{(\bf y_d,\pi)}$ unless it terminates. Upon termination, it either reports that $\mathbf{2-Stg \ RO}$ in (\ref{eq_2RO}) is infeasible, or converges to its optimal value with an exact solution.
\end{thm}

\begin{proof}
	We first consider the feasibility issue. Note that $\mathbf{2-Stg \ RO}$ is infeasible if $\bf{OMP}$ becomes infeasible.

	\noindent Claim 1: If $\widehat{(\bf y_d,\pi)^*}$ has been derived by one call of \textbf{ISF-C\&CG} in some outer iteration, it will not be produced in any following iterations.
	\begin{proof}[Proof of Claim 1:]
		Since $\widehat{(\bf y_d,\pi)^*}$ has already been derived by \textbf{ISF-C\&CG},  the formulation of $\mathbf{OMP}$ includes the  cutting set in any following iteration.
		\begin{align}
			\label{OuterCCGCut}
			\begin{split}
	\Big\{&\eta\geq \mathbf{c}_{2,c}\mathbf y^{\widehat{(\bf y_d,\pi)^*}}_c+\mathbf{c}_{2,d}\mathbf y^{\widehat{(\bf y_d,\pi)^*}}_d\\
	&\mathbf{B}_{2,c}\mathbf y^{\widehat{(\bf y_d,\pi)^*}}_c+\mathbf{B}_{2,d}\mathbf y^{\widehat{(\bf y_d,\pi)^*}}_d\geq \mathbf d-\mathbf B_1\mathbf x-\mathbf{E}\mathbf{u}\\
	&\mathbf{y}^{\widehat{(\bf y_d,\pi)^*}}_{c}\geq \mathbf{0},\ \mathbf{y}^{\widehat{(\bf y_d,\pi)^*}}_d\in \mathbb{Y}_d,\ (\mathbf u, \cdot)\in \mathcal{OU}\big(\mathbf x, \widehat{(\bf y_d,\pi)^*}\big) \Big\}.
\end{split}
		\end{align}
		Consider $\mathbf x^*$ obtained from computing $\mathbf{OMP}$ in one such iteration. We have
\begin{align*}
	0 = &\min\ \Big\{\tilde{ y}: \ \mathbf B_{2,c}\mathbf y_c + \mathbf B_{2,d}\mathbf y_d+ \mathbf{1}\tilde{y}\geq \mathbf d-\mathbf B_1\mathbf x^*-\mathbf E\mathbf u \\
	&\phantom{\min\ \Big\{\tilde{y}: \ \ \ } \mathbf{y}_{c}\geq \mathbf{0},\ \mathbf{y}_d\in \mathbb{Y}_d,\ \tilde y\geq0,\ (\mathbf u, \cdot)\in \mathcal{OU}\big(\mathbf x^*, \widehat{(\bf y_d,\pi)^*}\big)\Big\},
\end{align*}
which is $c_f\big(\bf x^*, \widehat{(\bf y_d,\pi)^*}\big)$ if $\widehat{(\bf y_d,\pi)^*}$ is output by \textbf{ISF-C\&CG}.   
Hence, once \textbf{ISF-C\&CG} is called, it either reports $\eta_f(\mathbf x^*)=0$ or yields set $\widehat{(\bf y_d,\pi)}$ different from $ \widehat{(\bf y_d,\pi)^*}$ to generate a new feasibility cutting set.
\end{proof}
 
\noindent Claim 2: Assume that  $\mathbf x^*$ generated by $\mathbf{OMP}$ is feasible, and $\widehat{(\bf y_d,\pi)^*}$ is then derived by subroutine \textbf{ISO-C\&CG}  in the current iteration. If $\widehat{(\bf y_d,\pi)^*}$ has been obtained by \textbf{ISO-C\&CG} in some previous iteration, we have $UB=LB$.
    
\begin{proof}[Proof of Claim 2:]
As $\widehat{(\bf y_d,\pi)^*}$  has already been previously identified by \textbf{ISO-C\&CG}, $\mathbf{OMP}$ has its corresponding cutting set (\ref{OuterCCGCut}) since then.  After computing $\bf{OMP}$ in the current iteration, we have
\begin{align*}
	\begin{split}
		LB\geq&\mathbf{c}_1\mathbf x^*+ \Big\{\mathbf{c}_{2,c}\mathbf y_c+\mathbf{c}_{2,d}\mathbf y_d: 		\ (\mathbf u, \cdot)\in \mathcal{OU}\big(\mathbf x^*, \widehat{(\bf y_d,\pi)^*}\big), \ (\bf y_c,\bf y_d)\in \mathcal{Y}(\bf x,\bf u)\Big\}\\
		=&\mathbf{c}_1\mathbf x^*+c_o\big(\bf x^*,\widehat{(\bf y_d,\pi)^*}\big)=\mathbf{c}_1\mathbf x^*+\eta_o(\bf x^*)\geq UB.
	\end{split}
\end{align*}
Hence, it leads to $UB=LB$.
	\end{proof}
Since both $\mathbb Y_d$ and the set of extreme points and rays of $\Pi$ are fixed and finite, it follows that the nested parametric C\&CG will terminate by either reporting the infeasibility of $\mathbf{2-Stg \ RO}$ or  converging to an exact solution after a finite number of iterations.
\end{proof}

Given the structure of $\Pi$ and Remark \ref{rmk_inner_feasibile}($iii$), the number of $\widehat{(\bf y_d,\pi)}$'s is bounded by $(|\mathcal P_\Pi|+|\mathcal R_\Pi|+1)^{|\mathbb{Y}_d|}$. So, we can easily bound the number of iterations of this algorithm.

\begin{cor}
\label{cor_mipddu_comp1}
The number of iterations of the algorithm before termination is bounded by $(|\mathcal P_\Pi|+|\mathcal R_\Pi|+1)^{|\mathbb{Y}_d|}$. Hence, the algorithm is of $O\big({\binom{n_{y_c}+\mu_y}{\mu_y}}^{|\mathbb{Y}_d|}\big)$ iteration complexity. $\hfill\square$
\end{cor}

Similar to other variants, this nested one converges  if the first-stage solution repeats.

\begin{prop}
	\label{prop_DupOpt2}
	If $\mathbf x^*$ is an optimal solution to $\mathbf{OMP}$ in both iterations $t_1$ and $t_2$ with $t_1<t_2$, the nested parametric C\&CG terminates and $\mathbf x^*$ is optimal to $\mathbf{2-Stg \ RO}$. Hence, if $\cal X$ is a finite set, the algorithm is of 	
	$O\big(\min\{\binom{n_{y_c}+\mu_y}{\mu_y}^{|\mathbb{Y}_d|}, |\cal X|\}\big)$ iteration complexity. $\hfill\square$
\end{prop}

\section{Extensions to More Complex Two-Stage RO}
\label{sect_MIDDUMIP}
As noted earlier, mixed integer DDU and recourse are much more powerful in modeling complex problems. In this section, we further our investigation to the most general case where both $\mathcal{U(\bf x)}$ and $\mathcal Y(\bf x,\bf u)$ are mixed integer sets. Basically, the aforementioned nested C\&CG can be extended to compute this general case exactly. Also, parametric C\&CG-MIU can be extended as an alternative method to derive approximate solutions. To minimize repetitions, we present the primary changes on top of those two procedures, while assure that readers should be able to retrieve the complete operations easily.

\subsection{Extending Nested C\&CG for Exact Computation}
\label{ssub_ExNeCCG}
Modifications on the nested parametric C\&CG are mainly made to its main outer C\&CG procedure. For two inner subroutines, almost no change needs to be made, except that mixed integer $\mathcal U(\bf x)$ is employed and two correction problems are revised accordingly. 

Let $(\bf u^*_d,\bf u^*_c)$ denote the optimal solution and  $\widehat{(\mathbf{ y}_d, \pi)}$ output from subroutine either for feasibility (i.e., $\bf{ISF-C\&CG}$) or for optimality (i.e., $\bf{ISO-C\&CG}$) in one iteration. Similar to \eqref{eq_LP_parametric_ed} and \eqref{eq_eta_inner_subset}, we introduce the following optimization  $\mathbf{LP}\big(\mathbf x, \mathbf u^*_d, \widehat{(\mathbf{y}_d, \pi)}\big)$ to support the development of the main outer procedure, recalling $\tilde{\mathcal U}$ has been defined in Section \ref{ssec_EnReform}.  
\begin{eqnarray*}
	\label{eq_MIDDUMIRE_parametric_ed}
\begin{split}
	\hat \eta\big(\mathbf x, \bf u^*_d, \widehat{(\mathbf  y_d, \pi)}\big) & = \max\Big\{\hat\eta -M(\bf 1^\intercal\tilde{\bf u}) :  \ (\bf u_c, \tilde{\bf u})\in \tilde{\mathcal U}(\bf x|\bf{u}^*_d)\\ 
	\hat \eta\leq & \mathbf{c}_{2,d}\mathbf y^t_d + (\mathbf d-\mathbf E_d\mathbf u^*_d-\mathbf B_1\mathbf x-\mathbf E_c\mathbf u_c-\mathbf B_{2,d}\mathbf y^{t}_d)^\intercal\pi^{t}, \ t=1,\dots, |\mathbb{\hat{Y}}_d|\Big\}
\end{split}
\end{eqnarray*}  

Accordingly, we can employ its KKT conditions or strong duality to characterize its optimal solution set. The KKT conditions based one is as the following.
	\begin{align*}
	\mathcal{OU}\big(\mathbf x, \bf u^*_d, \widehat{(\mathbf  y_d, \pi)}\big)\!=\!\left\{\!\begin{array}{l}
		\hat \eta\leq \mathbf{c}_{2,d}^\intercal\mathbf y^t_d + (\mathbf d-\mathbf E_d\mathbf u^*_d-\mathbf B_{2,d}\mathbf y^t_d-\mathbf B_1\mathbf x-\mathbf E_c\mathbf u_c)^\intercal\pi^{t},\\ 
		t=1,\dots, |\hat{\mathbb{Y}}_d|\\
		\mathbf F_c(\bf x)\mathbf u_c-\tilde{\bf u}\leq \mathbf h+\mathbf{Gx}-\mathbf F_d(\bf x)\mathbf u^*_d\\ 
		\mathbf F_c^\intercal(\bf x)\boldsymbol\lambda+\sum_{t=1}^{|\hat{\mathbb{Y}}_d|}\mathbf E_c^\intercal\pi^t\mu^t\geq\mathbf 0 \\
		M\bf 1-\bsm\lambda \geq \bf 0 \\
		\sum_{t=1}^{|\hat{\mathbb{Y}}_d|}\mu^t=1\\
		\mu^t\big(\mathbf{c}_{2,d}^\intercal\mathbf y^t_d + (\mathbf d-\bf E_d\bf u^*_d-\mathbf B_1\mathbf x-\mathbf E_c\mathbf u_c-\mathbf B_{2,d}\mathbf y^t_d)^\intercal\pi^t-\hat \eta\big)=0,\\
		t=1,\dots,|\hat{\mathbb{Y}}_d|\\
		\boldsymbol\lambda \circ \big(\mathbf h+\mathbf{Gx}- \bf F_d(\bf x)\bf u^*_d  -\mathbf F_c(\bf x)\mathbf{u}_c+\tilde{\bf u}\big)=\mathbf{0} \\
		\mathbf u_c\circ \big(\mathbf F_c^\intercal(\bf x)\boldsymbol\lambda+ \sum_{t=1}^{|\hat{\mathbb{Y}}_d|}\mathbf E^\intercal\pi^t\mu^t\big) = \mathbf 0 \\
		\tilde{\bf u}\circ(M\bf 1-\bsm\lambda) = \bf 0 \\
		\mathbf u_c\geq \mathbf{0},\ \tilde{\bf u}\geq \bf 0, \  \boldsymbol\lambda\geq \mathbf{0},\ \mu^t\geq0, \ t=1,\dots, |\hat{\mathbb Y}_d|
	\end{array}\!\right\}
\end{align*}

With the aforementioned $\mathcal{OU}$ set, the two correction problems for the feasibility and optimality subroutines, respectively,  are defined as follows.
\begin{eqnarray*}
	\begin{split}
		\mathbf{ICP_f}:\ & c_f\big(\bf x^*, \bf u^*_d, \widehat{(\bf y_d,\pi)}\big) =   \min \Big\{\tilde{y}: (\mathbf y_c,\mathbf y_d,\tilde{y})\in \tilde{\mathcal{Y}}(\mathbf x^*,\mathbf u_c,\mathbf u_d^*), \\ 
		& \phantom{c_f\big(\bf x^*, \bf u^*_d, \widehat{(\bf y_d,\pi)}\big) =   \min \Big\{\tilde{y}:\ } (\mathbf u_c, \cdot)\in \mathcal{OU}\big(\mathbf x^*, \bf u^*_d, \widehat{(\bf y_d,\pi)}\big)\Big\}\\
		\mathbf{ICP_o}:\ & c_o\big(\bf x^*,\bf u^*_d, \widehat{(\bf y_d,\pi)}\big) = \min\!\Big\{\bf{c}_{2,c}\bf y_c+\mathbf{c}_{2,d}\mathbf y_d: (\mathbf y_c,\mathbf y_d)\in \mathcal{Y}(\mathbf x^*,\mathbf u_c,\mathbf u_d^*), \\
		& \phantom{c_o\big(\bf x^*,\bf u^*_d, \widehat{(\bf y_d,\pi)}\big) = \min\!\Big\{\bf{c}_{2,c}\bf y_c+\mathbf{c}_{2,d}\mathbf y_d:\ } (\mathbf u_c, \cdot)\in \mathcal{OU}\big(\mathbf x^*,\bf u_d^*, \widehat{(\bf y_d,\pi)}\big)\Big\}\\
	\end{split}
\end{eqnarray*}

Let $\reallywidehat{\Big({\hat{\bb U}}_d, (\hat{\mathbb Y}_d,\hat{\bsm\pi})\Big)}$
denote the set of $(\bf u_d, \widehat{(\mathbf  y_d, \pi)})$'s obtained so far from executions of subroutines (and including those from initialization if applicable). Then,  the master problem $\bf{OMP}$ mainly consists of (unified) cutting sets defined with respect to each component of $\reallywidehat{\Big({\hat{\bb U}}_d, (\hat{\mathbb Y}_d,\hat{\bsm\pi})\Big)}$, in addition $\bf x\in \cal X$.  By making use of $\mathcal{OU}\big(\mathbf x, \bf u_d, \widehat{(\mathbf  y_d, \pi)}\big)$ and $\dot{u}(\tilde{\mathbf u})$ defined in~\eqref{eq_penVar_MIDDU}, the unified cutting sets are 
 	\begin{align*}
     \Big\{\eta&\geq \mathbf{c}_{2,c}\mathbf y^{(\bf u_d, \widehat{(\mathbf  y_d, \pi)})}_c+\mathbf{c}_{2,d}\mathbf y^{(\bf u_d, \widehat{(\mathbf  y_d, \pi)})}_d\\
 	& \ \ \mathbf{B}_{2,c}\mathbf y^{(\bf u_d, \widehat{(\mathbf  y_d, \pi)})}_c+\mathbf{B}_{2,d}\mathbf y^{(\bf u_d, \widehat{(\mathbf  y_d, \pi)})}_d+\mathbf{E}_c\mathbf{u}^{(\bf u_d, \widehat{(\mathbf  y_d, \pi)})}_c+\dot{u}(\tilde{\mathbf u}^{(\bf u_d, \widehat{(\mathbf  y_d, \pi)})})\geq \mathbf d-\mathbf{E}_d\mathbf{u}_d-\mathbf B_1\mathbf x\\\
 	& \ \ \mathbf{y}^{(\bf u_d, \widehat{(\mathbf  y_d, \pi)})}_{c}\geq \mathbf{0}, \  \ \mathbf{y}^{(\bf u_d, \widehat{(\mathbf  y_d, \pi)})}_d\in \mathbb{Y}_d\\
 	& \ \ (\mathbf u^{(\bf u_d, \widehat{(\mathbf  y_d, \pi)})}_c, \cdot )\in \mathcal{OU}\big(\mathbf x,  \bf u_d, \widehat{(\mathbf  y_d, \pi)} \big)\Big\}  \ \  \forall \ (\bf u_d, \widehat{(\mathbf  y_d, \pi)})\in \reallywidehat{\Big({\hat{\bb U}}_d, (\hat{\mathbb Y}_d,\hat{\bsm\pi})\Big)}.
 \end{align*}
 That is, whenever one call of some inner subroutine is completed and $\big(\mathbf u^*_d, \widehat{(\mathbf{y}_d, \pi)}\big)$ is output,  we expand set $\reallywidehat{\Big({\hat{\bb U}}_d, (\hat{\mathbb Y}_d,\hat{\bsm\pi})\Big)} =\reallywidehat{\Big({\hat{\bb U}}_d, (\hat{\mathbb Y}_d,\hat{\bsm\pi})\Big)}\bigcup\Big\{\big(\mathbf u^*_d, \widehat{(\mathbf{y}_d, \pi)}\big)\Big\}$ and then augment $\mathbf{OMP}$ in the  following iteration of the main outer procedure. 

Regarding the convergence of this variant, it can be proven that a new $\big(\mathbf u_d, \widehat{(\mathbf{y}_d, \pi)}\big)$ will be generated in every iteration by inner subroutines. Otherwise, we have $LB=UB$, leading to convergence. Since the total numbers of $\bf u_d$, $\bf y_d$ and $\pi$ are finite, respectively, the total number of distinct  $\big(\mathbf u_d, \widehat{(\mathbf{y}_d, \pi)}\big)$'s is finite. Such a result ensures the finite convergence, as well as can be used to bound the iteration complexity.

\subsection{Modifying Parametric C\&CG-MIU for Approximation}
\label{ssub_PCCGApprox}
As shown in \cite{zeng2022two}, standard parametric C\&CG can be applied, with rather simple changes, to approximately compute $\mathbf{2-Stg \ RO}$ with polytope DDU and MIP recourse. In particular, both lower and upper bounds on the objective function value of the approximate solution are available in the algorithm execution, providing a quantitative  evaluation for the approximation quality. Following that strategy, we modify parametric C\&CG-MIU accordingly so that it handles $\mathbf{2-Stg \ RO}$ with mixed integer DDU and MIP recourse. When mixed integer DDU reduces to a polytope,  it also reduces to that approximation variant of standard parametric C\&CG. We assume that the relatively complete recourse property holds, i.e., the continuous portion of the recourse
problem is feasible for $\big(\bf x, (\bf u_c, \bf u_d), \bf y_d\big)\in \mathcal X\times\mathcal U(\bf x)\times\mathbb Y_d$, allowing us to skip \textbf{Step 3} and \textbf{Step 4 (Case B)} of parametric C\&CG-MIU.

Modifications to parametric C\&CG-MIU are mainly made for constructing its subproblems. One change is to replace the minimization recourse problem of $\mathbf{SP2}$ in \eqref{eq_SP2} by the linear programming relaxation of the MIP recourse, bringing us the following $\mathbf{SP2}$. 
\begin{eqnarray}
	\label{eq_SP2-MIUMIP}
	\mathbf{SP2}: \ \max_{(\bf u_c, \bf u_d)\in \mathcal{U}(\mathbf x^*)} \min\{
	\mathbf c_2\mathbf y: \mathbf y\in \mathcal{Y}_r(\mathbf x^*, \mathbf u)\}
\end{eqnarray}
We also introduce a new subproblem, referred to as $\mathbf{SP4}$.
\begin{eqnarray}
	\label{eq_SP4-MIUMIP}
	\mathbf{SP4}: \ \tilde \eta_o(\mathbf x^*, \bf y^*_d)=\max_{(\bf u_c, \bf u_d)\in \mathcal{U}(\mathbf x^*)} \min\Big\{
	\mathbf c_{2,c}\mathbf y_c+\mathbf c_{2,d}\mathbf y_d^*: \mathbf y_c\in \mathcal Y(\mathbf x^*, \mathbf u, \mathbf y_d^*)\Big\}
\end{eqnarray}
Then, in \textbf{Step 4 (Case A)}, we compute $\mathbf{SP2}$ in \eqref{eq_SP2-MIUMIP}, instead of the one in \eqref{eq_SP2}, to derive optimal $(\mathbf u^*_{oc}, \mathbf u^*_{od})$ and corresponding extreme point $\pi^*$ of the dual problem associated with $\mathcal{Y}_r(\mathbf x^*, \mathbf u)$, and perform remaining operations in this step. Note that cutting set in the form of \eqref{eq_CCG_optimality} is defined with variables $(\bf y_c^{(\bf u^*_{od},\pi^*)}, \bf y_d^{(\bf u^*_{od},\pi^*)})$. Hence, integer restriction on $\bf y_d^{(\bf u^*_{od},\pi^*)}$ should be annexed to this cutting set.  Regardless of that $\mathcal Y_r$ replaces $\mathcal Y$ in $\mathbf{SP2}$, we have $\mathcal{OU}(\mathbf x, \mathbf u^*_{od}, \pi^*)\in \mathcal U(\bf x)$. The generated cutting sets are valid, and the resulting master problem keeps being an effective relaxation to $\mathbf{2-Stg \ RO}$. Therefore, the optimal value of master problem remains a valid lower bound.    

A critical change is made in $\mathbf{Step \ 5}$ where we compute the upper bound. Given $(\mathbf u^*_{oc}, \mathbf u^*_{od})$ output from $\mathbf{Step  \ 4}$, we solve the original MIP recourse problem defined on $\mathcal Y\big(\bf x^*, (\mathbf u^*_{oc}, \mathbf u^*_{od})\big)$ and derive its optimal solution $(\bf y^*_c, \bf y^*_d)$. With $(\bf x^*, \bf y^*_d)$, we then compute $\mathbf{SP4}$ in \eqref{eq_SP4-MIUMIP} and obtain its optimal value $\tilde \eta_o(\mathbf x^*, \bf y^*_d)$. Clearly, by fixing $\bf y_d$ to $\bf y^*_d$, 
$\tilde \eta_o(\mathbf x^*, \bf y^*_d)$ is larger than or equal to the actual $\eta_o(\bf x^*)$, i.e., $\bf x^*$'s worst case performance defined with the original MIP recourse. It indicates that $\bf c_1\bf x^*+\tilde \eta_o(\mathbf x^*, \bf y^*_d)$ is an effective upper bound to $\mathbf{2-Stg \ RO}$. Hence, we adopt $UB=\min\{UB, \bf c_1\bf x^*+\tilde \eta_o(\mathbf x^*, \bf y^*_d)\}$ to update $UB$ in $\mathbf{Step  \ 5}$.

With both lower bound $LB$ and upper bound $UB$  being available in the execution of this algorithm, a quantitative evaluation on the quality of the approximate solution is provided. As for the termination, we can stop the algorithm once $(i)$ a desired gap between $LB$ and $UB$ is achieved (which nevertheless may not be achievable); or $(ii)$ a repeated $\bf x^*$ is output from master problem; or $(iii)$ a pre-defined limit on the number of iterations is  reached. Actually,  according to our numerical 
study, this approximation scheme often produces strong solutions
with small amount of computational expenses, demonstrating a desirable advantage in computing practical-scale instances.

\section{Computational Study and Analysis}
\label{sec_MIPC}
In this section, we carry out a set of experiments to demonstrate the modeling power of mixed integer sets,
and the computational performance of the presented algorithms. Those algorithms are implemented by \texttt{Julia} with
\texttt{JuMP} and professional MIP solver \texttt{Gurobi~9.1} on a Windows PC with E5-1620 CPU and 32G RAM, and parameter $M$ is set to 10,000. Unless noted otherwise, the relative optimality tolerances of any algorithm is set to $.5\%$, with the time limit set to 3,600 seconds and solver's default settings. Test instances are generated for three variants of robust facility location problems. One is with a mixed integer DDU set and a linear program recourse, one is with a polytope DDU set and a mixed integer recourse, and the last one is with  mixed integer structures in both DDU set and  recourse.  

\subsection{Capturing Induced Demand by Mixed Integer DDU}
\label{FL_MIPDDU}
Consider the robust facility location model with uncertain induced demand. The first-stage decision determines the locations and service capacities of facilities. 
Then, after the demand is realized, service decision for each client-facility pair is made in the second  stage.
\begin{subequations}
	\begin{align}
	\mathbf{RFL-L}:	\min_{(\mathbf x_c,\mathbf x_d)\in\mathcal X}  \sum_{j\in J} (f_jx_{d,j}+a_jx_{c,j})+ \rho\max_{\mathbf u\in\mathcal U^I(\mathbf x)}\min_{\mathbf y\in\mathcal{Y}(\mathbf x,\mathbf u)} \sum_{i\in I}\sum_{j\in J}c_{ij}y_{ij}
	\end{align}
	 Set $I$ is the set of client sites, and $J\subseteq I$ consists of potential facility sites. Variable $x_{d,j}$ is binary that takes 1 if a facility is built at site $j$ and 0 otherwise; $x_{c,j}$ determines the service capacity of facility at site $j$; and $y_{ij}$ denotes the quantity of client $i$'s demand served by facility at $j$. Parameters $f_j$ and $a_j$ are construction and unit capacity cost, respectively;  $c_{ij}$ is the unit transportation or service cost between client $i$ and facility $j$; and coefficient $\rho$ integrates two stages' costs together. The first and second stages feasible sets are
	\begin{align}
		\mathcal X= \big\{(\mathbf x_c, \mathbf x_d)\in \mathbb{R}^{|J|}_+\times\{0,1\}^{|J|}: \underline A_jx_{d,j}\leq x_{c,j}\leq \overline A_jx_{d,j} \ \forall j\in J\big\}, \ \textrm{and}\\
		\mathcal Y(\mathbf x,\mathbf u)= \big\{\mathbf y\in \mathbb{R}^{|I|\times|J|}_+:  \sum_{j\in J}y_{ij}\geq u_i \ \forall i\in I, \  \sum_{i\in I}y_{ij}\leq x_{c,j} \ \forall j\in J\big\}
	\end{align}
\end{subequations}
with $\underline A_j$ and $\overline A_j$ being the lower and upper bounds for the service capacity if a facility is installed at site $j$.

Let $\mathbf u$ denote the actual demand with  $u_i=\underline{u}_i+\tilde{u}_i$, which are its nominal and induced demands for client $i$, respectively. It is often the case in practice that $\tilde u_i$ is related to facilities built at its neighboring sites $j\in J(i)$, while such a connection cannot be described exactly. Assume that two independent consulting services are employed to investigate such connection, each of which provides an interval estimation on the contribution of $j$'s facility on $\tilde u_i$. Denote those two intervals by $[\underline\zeta^i_j, \overline\zeta_j^i]$ and $[\underline\xi^i_j, \overline\xi_j^i]$, which are typically different, for $j\in J(i)$ and all $i$, respectively. Hence, we have
$$\sum_{j\in J(i)}\underline{\zeta}^i_jx_{d,j}\leq \tilde{u}_i\leq \sum_{j\in J(i)}\overline{\zeta}^i_jx_{d,j}, \ \ \textrm{and} \ \ \ \sum_{j\in J(i)}\underline{\xi}^i_jx_{d,j}\leq \tilde{u}_i\leq \sum_{j\in J(i)}\overline{\xi}^i_jx_{d,j}, \ \ \ \forall i \in I. $$

\begin{rem}
	One traditional approach  is to construct an interval that just contains both estimation intervals to build a safe one, e.g., 
	\begin{equation}
	\label{eq_polytopeDDU_FL}
		\min\{\sum_{j\in J(i)}\underline\zeta^i_jx_{d,j}, \sum_{j\in J(i)}\underline\xi^i_jx_{d,j}\}\leq \tilde u_i\leq \max\{\sum_{j\in J(i)}\overline \zeta^i_jx_{d,j}, \sum_{j\in J(i)}\overline\xi_j^ix_{d,j}\}, \ \forall i \in I,
	\end{equation}
	which nevertheless could be an overestimation. As another strategy, we believe that one of those two estimations should be correct, but we just do not know which is that one. Moreover, we probably do not want to rely on any particular consulting service for all clients, which can be achieved by limiting the dependence on either one.  
\end{rem}

Following the aforementioned remark, we construct next a mixed integer DDU set.  
\begin{subequations}
	\begin{align}
		\mathcal U^I(\mathbf x)=\Big\{(\mathbf u, \boldsymbol{\delta}):\ &u_i=\underline u_i+\tilde u_i \quad\forall i\in I \label{eq_RFL_induced}\\
		&\tilde u_i\geq \delta_i^1\sum_{j\in J(i)}\underline{\zeta}^i_jx_{d,j}+\delta_i^2\sum_{j\in J(i)}\underline{\xi}^i_jx_{d,j} \quad\forall i\in I \\
		& \tilde u_i\leq\delta_i^1\sum_{j\in  J(i)}\overline{\zeta}^i_jx_{d,j}+\delta_i^2\sum_{j\in  J(i)}\overline{\xi}^i_jx_{d,j} \quad\forall i\in I \\
		 &\delta_i^1+\delta_i^2=1 \quad\forall i\in I, \ \  
	     \sum_{i\in I}\delta_i^h\leq k_h \ \ h=1,2 \label{eq_RFL_card}\\
		&\sum_{i\in I}\tilde u_i\leq \alpha \sum_{j\in J} x_{c,j}, \ \  
		\delta_i^1,\delta_i^2\in\{0,1\}\quad\forall i\in I\Big\} \label{eq_RFL_induced_upper}
	\end{align}
\end{subequations} 
Using binary variables $\delta_i^1$ and $\delta_i^2$,  constraints in \eqref{eq_RFL_card}  indicate that estimation from only one consulting service will be adopted for every client, while reliance on a particular consulting service across all $i$'s is restricted by some limit. The constraint in \eqref{eq_RFL_induced_upper} bounds the total induced demand by a multiple of the total installed capacities. To demonstrate the impact of mixed integer DDU,  we also consider the following polytope DDU for comparison. 
\begin{eqnarray}
\label{eq_FL_DDU_polytope}
\mathcal U^C(\mathbf x) =\{\mathbf u: \eqref{eq_RFL_induced},\eqref{eq_polytopeDDU_FL}, \eqref{eq_RFL_induced_upper}\}	
\end{eqnarray}

We study a problem with 40 sites, with data on their locations, distances and nominal demands adopted from \cite{snyder2005reliability}. The set $J(i)$ contains all sites within 6 distance units to $i$. Construction cost $f_j$ is set to the product between total demand of sites in $J(i)$ and a random number between 0.4 and 2.4, unit capacity cost $a_j$ is a random number between 0.3 and 0.5; capacity bounds $\underline A_j$ and $\overline A_j$ are set as 10\% and 100\% of the total demand in set $J(j)$, respectively. In the DDU sets, $\underline{\zeta}^i_j$ and $\overline{\zeta}^i_j$ are set to  1\% and 3\% of the total demand in set $J(j)$, $\overline \xi_j$ and $\underline \xi_j$ are $(1+r)$ of their $\zeta_j$ counterparts with $r$ being a changing parameter, and $\alpha$ is chosen to be 0.5. Also, coefficient $\rho$ is set to 1. 

Instances of \textbf{RFL-L} model with uncertainty sets $\mathcal U^C(\mathbf x)$ and $\mathcal U^I(\mathbf x)$ under different $r$ and $k_2$ values (with $k_1$ fixed to 40) are computed. Results are shown in Table \ref{TB_FL_MIPDDU}, including lower and upper bounds upon termination, relative gaps between those bounds, number of iterations, and computational time in seconds.

\begin{table}[htbp]
	\centering
	\caption{Computational Results of \textbf{RFL-L} with $\mathcal U^C(\mathbf x)$ or $\mathcal U^I(\mathbf x)$}
	\resizebox{\textwidth}{!}{
	\begin{tabular}{|c|ccccc|cccccr|}
    \hline
& \multicolumn{5}{c|}{$\mathcal U^C(\mathbf x)$ {\textbf/} Standard Parametric C\&CG}               & \multicolumn{6}{c|}{$\mathcal U^I(\mathbf x)$ {\textbf/} Parametric C\&CG-MIU} \\
\cline{2-12}    $r$     & LB    & UB    & Gap   & Iter  & Time(s) & $k_2$     & LB    & UB    & Gap   & Iter  & \multicolumn{1}{c|}{Time(s)} \\
\hline
          &       &       &       &       &       & 1     & 22859.55 & 22859.80 & 0.00\% & 3     & 118.88 \\
&       &       &       &       &       & 2     & 22861.66 & 22863.50 & 0.01\% & 3     & 80.38 \\
0.1   & 23410.39 & 23410.39 & 0.00\% & 3     & 114.30 & 3     & 22865.35 & 22867.15 & 0.01\% & 3     & 81.58 \\
&       &       &       &       &       & 4     & 22866.63 & 22870.41 & 0.02\% & 3     & 130.12 \\
&       &       &       &       &       & 5     & 22869.58 & 22873.55 & 0.02\% & 3     & 76.59 \\
\hline
&       &       &       &       &       & 1     & 22863.48 & 22863.98 & 0.00\% & 3     & 108.31 \\
&       &       &       &       &       & 2     & 22867.70 & 22871.37 & 0.02\% & 3     & 129.03 \\
0.2   & 24067.47 & 24067.47 & 0.00\% & 3     & 90.80 & 3     & 22875.07 & 22878.67 & 0.02\% & 3     & 150.62 \\
&       &       &       &       &       & 4     & 22876.81 & 22885.19 & 0.04\% & 3     & 75.25 \\
&       &       &       &       &       & 5     & 23204.73 & 23205.80 & 0.00\% & 5     & 261.47 \\
\hline
&       &       &       &       &       & 1     & 22867.41 & 22868.15 & 0.00\% & 3     & 107.07 \\
&       &       &       &       &       & 2     & 22873.74 & 22879.24 & 0.02\% & 3     & 78.50 \\
0.3   & 24169.91 & 24169.92 & 0.00\% & 2     & 99.85 & 3     & 23196.41 & 23204.24 & 0.03\% & 5     & 192.13 \\
&       &       &       &       &       & 4     & 23259.56 & 23264.87 & 0.02\% & 5     & 258.62 \\
&       &       &       &       &       & 5     & 23414.51 & 23417.67 & 0.01\% & 5     & 114.75 \\
\hline
&       &       &       &       &       & 1     & 22871.34 & 22872.33 & 0.00\% & 3     & 110.65 \\
&       &       &       &       &       & 2     & 23195.18 & 23201.48 & 0.03\% & 4     & 239.47 \\
0.4   & 24640.30 & 24640.31 & 0.00\% & 2     & 167.50 & 3     & 23406.38 & 23408.18 & 0.01\% & 5     & 141.45 \\
&       &       &       &       &       & 4     & 23433.39 & 23447.23 & 0.06\% & 4     & 171.22 \\
&       &       &       &       &       & 5     & 23452.71 & 23460.94 & 0.04\% & 4     & 91.64 \\
\hline
&       &       &       &       &       & 1     & 22875.27 & 22876.50 & 0.01\% & 3     & 114.73 \\
&       &       &       &       &       & 2     & 23210.01 & 23211.20 & 0.01\% & 3     & 154.43 \\
0.5   & 25636.39 & 25636.4 & 0.00\% & 2     & 281.45 & 3     & 23437.45 & 23444.53 & 0.03\% & 4     & 153.94 \\
&       &       &       &       &       & 4     & 23442.98 & 23462.16 & 0.08\% & 3     & 134.08 \\
&       &       &       &       &       & 5     & 23624.46 & 23649.00 & 0.10\% & 3     & 155.10 \\
\hline
	\end{tabular}}%
	\label{TB_FL_MIPDDU}%
\end{table}%

From Table \ref{TB_FL_MIPDDU}, it is obvious that the cost associated with $\mathcal U^C (\bf x)$ is always higher than that with $\mathcal U^I (\bf x)$. It is sensible as $\mathcal U^C (\bf x)$ is defined on top of \eqref{eq_polytopeDDU_FL}, which is a larger decision-dependent  interval containing both estimation intervals and should yield a more conservative solution. A significant difference  can often be seen between the costs with $\mathcal U^C (\bf x)$ and $\mathcal U^I (\bf x)$, respectively. Note that when two consulting services provide quite different estimations, e.g., $r=0.5$, the cost with $\mathcal U^C(\bf x)$ could be 12\% more than that with $\mathcal U^I(\bf x)$, indicating that \textbf{RFL-L} with $\mathcal U^I(\bf x)$ derives a less conserve solution. Finally, regarding the computational expenses for those two types of DDU sets,  mixed integer DDU does not necessarily demand for a much longer computational time. Nevertheless, due to the more sophisticated structure underlying $\mathcal U^I (\bf x)$, parametric C\&CG-MIU typically needs more iterations to converge.   

\begin{figure}[htp]
\textrm{ \ \ \ \ \ \ }	\begin{subfigure}[t]{0.45\textwidth}
		\pgfplotstableread{Sp7.dat}{\Sp}
		\pgfplotsset{tick label style={font=\small\bfseries},
			label style={font=\small},
			legend style={font=\tiny}
		}
		\begin{tikzpicture}[scale=1]
			\begin{axis}[ymajorgrids=true,legend style={at={(1.1,1.4)},anchor=north,legend columns=-1,/tikz/every even column/.append style={column sep=5pt}},height=5cm,width=7.5cm,
				xmin=0,xmax=6,ymin=20000,ymax=27500,
				xtick={0,1,2,3,4,5,6},
				ytick={20000,21000,22000,23000,24000,25000,27000},
				yticklabels={$2.0$,$2.1$,$2.2$,$2.3$,$2.4$,$2.5$,$+\infty$},
				extra y ticks={26400},
				extra y tick label={$\vdots\ $},
				xlabel= Iterations]
				\addplot [blue,semitransparent,very thick,dashed] table [x={k5it}, y={k5lb}] {\Sp};
				\addplot [red,semitransparent,very thick,dashed] table [x={cit}, y={clb}] {\Sp};
				\addplot [blue,semitransparent,very thick] table [x={k5it}, y={k5ub}] {\Sp};
				\addlegendentry{$\mathcal U^I(\mathbf x),\ r=0.4, \ k_2=5$,}
				\addplot [red,semitransparent,very thick] table [x={cit}, y={cub}] {\Sp};
				\addlegendentry{$\mathcal U^C(\mathbf x)$};
				\node at (300,406.7) [circle, scale=0.4, draw=red!50,fill=red!50] {};
				\node at (500,320.5) [circle, scale=0.4, draw=blue!50,fill=blue!50] {};
			\end{axis}
		\end{tikzpicture}
		\caption{Convergence over Iterations}
		\label{ConIter}
	\end{subfigure}
	\begin{subfigure}[t]{0.45\textwidth}
		\pgfplotstableread{Sp7.dat}{\Sp}
		\pgfplotsset{tick label style={font=\small\bfseries},
			label style={font=\small},
			legend style={font=\tiny}
		}
		\begin{tikzpicture}[scale=1]
			\begin{axis}[ymajorgrids=true,height=5cm,width=7.5cm,
				xmin=0,xmax=300,ymin=20000,ymax=27500,
				xtick={0,50,100,150,200,250,300},
				ytick={20000,21000,22000,23000,24000,25000,27000},
				yticklabels={$2.0$,$2.1$,$2.2$,$2.3$,$2.4$,$2.5$,$+\infty$},
				extra y ticks={26400},
				extra y tick label={$\vdots\ $},
				xlabel= Time (s)]
				\addplot [blue,semitransparent,very thick,dashed] table [x={k5ti}, y={k5lb}] {\Sp};
				\addplot [red,semitransparent,very thick,dashed] table [x={cti}, y={clb}] {\Sp};
				\addplot [blue,semitransparent,very thick] table [x={k5ti}, y={k5ub}] {\Sp};
				\addplot [red,semitransparent,very thick] table [x={cti}, y={cub}] {\Sp};
				\node at (90.8,406.7) [circle, scale=0.4, draw=red!50,fill=red!50] {};
				\node at (261.5,320.5) [circle, scale=0.4, draw=blue!50,fill=blue!50] {};
			\end{axis}
		\end{tikzpicture}
		\caption{Convergence over Time}
		\label{ConTime}
	\end{subfigure}
	\caption{Convergence of \textbf{RFL-L} for $r=0.4$ and $k_2=5$}
	\label{ConFL}
\end{figure}
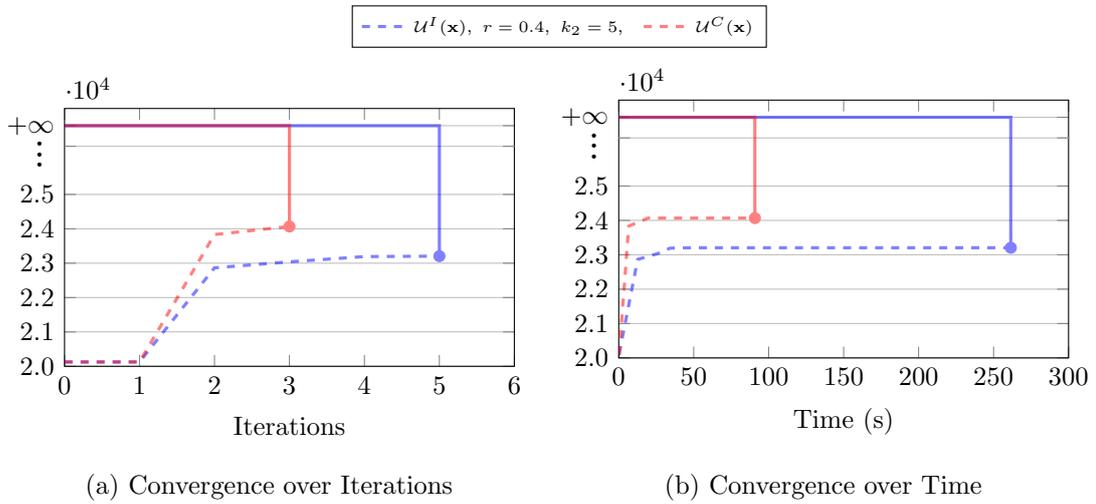

In Figure \ref{ConFL} we show the convergence behaviors of algorithms over iterations and time for the case $r=0.2$ with $k_2=5$ for  $\mathcal U^C(\mathbf x)$ and $\mathcal U^I(\mathbf x)$, with upper and lower bounds represented by solid and dashed lines, respectively. Their starting lower bounds are $w_R$ calculated according to equation (\ref{eq_MIP_RO}). It can be observed that algorithms for $\mathcal U^C(\mathbf x)$ and $\mathcal U^I(\mathbf x)$  are efficient as both of them have a quick converge behavior. Between them, parametric C\&CG-MIU for $\mathcal U^I(\mathbf x)$ requires more iterations to converge. Especially, more subproblems are needed to ensure a feasible (actually optimal) solution.

\subsection{Modeling Decisions of Temporary Facility  by MIP Recourse}
\label{ssec_MIPRecourse}
Assuming that temporary facilities of fixed capacities can be built in the second stage to serve clients, 
we modify $\mathbf{RFL-L}$ in Section \ref{FL_MIPDDU} to $\mathbf{RFL-I}$ as in the following
\begin{subequations}
	\begin{align*}
	\mathbf{RFL-I} :	&\min_{(\mathbf x_c,\mathbf x_d)\in\mathcal X} \sum_{j\in J} (f_jx_{d,j}+a_jx_{c,j}) + \max_{\mathbf u\in\mathcal U(\mathbf x)} \ \min_{(\mathbf y, \mathbf z)\in\mathcal{Y}^I(\mathbf x,\mathbf u)} \ \sum_{i\in I}\sum_{j\in J}c_{ij}y_{ij}+\sum_{j\in J}f^+_jz_j
	\end{align*}
    where   binary variable $z_j$ denotes the construction decision of a temporary facility at site $j$, and $f^+_j$ is its fixed cost. Hence, the feasible set of the recourse problem is changed to 
    \begin{align*}
		\mathcal Y^I(\mathbf x,\mathbf u)=\Big\{(\mathbf y, \mathbf z)\in \mathbb{R}^{|I|\times|J|}_+\times\mathbb B^{|J|}:\ \sum_{j\in J}y_{ij} \geq u_i \ \forall i\in I,\ \sum_{i\in I}y_{ij}\leq x_{c,j}+h_jz_j \ \forall j\in J\Big\}
	\end{align*}
\end{subequations}
with $h_j$ denoting the capacity associated with the temporary facility at $j$.

Modifying the set of instances adopted in Section \ref{FL_MIPDDU}, we set  $h_j$ to some number between 30 to 80 depending on the nominal demand of site $j$. Specifically, the higher $\underline u_j$ is, the larger $h_j$ we have. Let $\bar a$ denote the largest  unit capacity cost among all $a_j$'s. Then, we set fixed cost $f^+_j=5h_j\bar a$. The rest parameters are the same as those in Section \ref{FL_MIPDDU}. Instances of $\mathbf{RFL-I}$ for $\mathcal U(\bf x)=\mathcal{U}^C(\bf x)$ and $\mathcal U(\bf x)=\mathcal{U}^I(\bf x)$ are computed by nested and extended nested parametric C\&CG, respectively, with the naive initialization strategy. 
All results are shown in Table \ref{TB_FL_MIPRecourse}, where column  ``O-iter'' reports the total number of iterations for outer C\&CG procedure and ``I-iter''  the total number of iterations for inner C\&CG ones.

\begin{table}[htbp]
	\centering
	\caption{Computational Results of \textbf{RFL-I} with Temporary Facilities}
	\resizebox{\textwidth}{!}{
		\begin{tabular}{|c|cccccc|ccccccc|}
			\hline
			& \multicolumn{6}{c|}{\textbf{RFL-I} with $\mathcal U^C(\bf x)$}             & \multicolumn{7}{c|}{\textbf{RFL-I} with $\mathcal U^I(\bf x)$} \\
			\cline{2-14}    $r$     & LB    & UB    & Gap   & O-iter & I-iter & Time(s) & $k_2$     & LB    & UB    & Gap   & O-iter & I-iter & Time(s) \\
			\hline
          &       &       &       &       &       &       & 1     & 19453.79 & 19543.11 & 0.46\% & 2     & 5     & 69.51 \\
&       &       &       &       &       &       & 2     & 19449.94 & 19544.91 & 0.49\% & 2     & 5     & 59.97 \\
0.1   & 19510.41 & 19510.41 & 0.00\% & 2     & 10    & 1771.99 & 3     & 19460.59 & 19546.57 & 0.44\% & 2     & 5     & 58.82 \\
&       &       &       &       &       &       & 4     & 19462.99 & 19548.24 & 0.44\% & 2     & 5     & 64.84 \\
&       &       &       &       &       &       & 5     & 19466.59 & 19560.05 & 0.48\% & 2     & 5     & 57.91 \\
\hline
&       &       &       &       &       &       & 1     & 19458.07 & 19545.10 & 0.45\% & 2     & 5     & 61.89 \\
&       &       &       &       &       &       & 2     & 19465.34 & 19558.89 & 0.48\% & 2     & 5     & 59.41 \\
0.2   & 19563.92 & 19597.24 & 0.17\% & 2     & 10    & 121.30 & 3     & 19471.52 & 19484.03 & 0.06\% & 2     & 10    & 462.75 \\
&       &       &       &       &       &       & 4     & 19476.32 & 19490.25 & 0.07\% & 2     & 10    & 314.58 \\
&       &       &       &       &       &       & 5     & 19476.11 & 19498.26 & 0.11\% & 2     & 10    & 362.53 \\
\hline
&       &       &       &       &       &       & 1     & 19462.36 & 19547.10 & 0.43\% & 2     & 5     & 61.56 \\
&       &       &       &       &       &       & 2     & 19473.04 & 19485.07 & 0.06\% & 2     & 10    & 879.67 \\
0.3   & 19616.90 & 19664.66 & 0.24\% & 2     & 10    & 134.44 & 3     & 19471.35 & 19502.96 & 0.16\% & 2     & 10    & 1099.01 \\
&       &       &       &       &       &       & 4     & 19489.97 & 19506.85 & 0.09\% & 2     & 10    & 267.01 \\
&       &       &       &       &       &       & 5     & 19507.21 & 19520.53 & 0.07\% & 2     & 10    & 465.45 \\
\hline
&       &       &       &       &       &       & 1     & 19466.60 & 19559.27 & 0.47\% & 2     & 5     & 58.38 \\
&       &       &       &       &       &       & 2     & 19478.80 & 19493.07 & 0.07\% & 2     & 10    & 335.22 \\
0.4   & 19706.25 & 19760.07 & 0.27\% & 2     & 10    & 116.48 & 3     & 19488.39 & 19507.34 & 0.10\% & 2     & 10    & 341.16 \\
&       &       &       &       &       &       & 4     & 19510.48 & 19526.26 & 0.08\% & 2     & 10    & 290.12 \\
&       &       &       &       &       &       & 5     & 19510.21 & 19545.91 & 0.18\% & 2     & 10    & 483.04 \\
\hline
&       &       &       &       &       &       & 1     & 19463.22 & 19480.34 & 0.09\% & 2     & 10    & 454.08 \\
&       &       &       &       &       &       & 2     & 19478.74 & 19507.89 & 0.15\% & 2     & 10    & 506.51 \\
0.5   & 19837.78 & 19841.37 & 0.02\% & 3     & 10    & 115.86 & 3     & 19505.48 & 19535.94 & 0.16\% & 2     & 10    & 273.31 \\
&       &       &       &       &       &       & 4     & 19510.10 & 19548.27 & 0.20\% & 2     & 10    & 236.64 \\
&       &       &       &       &       &       & 5     & 19550.04 & 19562.41 & 0.06\% & 2     & 11    & 307.42 \\
\hline
	\end{tabular}}%
	\label{TB_FL_MIPRecourse}%
\end{table}%

From Table \ref{TB_FL_MIPRecourse}, it is interesting to note that the benefit of adopting mixed integer DDU becomes much less significant, although the cost with $\mathcal U^C(\bf x)$ is consistently  larger than that with $\mathcal U^I(\bf x)$ across all instances. This can be explained by the introduction of temporary facilities in the recourse stage. Those temporary facilities lessen the dependence on the permanent ones that are determined before demand realization, and help to mitigate the demand variability on the total cost. Hence,we believe that a strong and cost-effective recourse capacity is highly desired in a complex uncertain environment, especially when the randomness is hard to have a sound understanding.   Moreover, almost all instances can be solved in two to three outer iterations.  In general, the computational time for $\mathbf{RFL-I}$ with $\mathcal U^I(\bf x)$ could be much longer than that for $\mathbf{RFL-I}$ with $\mathcal U^C(\bf x)$ and than that for $\mathbf{RFL-L}$ with $\mathcal U^I(\bf x)$. Additionally, we note that numerical gaps are often noticeable among instances of $\mathbf{RFL-I}$ with $\mathcal U^I(\bf x)$.  Those observations indicate the computational challenge arising from interacting mixed integer structures in recourse and DDU sets.

\subsection{Approximations for DDU and MIP Recourse}
We also compute all instances of $\mathbf{RFL-I}$ using the approximation variant of parametric C\&CG-MIU presented in Section \ref{sect_MIDDUMIP}. Regarding the termination condition,  we stop the algorithm when the master problem outputs a previously generated first-stage decision. All computational results are reported in Table \ref{TB_FL_MIPApprox}.

We observe that non-zero gaps between lower and upper bounds exist when the algorithm terminates, which reflects the algorithm's approximation nature. Comparing  \textbf{RFL-I} with $\mathcal U^C(\bf x)$ and with $\mathcal U^I(\bf x)$, it is interesting to see that the relative gaps for the latter case are significantly smaller than those for the former one. Basically, the majority of relative gaps with $\mathcal U^I(\bf x)$ are less than 10\%, while such gaps with $\mathcal U^C(\bf x)$ are generally between 10\% and 25\%. It shows that this approximation variant works very well with mixed integer DDU sets. Indeed,  comparing Tables \ref{TB_FL_MIPRecourse}  and \ref{TB_FL_MIPApprox}, it is worth highlighting that this approximation variant generally produces high quality solutions in a fast fashion, and it can scale well to handle challenging instances.

\begin{table}[htbp]
	\centering
	\caption{Approximation Results for \textbf{RFL-I} with Temporary Facilities}
	\resizebox{\textwidth}{!}{
		\begin{tabular}{|c|ccccc|cccccc|}
			\hline
          & \multicolumn{5}{c|}{\textbf{RFL-I} with $\mathcal U^C(\bf x)$}     & \multicolumn{6}{c|}{\textbf{RFL-I} with $\mathcal U^I(\bf x)$} \\
\cline{2-12}    $r$     & LB    & UB    & Gap   & Iter & Time(s) & $k_2$     & LB    & UB    & Gap   & Iter  & Time(s) \\
\hline
          &       &       &       &       &       & 1     & 19416.93 & 20189.23 & 3.83\% & 3     & 215.86 \\
&       &       &       &       &       & 2     & 19416.93 & 20278.05 & 4.25\% & 3     & 200.97 \\
0.1   & 19452.8 & 21266.37 & 8.53\% & 3     & 229.52 & 3     & 19416.93 & 20338.76 & 4.53\% & 3     & 246.52 \\
&       &       &       &       &       & 4     & 19416.93 & 20399.11 & 4.81\% & 3     & 318.36 \\
&       &       &       &       &       & 5     & 19416.93 & 20459.45 & 5.10\% & 3     & 180.08 \\
\hline
&       &       &       &       &       & 1     & 19429.06 & 20083.26 & 3.26\% & 3     & 240.58 \\
&       &       &       &       &       & 2     & 19429.06 & 20260.91 & 4.11\% & 3     & 292.58 \\
0.2   & 19506.10 & 22765.87 & 14.32\% & 3     & 73.06 & 3     & 19429.06 & 20382.26 & 4.68\% & 3     & 280.08 \\
&       &       &       &       &       & 4     & 19429.06 & 20502.88 & 5.24\% & 3     & 178.69 \\
&       &       &       &       &       & 5     & 19429.06 & 20623.38 & 5.79\% & 3     & 222.33 \\
\hline
&       &       &       &       &       & 1     & 19440.38 & 20722.05 & 6.19\% & 3     & 255.38 \\
&       &       &       &       &       & 2     & 19440.38 & 20988.53 & 7.38\% & 3     & 213.98 \\
0.3   & 19560.09 & 23682.97 & 17.41\% & 3     & 73.41 & 3     & 19440.38 & 21170.91 & 8.17\% & 3     & 160.41 \\
&       &       &       &       &       & 4     & 19440.38 & 21352.17 & 8.95\% & 3     & 280.86 \\
&       &       &       &       &       & 5     & 19440.38 & 21533.85 & 9.72\% & 3     & 213.83 \\
\hline
&       &       &       &       &       & 1     & 19452.07 & 20615.67 & 5.64\% & 3     & 224.38 \\
&       &       &       &       &       & 2     & 19452.07 & 20970.97 & 7.24\% & 3     & 279.91 \\
0.4   & 19650.87 & 25277.39 & 22.26\% & 3     & 78.88 & 3     & 19452.07 & 21214.15 & 8.31\% & 3     & 124.14 \\
&       &       &       &       &       & 4     & 19452.07 & 21456.11 & 9.34\% & 3     & 217.47 \\
&       &       &       &       &       & 5     & 19452.07 & 21698.34 & 10.35\% & 3     & 160.34 \\
\hline
&       &       &       &       &       & 1     & 19463.75 & 21035.65 & 7.47\% & 3     & 255.36 \\
&       &       &       &       &       & 2     & 19463.75 & 21480.46 & 9.39\% & 3     & 135.19 \\
0.5   & 19717.20 & 26294.94 & 25.02\% & 3     & 108.36 & 3     & 19463.75 & 21785.53 & 10.66\% & 3     & 157.75 \\
&       &       &       &       &       & 4     & 19463.75 & 22088.36 & 11.88\% & 3     & 106.12 \\
&       &       &       &       &       & 5     & 19463.75 & 22391.15 & 13.07\% & 3     & 102.37 \\
\hline
	\end{tabular}}%
		\label{TB_FL_MIPApprox}%
\end{table}%

\section{Conclusions}
\label{sect_conclusion}
In this paper, we present a systematic study on DDU-based two-stage RO with mixed integer structures, which appear in either or both DDU set and recourse problem. Exact solution methods based on parametric C\&CG are designed and analyzed to produce accurate solutions for those challenging problems, which also address some limitations presented in their existing counterparts. Also, a computationally friendly approximation variant is developed for the most complex two-stage RO that has both mixed integer DDU and recourse. Finally, we conduct a numerical study on instances of robust facility location problem with DDU demands to appreciate the modeling strength of mixed integer structures and the solution strength of the developed algorithms. 

As a future research direction, it would be interesting to develop modeling techniques to capture complex decision dependence by a mixed integer DDU that leads to desirable computational behavior. Also, enhancing the developed algorithms to deal with large-scale instances are expected. Finally, employing the developed algorithms to solve challenging decision making problems should be carried out to support various real systems.

\bibliographystyle{plain}
\bibliography{references_MIP,robustUC}

\newpage
\begin{center}
	\textbf{\Large{Appendix}}
\end{center}
\setcounter{section}{0}

\setcounter{equation}{0}
\setcounter{thm}{0}
\setcounter{rem}{0}
\def\theequation{A.\arabic{equation}}
\def\thethm{A.\arabic{thm}}
\def\therem{A.\arabic{rem}}

Consider the general bilevel linear optimization problem as in the following 
\begin{subequations}
\label{eq_bilevel_general}
\begin{align}
\mathbf{BLO}: \	v^*= \max \quad& \bf c_{\bf x}\bf x+\bf c_{\bf y} \bf y^* \\
	     \mbox{s.t.} \quad &\bf x\in \cal{X}, \ \ \bf y^*\in \arg\min \Big\{\mathbf f\mathbf y:  \mathbf B_{\mathbf y}\mathbf y\geq \mathbf d-\mathbf{B}_{\mathbf x}\mathbf x, \ \mathbf{y}\geq \mathbf{0}\Big\}, \label{ep_blp_low}	
\end{align}
\end{subequations}
where $\mathcal{X}$ is a non-empty mixed integer set, $\bf x$ is referred to as the upper-level decision, the optimization problem in \eqref{ep_blp_low} the lower-level decision making problem, and $\bf y^*$ an associated optimal solution.

Suppose that $\mathbf{BLO}$ has a finite optimal value, which is the case for all bilevel optimization formulations involved in all algorithmic procedures, except $\mathbf{IMP_o}$ of \eqref{eq_MP-io} for the inner C\&CG subroutine for optimality. Then, we can replace the lower-level problem by its optimality conditions. For example, if KKT conditions are employed, $\mathbf{BLO}$ can be converted into the following single-level formulation
\begin{subequations}
\label{eq_bilevel_general_KKT}
\begin{align}
		v^*=  \max\quad&\mathbf c_{\mathbf y}\mathbf y\\
		\mathrm{s.t.}\quad&\mathbf x\in \mathcal X, \ \ \ \mathbf B_{\mathbf y}\mathbf y\geq \mathbf d-\mathbf{B}_{\mathbf x}\mathbf x, \ \ \ \mathbf B_{\mathbf y}^\intercal\boldsymbol\lambda\leq \mathbf	c^\intercal_{\mathbf y}\\
		& \boldsymbol\lambda\circ(\mathbf B_{\mathbf y}\mathbf y-\mathbf d+\mathbf{B}_{\mathbf x}\mathbf x) = \mathbf 0, \ \ \ \mathbf y\circ (\mathbf	c^\intercal_{\mathbf y}- \mathbf B_{\mathbf y}^\intercal\boldsymbol\lambda)=\mathbf 0 \label{Aeq_bilevel_complementary}\\
		& \mathbf{y}\geq \mathbf{0}, \ \boldsymbol\lambda\geq \mathbf 0
\end{align}
\end{subequations}
where $\boldsymbol\lambda$ denotes the dual variables of the lower-level problem.  Note that constraints in \eqref{Aeq_bilevel_complementary} are referred to as the complementarity ones. Those constraints actually can be linearized, converting \eqref{eq_bilevel_general_KKT} into an MIP. Specifically, consider the $j$-th constraint in the  first set of complementarity constraints in \eqref{Aeq_bilevel_complementary}, i.e.,
$$\lambda_j (\mathbf d-\mathbf{B}_{\mathbf x}\mathbf x-\mathbf{B}_{\mathbf y}\mathbf y)_j=0.$$
Let $\delta_j$ be a binary variable and recall that $M$ denotes a sufficiently large number. Then, this complementarity constraint can be replaced by the next two linear constraints.
\begin{equation*}
	\lambda_j\leq M\delta_j,  \  (\mathbf{B}_{\mathbf y}\mathbf y - \mathbf d+\mathbf{B}_{\mathbf x}\mathbf x)_j\leq M(1-\delta_j)
\end{equation*}
Note that when $\delta_j=0$, together with the nonnegativity of $\lambda_j$, we have
$\lambda_j=0$. When $\delta_j=1$, together with the primal constraint $(\mathbf B_{\mathbf y}\mathbf y)_j\geq (\mathbf d-\mathbf{B}_{\mathbf x}\mathbf x)_j$, we have $(\mathbf B_{\mathbf y}\mathbf y - \mathbf d-\mathbf{B}_{\mathbf x}\mathbf x)_j=0$. Hence, such binary variables and linear constraints help us achieve the same effect as those complementarity constraints.
Recall that the only nonlinear constraints in $\mathcal{OU}$ and $\mathcal{OV}$ are complementarity ones. Hence, this technique can be used to convert them into MIP sets. 

Similarly, duality-based reformulation method can also be used to convert $\mathbf{BLO}$ into a single-level bilinear  problem.  
Currently, it is a formulation computable (with a rather restricted capacity) by some professional solvers. Indeed, this strategy works very well when every bilinear term is a product between a binary variable and a continuous variable, since it simply can be linearized. For example, as stated in Remark \ref{rem_inner_optimality},  $\mathbf{IMP_o}$ of \eqref{eq_MP-io} is converted into such a formulation, where not every lower-level problem is required to be feasible.

\end{document}